\documentclass[12pt,twoside,a4paper]{amsart}
\usepackage{amssymb}
\usepackage{latexsym}
\usepackage{amsmath}
\usepackage{euscript}
\newcommand{\RR}{\mathbb{R}}
\newcommand{\CC}{\mathbb{C}}

\newcommand{\DD}{\mathbb{D}}

\newcommand{\ZZ}{\mathbb{Z}}
\newcommand{\TT}{\mathbb{T}}
\font\gothic=eufm10

\def\gD{\mathfrak{D}}

\def\gR{\mbox{\gothic\char'122}}

\def\gb{\mbox{\gothic\char'142}}

\def\filbcirc#1{\hbox{\rlap{\hbox to 25pt {$\bigcirc$\hfill}}
\hbox to 2pt{\hfill \small #1\hfill}}}

\def\sD{{\mathfrak D}}      \def\sF{{\mathfrak F}}
   \def\sH{{\mathfrak H}}

\def\gh{{\mathfrak h}}
\def\gs{{\mathfrak s}}

      \def\dC{{\mathbb C}}
\def\dD{{\mathbb D}}

   \def\dN{{\mathbb N}}   
      
   \def\dT{{\mathbb T}}

\def\cD{{\mathcal D}}      \def\cF{{\mathcal F}}
   \def\cH{{\mathcal H}}   
   \def\cK{{\mathcal K}}   \def\cL{{\mathcal L}}
\def\cM{{\mathcal M}}   \def\cN{{\mathcal N}}   
      \def\cR{{\mathcal R}}
\def\cS{{\mathcal S}}   \def\cT{{\mathcal T}}   \def\cU{{\mathcal U}}

\newcommand{\IPi}{\hbox{$\,$I\vbox{\moveleft 2pt\hbox{$\Pi$}}}}

\newcommand{{\bx}}{\bf x}

\newcommand{{{\bff}}}{\bf f}

\newcommand{{\bw}}{\bf w}
\newcommand{{\bg}}{\bf g}
\newcommand{{\bb}}{\bf b}
\newcommand{{\ba}}{\bf a}
\newcommand{{\bq}}{\bf q}
\newcommand{{\bC}}{\bf C}

\newcommand{\lam}{\lambda}
\newcommand{\om}{\omega}
\newcommand{\ptp}{p\times p}
\newcommand{\ptq}{p\times q}
\newcommand{\qtq}{q\times q}

\newcommand{\qtp}{q\times p}
\newcommand{\mtm}{m\times m}
\newcommand{\ntn}{n\times n}

\newcommand{\wtilde}{\widetilde}

\def\h#1{{{\hat #1} }}
\def\wt#1{{{\widetilde #1} }}
\def\wh#1{{{\widehat #1} }}

\def\bm\chi{\mbox{\boldmath$\chi$}}

\def\ran{{\rm rng\,}}

\def\dim{{\rm dim\,}}

\def\rank{{\rm rank\,}}

\let\xker=\ker \def\ker{{\xker\,}}

\newtheorem{thm}{Theorem}[section]
\newtheorem{proposition}[thm]{Proposition}
\newtheorem{corollary}[thm]{Corollary}
\newtheorem{lem}[thm]{Lemma}
\newtheorem{definition}[thm]{Definition}

\theoremstyle{definition}
\newtheorem{example}{Example}
\newtheorem{remark}[thm]{Remark}

\numberwithin{equation}{section}

\begin{document}

\title
{Bitangential interpolation in generalized Schur classes}
\author{Vladimir Derkach}
\address{Department of Mathematics \\
Donetsk State University \\
Universitetskaya str. 24 \\
83055 Donetsk \\
Ukraine} \email{derkach.v@gmail.com}
\author{Harry Dym}
\address{Department of Mathematics \\
The Weizmann Institute of Science\\
Rehovot 76100, Israel}
\email{harry.dym@weizmann.ac.il} \dedicatory{}
\date{\today}
\thanks{V. Derkach  wishes to thank the Weizmann Institute of Science
for hospitality and support}
\subjclass{Primary 30E05; Secondary 46C20, 47A48, 47B50.}
\keywords{ Bitangential interpolation, generalized Schur class,
Kre\u{\i}n-Langer factorization, resolvent matrix, linear fractional transformation, coprime
factorization}

\begin{abstract}
Bitangential interpolation problems in the class of matrix valued functions
in the generalized Schur class are considered in both the open unit disc and
the open right half plane, including problems in which the solutions is not
assumed to be holomorphic at the interpolation points. Linear
fractional representations of the set of solutions to these problems are
presented for invertible and singular Hermitian Pick matrices. These
representations make use of a description of the ranges of linear fractional
transformations with suitably chosen domains that was developed in
\cite{DerDym08}.
\end{abstract}

\maketitle



\section{Introduction }
The main objective of this paper is to study bitangential
interpolation problems in the {\it generalized Schur class}
$\cS_{\kappa}^{p\times q}(\Omega_+)$ of $\ptq$ matrix valued
functions that are meromorphic in $\Omega_+$ and for which the
kernel
\begin{equation}\label{kerLambda}
{\mathsf \Lambda}_\omega^s(\lambda)=
\frac{I_{p}-s(\lambda)s(\omega)^*}{\rho_\omega(\lambda)}
\end{equation}
has $\kappa$ negative squares in ${\mathfrak h}_s^+\times{\mathfrak
h}_s^+$ (see~\cite{KL}), where ${\mathfrak h}_s^+$ denotes the
domain of holomorphy of $s$ in $\Omega_+$, $\Omega_+$ is either
equal to $\dD=\{\lam\in\dC:\,|\lam|<1\}$ or
$\Pi_+=\{\lam\in\dC:\,\lam+\overline{\lam}>0\}$, and
\[
\rho_{\om}(\lam)=\left\{\begin{array}{ll}
                          1-\lam\overline{\om},     &
\mbox{ if }\Omega_+=\dD; \\
\\
\lam+\overline{\om} & \mbox{ if }\  \Omega_+=\Pi_+.
\end{array}\right.
\]
Thus, in both cases
\[
\Omega_+=\{\omega\in\dC:\rho_\omega(\omega)>0\}
\] and
$\Omega_0=\{\omega\in\dC:\rho_\omega(\omega)=0\}$ is the boundary of
$\Omega_+$. Correspondingly we set
\begin{equation}
\label{eq:omegaminus} \Omega_- =\dC \setminus (\Omega_+
\cup\Omega_0) = \{\om\in\dC: \rho_{\om}(\om) < 0\}.
\end{equation}
The normalized standard inner product $\langle f, g\rangle_{nst}$ is defined as
$$
\langle f, g\rangle_{nst}=\left\{\begin{array}{l}\ \frac{1}{2\pi}
\int_0^{2\pi}g(e^{i\theta})^*f(e^{i\theta})d\theta\quad\text{if} \
\Omega_0=\DD\\
{}\\
\ \frac{1}{2\pi}\int_{-\infty}^\infty g(i\mu)^*f(i\mu)d\mu \quad\text{if}\
\Omega_0=i\RR.\end{array}\right.
$$
The symbol $H_2^{\ptq}$ (resp., $H_\infty^{\ptq}$)  stands for the
class of $\ptq$ mvf's with entries in the Hardy space $H_2$ (resp.,
$H_\infty$); $H_2^p$  is short for $H_2^{p\times 1}$ and
$(H_2^p)^\perp$ is the orthogonal complement of $H_2^p$ in $L_2^p$
with respect to the standard inner product on $\Omega_0$.

Most of the other notation that we use will be fairly standard:
\bigskip

\noindent mvf for matrix valued function, vvf for vector valued function,
and $\cR$ for rational mvf's;
$\ker A$ and $\ran A$ for the kernel and range of a matrix $A$, and, if $A$
is square,
$\sigma (A)$ for its spectrum and
$\nu_-(A)$  (resp., $\nu_+(A)$) for the number of its negative
(resp., positive) eigenvalues
(counting multiplicities). If $f(\lam)$ is a mvf, then
$$
f^{\#}(\lam)=f(\lam^\circ)^*, \mbox{ where }\lam^\circ=\left\{\begin{array}{ll}
                                                        1/\overline{\lam}  & :\,\mbox{ if }\Omega_+=\dD,\,\lam\ne 0; \\
                                                        -\overline{\lam}  & :\,\mbox{ if }\Omega_+=\Pi_+ \\
                                                      \end{array}\right.
$$
$$
\gh_f=\{\lambda\in\CC\ \textrm{at which $f(\lambda)$ is
holomorphic}\}\quad\textrm{and}\quad \gh_f^\pm=\gh_f\cap\Omega_\pm.
$$

By a fundamental result of  Kre\u{\i}n and Langer~\cite{KL}, every
generalized Schur function $s\in
\cS_{\kappa}^{\ptq}:=\cS_{\kappa}^{\ptq}(\Omega_+)$ admits a pair of
coprime factorizations
\begin{equation}\label{KL}
s(\lam)=b_\ell(\lam)^{-1}s_\ell(\lam)=s_r(\lam)b_r(\lam)^{-1} \quad (\lam\in{\mathfrak
H}_s^+),
\end{equation}
where $b_\ell$ and $b_r$ are Blaschke--Potapov products of sizes
$\ptp$ and $\qtq$, respectively, of degree $\kappa$, and the mvf's
$s_\ell$ and $s_r$ both belong to the {\it Schur class}
$\cS^{\ptq}:=\cS_0^{\ptq}(\Omega_+)$.

Various interpolation problems in the class of generalized Schur
mvf's were considered in~\cite{T24}, \cite{Nud81}, \cite{Gol83},
 \cite{BH83}, \cite{BH86}, \cite{BGR}, \cite{ACDR}, \cite{Der03},
\cite{Bol04}.

The data  for the interpolation problem we consider in this paper for $\Omega_+=\DD$ is
coded into a set of four matrices, $M$, $N$, $P\in\CC^{\ntn}$ and $C\in\dC^{m\times n}$,
that are subject to the following constraints:
\begin{enumerate}
\item[\rm(A1)] The pencil ${M}-\lam{N}$ is invertible on
$\Omega_0=\TT$, i.e., the resolvent set of
this pencil,
$$
\rho({M},{N})\stackrel{def}{=}\{\lambda\in\CC:\, \det({M}-\lam{N})\ne 0\}
\quad\textrm{contains}\quad \Omega_0.
$$
\vspace{1mm}
\item[\rm (A2)] $P$ is  a Hermitian solution of the Lyapunov-Stein equation
\begin{equation}\label{eq:0.LS}
    {M}^* P{M}-{N}^*P{N}=C^*j_{pq}C,
\end{equation}
where
\begin{equation}\label{eq:jpq}
j_{pq}=\begin{bmatrix} I_p & 0\\
                        0 & -I_q
        \end{bmatrix},  \quad p>0, \quad q>0 \quad\textrm{and}\ m=p+q.
\end{equation}
\item[\rm(A3)] The matrices $C_1$, $C_2$ determined by the decomposition
\begin{equation}\label{eq:C}
    C=\begin{bmatrix} C_1\\
                        C_2
        \end{bmatrix},\quad C_1\in\dC^{p\times
n},\quad C_2\in\dC^{q\times n},
\end{equation}
satisfy the rank conditions
\begin{equation}\label{eq:A3}
\textup{rank}
\begin{bmatrix}M-\lambda N\\
C_2\end{bmatrix}=n\quad\textrm{and} \quad \textup{rank}\begin{bmatrix}\lambda  M-N\\
C_1\end{bmatrix}=n \quad\text{for every}\ \lambda\in\Omega_+.
\end{equation}
\item [\rm(A4)]
There exists an $n\times n$ Hermitian matrix $X$ that meets the
conditions specified in (B4) below.
\end{enumerate}
In particular, it follows from (A3) that the triple $(C,M,N)$ is observable:
\begin{equation}\label{eq:A4}
       \bigcap_{\lam\in\rho({M},{N})}\ker C({M}-\lam
        {N})^{-1}=\{0\},
\end{equation}
see Proposition \ref{ChatJ}, and for additional discussion of observability
and controllabilty of a fairly general class of pencils, Theorem 3.5
of  \cite{AD96}.

The basic {\bf bitangential interpolation} problem under consideration
for $\Omega_+=\DD$
corresponds to the decompositions:
\begin{equation}\label{eq:0.MN}
    {M} =\left[
   \begin{array}{cc}
   A_1&0\\
   0& I_{n_2}
   \end{array}\right],\quad
  N=\left[
   \begin{array}{cc}
   I_{n_1}&0\\
   0& A_2
   \end{array}\right]
\end{equation}
and
\[
C=\begin{bmatrix}
  C_{11} & C_{12} \\
  C_{21} & C_{22} \\
\end{bmatrix}:
\begin{bmatrix}
  \dC^{n_1} \\
  \dC^{n_2} \\
\end{bmatrix}\to
\begin{bmatrix}
  \dC^{p} \\
  \dC^{q} \\
\end{bmatrix},\quad\textrm{where } n=n_1+n_2,\quad n_1>0,
\quad n_2>0,
\]
$A_1\in\dC^{n_1\times n_1}$ and  $A_2\in\dC^{n_2\times n_2}$. These matrices
are subject to the following constraints:
\begin{enumerate}
\item[(B1)] $\sigma(A_1)\cup\sigma(A_2)\subset \dD$.
\item[(B2)] 
$P$ is  a Hermitian solution of the Lyapunov-Stein equation (\ref{eq:0.LS})
 with $M$ and $N$ as in (\ref{eq:0.MN}).
\item[(B3)]
The pairs $(C_{12}, A_2)$ and $(C_{21}, A_1)$ are observable.
\item [(B4)]
There exists an $n\times n$ Hermitian matrix $X$ such that:
\begin{enumerate}
    \item [(i)] $XPX=X$.
\vspace{2mm}
    \item [(ii)] $PXP=P$.
\vspace{2mm}
    \item [(iii)] $\ran\ X$ is invariant for $M$ and $N$, i.e.,
$Mx\in\ran X$ and $Nx\in\ran X$ if $x\in\ran X$.
 \end{enumerate}
\end{enumerate}
If $P$ is invertible, then (B4) is superfluous, since it is automatically
satisfied by $X=P^{-1}$.

The {\bf one sided tangential interpolation} problem corresponds to the case
when either $n_1=0$ or $n_2=0$ (but not both). In these two cases the formulations in (B1) and (B3) must be interpreted properly.  For ease of future reference
we summarize the changes in the following remark:

\begin{remark}
\label{rem:oct29a8} If $n_2=0$, then $n_1=n>0$, $M=A_1$, $N=I_n$, (B1) reduces to
$\sigma(A_1)\subset\DD$ and (B3) to $(C_2, A_1)$ is observable.

If $n_1=0$, then $n_2=n>0$, $M=I_n$, $N=A_2$, (B1) reduces to
$\sigma(A_2)\subset\DD$ and (B3) to $(C_1, A_2)$ is observable.

In both of these cases $C_1$ and $C_2$ are determined by the decomposition~\eqref{eq:C}.
\end{remark}

The {\bf bitangential interpolation problem} corresponding to the
data set $M,N,C,P$ is formulated in terms of the Krein-Langer
factorizations (\ref{KL}) of the mvf $s\in \cS^{\ptq}_\kappa$, the
mvf
$$F(\lambda)=C({M}-\lam {N})^{-1}$$
and the Hermitian
matrix $P_s\in \dC^{n\times n}$ that is defined by
formula~\eqref{eq:3.7} as follows:\\
Describe the set $\widehat{\cS}_\kappa(M, N, C,P)$ of mvf's $s\in
\cS^{\ptq}_\kappa$ which satisfy the three conditions:
\begin{enumerate}
    \item [(C1)] $\begin{bmatrix}
      b_\ell & -s_\ell \
    \end{bmatrix}Fu\in H^p_2$ for every $u\in \dC^n$;
\vspace{2mm}
    \item [(C2)]
    $\begin{bmatrix}
      -s_r^* & b_r^* \
    \end{bmatrix}Fu\in (H_2^q)^\perp$ for every $u\in \dC^n$;
\vspace{2mm}
    \item [(C3)] $P_s\leq P$.
\end{enumerate}

The set of mvf's $s\in\widehat{\cS}_\kappa(M, N, C,P)$ for which the
equality $P_s=P$ prevails in (C3) will be denoted by
${\cS}_\kappa(M, N, C,P)$. We will show in Theorem~\ref{thm:3.1}
that $P_s$ is a solution of the Lyapunov-Stein
equation~\eqref{eq:0.LS} and  $P_s=P$ for every
$s\in\widehat{\cS}_\kappa(M, N, C,P)$ and hence
$$
\widehat{\cS}_\kappa(M, N, C,P)={\cS}_\kappa(M, N, C,P).
$$
We shall also write
$$
{\cS}_\kappa(M, N, C,P)=\left\{\begin{array}{ll}{\cS}_\kappa(A_1,
A_2, C,P) &
\quad\textrm{if}\quad n_1>0\quad\textrm{and}\quad n_2>0,\\
{\cS}_\kappa(A_1, C,P) & \quad\textrm{if}\quad n_1=n\quad\textrm{and}\quad n_2=0,\\
{\cS}_\kappa(A_2, C,P) & \quad\textrm{if}\quad
n_1=0\quad\textrm{and}\quad n_2=n.
\end{array}\right.
$$
If $P$ is invertible, then
\begin{equation}\label{eq:0.6}
\cS_\kappa(M, N,C,P)\ne\emptyset\Longleftrightarrow     \kappa \ge
\nu_-(P).
\end{equation}
Moreover, if $P$ is invertible and $\nu_-(P)=\kappa_1$, then
\begin{equation}\label{eq:descr}
    {\cS}_\kappa(M, N, C,P)=T_W[\cS_{\kappa-\kappa_1}^{\ptq}]\cap
\cS_{\kappa}^{\ptq},
\end{equation}
where
\begin{equation}\label{eq:0.2}
    T_W[\cS_{\kappa_2}^{\ptq}]=\left\{T_W[\varepsilon]:\,\varepsilon\in \cS_{\kappa_2}^{\ptq}\right\}
\end{equation}
is the range of the linear fractional transformation
\begin{equation}\label{eq:0.3}
    T_W[\varepsilon]=(w_{11}\varepsilon+w_{12})(w_{21}\varepsilon+w_{22})^{-1}
\end{equation}
based on the $m\times m$ mvf
\begin{equation}\label{eq:0.5}
    W(\lam)=I_m-\rho_\mu(\lam)F(\lam)XF(\mu)^*j_{pq},
\end{equation}
with $\mu\in\Omega_0$, $X=P^{-1}$ and block decomposition
\[
W(\lam)=\begin{bmatrix} w_{11} & w_{12}\\
                        w_{21} & w_{22}
        \end{bmatrix}
\]
that is conformal with $j_{pq}$ in~\eqref{eq:jpq}. The mvf $W$  defined by
(\ref{eq:0.5}) belongs to the class $\cU_{\kappa_1}(j_{pq})$ of $m\times m$ mvf's $W$
meromorphic in $\Omega_+$ such that the kernel
\begin{equation}\label{eq:0.4}
   {\mathsf K}_\omega^W(\lambda):=\frac{j_{pq}-W(\lambda)j_{pq}
W(\omega)^*}{\rho_\omega(\lambda)}
\end{equation}
has $\kappa_1$ negative squares on $\gh_W^+\times\gh_W^+$ and
$ W(\mu)j_{pq} W(\mu)^*
=j_{pq}$ a.e. on $\Omega_0$.
The class of mvf's $W\in\cU_\kappa(j_{pq})$ which satisfy
\begin{equation}\label{eq:0.3a}
    s_{21}:=-w_{22}^{-1}w_{21}\in \cS_\kappa^{q\times p}.
\end{equation}
will be denoted $\cU_\kappa^\circ(j_{pq})$. Here $s_{21}$ is the lower left
hand corner block in the Potapov-Ginzburg transform
$$
S=\begin{bmatrix}s_{11}&s_{12}\\ s_{21}&s_{22}\end{bmatrix}
$$
of $W$ that is defined in (\ref{PGtrans}).

The characterization of the set $T_W[\cS_{\kappa_2}^{\ptq}]$ given
in Theorem~\ref{thm:2.2} was obtained in~\cite{DerDym08}; it is a
generalization to an indefinite setting of a result from~\cite{D8}.
The formulation of the result requires some facts about the
reproducing kernel de Branges-Kre\u{\i}n space $\cK(W)$ associated
with the kernel ${\mathsf K}_\omega^W(\lambda)$ and an indefinite
analog of the de Branges-Rovnayk space $\sD(s)$, developed
in~\cite{AD86},  \cite{ADRS} and \cite{DerDym08}. The needed facts
are reviewed in the next section  for the convenience of the reader.

Then, as follows from Theorem~\ref{thm:2.2}, the set of solutions of the
problem $(C1)-(C3)$ when $P$ is not invertible can be
described via a formula that is
similar to \eqref{eq:descr} (see Theorem~\ref{thm:3.13}).

Finer analysis of the set of solutions is connected with the
factorization of the resolvent matrix $W(\lam)$ that is presented in
Theorem \ref{factofW}. In Lemma~\ref{lem:7.6a} we show that if
(B1)-(B4) are in force, then  $s_{21}\in \cS_\kappa^{q\times p}$ and
hence $W\in\cU_\kappa^\circ(j_{pq})$. Let  $s_{21}$ have
Kre\u{\i}n-Langer factorizations
\begin{equation}\label{eq:0.4a}
    s_{21}(\lambda):=\gb_{\ell}(\lambda)^{-1}\gs_{\ell}(\lambda)
=\gs_r(\lambda)\gb_r(\lambda)^{-1},\quad\lambda\in\gh_s^+,
\end{equation}
where $\gb_{\ell}$, $\gb_{r}$ are Blashke-Potapov products of degree
$\kappa$ and $\gs_{\ell},\gs_r\in \cS^{\qtp}$. Then, as was shown
in~\cite{DerDym08}, the mvf's  $s_{11}\gb_r$ and $\gb_{\ell}s_{22}$
belong to the classes $\cS^{\ptp}$ and $\cS^{\qtq}$, respectively.
Therefore, they admit inner-outer and outer-inner factorizations
\begin{equation}\label{eq:0.5a}
    s_{11}\gb_r=b_1\varphi_1,\quad\gb_{\ell}s_{22}=\varphi_2b_2,
\end{equation}
where $b_1\in \cS_{in}^{\ptp}$, $b_2\in \cS_{in}^{\qtq}$,
$\varphi_1\in \cS_{out}^{\ptp}$, $\varphi_2\in \cS_{out}^{\qtq}$. In
keeping with the terminology used in \cite{ArovD08}, the pair
$\{b_1,b_2\}$ is called an {\it associated pair} of the mvf
$W\in\cU_\kappa^\circ(J)$ and denoted as $\{b_1,b_2\}\in ap(W)$. In
the definite case the formulas in \eqref{eq:0.5a} simplify to the
inner-outer factorization of $w_{11}^\#$ and the outer-inner
factorization of $w_{22}^{-1}$ (see~\cite{Arov93},~\cite{ArovD97}).

Using the results of~\cite{DerDym08}  we show that the
matrices $\begin{bmatrix}w_{11} & w_{12}\end{bmatrix}$ and
 $\begin{bmatrix}w_{21} & w_{22}\end{bmatrix}$ admit coprime
 factorization over $\Omega_-$ and  $\Omega_+$, respectively:
\begin{equation}\label{eq:0.8}
 \begin{bmatrix}w_{11} & w_{12}\end{bmatrix}
=b_1\begin{bmatrix}\wt\varphi_{11}& \wt\varphi_{12}\end{bmatrix}
 =(b_1^\#)^{-1}\begin{bmatrix}\wt\varphi_{11}& \wt\varphi_{12}\end{bmatrix},
\end{equation}
\begin{equation}\label{eq:0.7}
 \begin{bmatrix}w_{21}& w_{22}\end{bmatrix}=
b_2^{-1}\begin{bmatrix}\varphi_{21}&\varphi_{22}\end{bmatrix},
\end{equation}
where $\{b_1,b_2\}\in\mbox{ap}(W)$,  $\wt\varphi_{11}\in \cR\cap H_\infty^{p\times
p}(\Omega_-)$,  $\wt\varphi_{12}\in \cR\cap H_\infty^{p\times q}(\Omega_-)$,
$\varphi_{21}\in \cR\cap H_\infty^{q\times p}$ and $\varphi_{22}\in \cR\cap
H_\infty^{q\times q}$.

Then applying the Kre\u{\i}n-Langer generalization of Rouche's
theorem we show that for every mvf $\varepsilon\in
\cS_{\kappa-\kappa_1}^{p\times q}$ with the Kre\u{\i}n-Langer
factorizations
\begin{equation}\label{eq:0.9}
\varepsilon=\theta_\ell^{-1}\varepsilon_\ell=\varepsilon_r\theta_r^{-1}
\end{equation}
 the mvf's $\theta_\ell
\wt\varphi_{11}^\#+\varepsilon_\ell\wt\varphi_{12}^\#$ and
$\varphi_{21}\varepsilon_r+\varphi_{22}\theta_r$ have exactly
$\kappa$ zeros in $\Omega_+$.

The main results of this paper is the following description of
$\cS_\kappa(M, N,C,P)$ when $\Omega_+=\DD$ and the analogous
description for $\Omega_+=\Pi_+$ that is presented in Section
\ref{sec:bipinrhp}:

\begin{thm}\label{thm:0.2}
Let (B1)--(B4) be in force, let $\nu_-(P)=\kappa$ and let
\begin{equation}
\label{nu}
    \nu=\rank(M^*P^2M+N^*P^2N+C^*C)-\rank P.
\end{equation}
Then there are unitary matrices $U\in \dC^{p\times p}$, $V\in
\dC^{q\times q}$, such that
$s\in \cS_\kappa(M, N,C,P)$ if and only if $s$ belongs to
$\cS_\kappa^{\ptq}$ and is of the form $s=T_W[\varepsilon]$, where
\begin{equation}\label{epsilon}
    \varepsilon=U\left[\begin{array}{cc}
                                  \wt \varepsilon & 0 \\
                                        0 & I_\nu \\
                                      \end{array}  \right]V^*,\quad
\textrm{and}\quad \wt\varepsilon \in \cS^{(p-\nu)\times(q-\nu)}.
\end{equation}
If $\varepsilon\in\cS^{\ptq}$, then $T_W[\varepsilon]\in
\cS_\kappa^{\ptq}$ if and only if
\begin{enumerate}
    \item [(a)]
     the factorization
    $w^\#_{11}+\varepsilon w^\#_{12}
    =( \wt\varphi_{11}^\#+\varepsilon\wt\varphi_{12}^\#)b_1^{-1}$ is
    coprime over $\Omega_+$  and
\vspace{2mm}
    \item [(b)] the factorization
    $w_{21}\varepsilon+ w_{22}=b_2^{-1}(\varphi_{21}\varepsilon
+\varphi_{22})$ is
    coprime over $\Omega_+$.
\end{enumerate}
\end{thm}

\begin{proof}
See subsection \ref{subsec:prfthm:0.2}.
\end{proof}

If $P$ is invertible, the statement of this theorem is simpler, since
$\nu=0$ (see Corollary~\ref{cor:C13}), and
we can also treat the case when $\nu_-(P)<\kappa$.

\begin{thm}\label{thm:0.1}
Let the data set $(M, N,C,P)$ satisfy the assumptions (B1)-(B3), let $P$ be invertible,
$\kappa_1=\nu_-(P)\le\kappa$, and let the mvf's $W$, $b_1$, $b_2$ be defined
by~\eqref{eq:0.5}, \eqref{eq:0.7}, \eqref{eq:0.8}. Then:
\begin{enumerate}
\item[\rm(I)] $s\in S_\kappa(M, N,C,P)$
if and only if $s=T_W[\varepsilon]$, where $\varepsilon\in
\cS_{\kappa-\kappa_1}^{p\times q}$ and
$T_W[\varepsilon]\in\cS_\kappa^{\ptq}$.
\item[\rm(II)] If  $\varepsilon\in
\cS_{\kappa-\kappa_1}^{p\times q}$ and $\theta_\ell$, $\theta_r$,
$\varepsilon_\ell$, $\varepsilon_r$ are a choice of its
Kre\u{\i}n-Langer factorizations
as in \eqref{eq:0.9},
  then $T_W[\varepsilon]\in
\cS_\kappa^{p\times q}$  if and only if the factorizations
\begin{equation}\label{Reg1}
    \theta_\ell w^\#_{11}+\varepsilon_\ell w^\#_{12}
    =(\theta_\ell
    \wt\varphi_{11}^\#+\varepsilon_\ell\wt\varphi_{12}^\#)b_1^{-1},
\end{equation}
\begin{equation}\label{Reg2}
    w_{21}\varepsilon_r+ w_{22}\theta_r=b_2^{-1}(\varphi_{21}\varepsilon_r
+\varphi_{22}\theta_r)
\end{equation}
are coprime over $\Omega_+$.
\end{enumerate}
\end{thm}
\begin{proof}
See subsection \ref{subsec:prfthm:0.1}.
\end{proof}

The set of mvf's $s\in \cS^{\ptq}_\kappa$ which satisfy (C1)-(C3)
and the supplementary condition
 \begin{enumerate}
    \item [(C4)] $s$ is holomorphic in $\sigma(A_1)\cup\overline{\sigma(A_2)}$,
\end{enumerate}
is denoted by ${\cT\cN}_\kappa(M, N, C,P)$, and $s\in {\cT\cN}_\kappa(M, N, C,P) $ is
said to be a
solution of the {\it Takagi-Nudelman problem}.

For the Takagi-Nudelman problem the conditions (\ref{Reg1}) and (\ref{Reg2})
of Theorem~\ref{thm:0.1} are replaced by the
single condition:
\begin{equation}
\label{eq:sep12b8} (\varphi_{21}\varepsilon_r+ \varphi_{22}\theta_r)^{-1}
\quad\textrm{is
holomorphic in}\  \sigma(A_1)\cup\overline{\sigma(A_2)}.
\end{equation}
which reduces to a condition in \cite{Nud81} when $\varepsilon\in\cS^{\ptq}$
and reduces to (19.2.6) in
\cite{BGR} when $\varepsilon$ is a rational mvf in $\cS^{\ptq}$.

In the scalar case ($p=q=1$, $n_2=0$, $\nu_-(P)=\kappa$) the two
conditions (\ref{Reg1}) and (\ref{Reg2}) are in force if and only if
(\ref{eq:sep12b8}) holds. Consequently, every solution $s$ from the
set $ \cS_\kappa(A_1,C,P)$ is holomorphic on $\sigma(A_1)$. A new
effect which is revealed in the matrix case is that there are mvf's
$s$ which belong to $\cS_\kappa(A_1,A_2,C,P)$ but are not
holomorphic on $\sigma(A_1)\cup\overline{\sigma(A_2)}$ (see
Example~\ref{ex:6.8}). Thus, there are mvf's that satisfy the
conditions (\ref{Reg1}) and (\ref{Reg2}) but do not satisfy
(\ref{eq:sep12b8}). Therefore, the inclusion
\[
{\cT\cN}_\kappa(M,N,C,P)\subseteq \cS_\kappa(M,N,C,P)
\]
can be proper.\\
The parameter $\varepsilon\in \cS^{\ptq}_{\kappa-\kappa_1}$ in
\eqref{eq:0.3} is said to be {\it excluded} for the problem
(C1)-(C4) 
if $s=T_W[\varepsilon]\not\in{\cT\cN}_\kappa(M,N,C,P)$. If
$\sigma(A_1)\cap\overline{\sigma(A_2)}=\emptyset$, the condition
\[
{\mathsf K}_\omega^W(\omega)\geqslant 0
\ \text{for all}\ \omega\in \sigma(A_1)^o\cup
\overline{\sigma{(A_2 )}}\,,
\]
which is formulated in terms of the
kernel~\eqref{eq:0.4}, suffices to guarantee that
the problem (C1)-(C4) has no excluded
parameters. In the scalar case
sufficient
conditions for the Nevanlinna-Pick problem in the generalized Nevanlinna
class to have no excluded parameters
were found
in~\cite{DL97} (see also~\cite{ABD98} and \cite{AmDer} for the matrix
case).

In subsection~\ref{TSproblem}, given a data set $(M,N,C,P)$ satisfying the assumptions
(B1)-(B3),  we will consider an associated pair $\{b_1,b_2\}$ for the mvf $W$ and a
rational mvf $K$ holomorphic in $\Omega_+$ such that
\begin{equation}\label{eq:0.11}
    \mbox{the mvf }b_1^{-1}(s-K)b_2^{-1}\mbox{ has exactly $\kappa$ poles
(counting multiplicities) in
    $\Omega_+$}
\end{equation}
for every $s\in \cS_\kappa(M,N,C,P)$. Every mvf $s\in
\cS_{\kappa'}^{p\times q}$ $(\kappa'\le \kappa)$ which
satisfies~\eqref{eq:0.11} is called a solution of the Takagi-Sarason
problem with data set $(b_1,b_2,K)$, and the symbol
${\cT\cS}_\kappa(b_1,b_2,K)$ is used to denote the set of such  solutions
$s$. We shall show that for $(b_1,b_2,K)$ corresponding to the
problem (C1)-(C3)
${\cT\cS}_\kappa(b_1,b_2,K)=T_W[\cS_{\kappa-\nu_-(P)}^{p\times q}]$ and,
hence, that
\[
{\cT\cN}_\kappa(M,N,C,P)\subseteq \cS_\kappa(M,N,C,P)\subseteq
{\cT\cS}_\kappa(b_1,b_2,K).
\]
Rational solutions of the Takagi-Nudelman and Takagi-Sarason
problems in the case of invertible $P$ have been described earlier
in~\cite{BGR}.

 The paper is organized as follows. In Section
\ref{preli} the basic notions are introduced. Section 3 focuses on the
bitangential interpolation problem in the open unit disc $\DD$. Pole and
zero multiplicities and a factorization formula for the resolvent matrix for
the interpolation problem considered in Section 3 are developed in Section 4.
Theorems \ref{thm:0.2} and \ref{thm:0.1}, which provide parametrizations of
the set of all
solutions to this problem both when $P$ is invertible and $P$ is not
invertible are completed in the first part of Section 5. The latter part
discusses the Tagaki-Nudelman problem and excluded parameters. The
overall strategy for analyzing the
bitangential interpolation problem in the open right half plane $\Pi_+$
(and the open upper half plane $\CC_+$) is much the same as for $\DD$. The
changes in the formulas and the main conclusions for $\Omega_+=\DD$ are
discussed briefly in Section 6, without proof.

\section{Preliminaries \label{preli}}

\subsection{The generalized Schur class}
Recall that a Hermitian kernel  ${\mathsf
K}_\omega(\lambda):\Omega\times\Omega\to\dC^{m\times m}$ is said to
have $\kappa $ negative squares and is written as
\[
\mbox{sq}_-{\mathsf K}=\kappa,
\]
if for every positive integer $n$ and every choice of $\omega_j\in\Omega$
and $u_j\in\dC^m$
$(j=1,\dots,n)$ the matrix
\[
\left(\left<{\mathsf
K}_{\omega_j}(\omega_k)u_j,u_k\right>\right)_{j,k=1}^n
\]
has at most $\kappa$ negative eigenvalues and for some choice of
$\omega_j\in\Omega$ and $u_j\in\dC^m$ exactly $\kappa$ negative
eigenvalues.

The class $\cS^{\ptq}:=\cS_{0}^{\ptq}(\Omega_+)$ is the usual Schur
class. Recall that a mvf $s\in \cS^{p\times q}$ is called inner
(resp., $*$-inner), if $s(\mu)$ is an isometry (resp., co-isometry)
for a.e. $\mu\in\Omega_0$, that is
\[
I_q-s(\mu)^*s(\mu)=0 \quad (\mbox{resp., } I_p-s(\mu)s(\mu)^*=0)
\quad\textrm{for almost all points}\  \mu\in\Omega_0.
\]
Let $\cS_{in}^{p\times q}$ ($\cS_{*in}^{p\times q}$) denote the set
of all inner (resp., $*$-inner) mvf's $s\in \cS^{p\times q}$. An
example of an inner square mvf is provided by the {\it
Blaschke--Potapov product}, that in the case of the unit disc
($\Omega_+=\dD$) is given by
\begin{equation}\label{BPprod}
b(\lam)=\prod_{j=1}^{\kappa}b_j(\lam),\quad
b_j(\lam)=I-P_j+\frac{\lam-\alpha_j}{1-\bar\alpha_j \lam}P_j,
\end{equation}
where $\alpha_j\in\dD$, $P_j$ are orthogonal projections in $\dC^p$
$(j=1,\dots ,n)$. The factor $b_j$ is called {\it simple} if $P_j$
has rank one. Although $b(\lam)$ in not uniquely represented as a
product of simple factors the number $\kappa$ of these factors is
the same for every representation~\eqref{BPprod}. It is called the
{\it degree} of the Blaschke--Potapov product
$b(\lam)$~\cite{Pot}.

A theorem of  Kre\u{\i}n and Langer~\cite{KL} guarantees that every
generalized Schur function $s\in \cS_{\kappa}^{\ptq}(\Omega_+)$
admits a factorization of the form
\begin{equation}\label{KLleft}
s(\lam)=b_\ell(\lam)^{-1}s_\ell(\lam) \quad \textrm{for}\
\lam\in{\mathfrak h}_s^+,
\end{equation}
where $b_\ell$ is a Blaschke--Potapov product of degree $\kappa$,
$s_\ell$ is in the Schur class $\cS^{\ptq}(\Omega_+)$ and
\begin{equation}\label{KLcanon}
\ker s_\ell(\lam)^*\cap \ker b_\ell(\lam)^*=\{0\}\quad\textrm{for}\ \lam\in\Omega_+.
\end{equation}
The representation~\eqref{KLleft} is called a {\it left
Kre\u{\i}n--Langer factorization}. The assumption~\eqref{KLcanon}
can be rewritten in the equivalent form
\begin{equation}\label{KLcanon2}
\rank \left[\begin{array}{cc}
 b_\ell(\lam) & s_\ell(\lam)
\end{array}\right]
=p\quad \textrm{for}\ \lam\in\Omega_+.
\end{equation}
If $\alpha_j\in\dD$ $(j=1,\dots ,n)$ are all the zeros of $b_\ell$ in $\Omega_+$, then
the noncancellation condition~\eqref{KLcanon} ensures that ${\mathfrak
H}_s^+=\Omega_+\setminus\{\alpha_1,\dots,\alpha_n\}$. The left Kre\u{\i}n--Langer
factorization~\eqref{KLleft} is essentially unique in a sense that $b_\ell$ is defined
uniquely up to a left unitary factor $V\in\dC^{\ptp}$.

 Similarly, every generalized Schur
function $s\in \cS_{\kappa}^{\ptq}(\Omega_+)$ admits a {\it right
Kre\u{\i}n--Langer factorization}
\begin{equation}\label{KLright}
s(\lam)=s_r(\lam)b_r(\lam)^{-1}\quad\textrm{for}\ \lam\in{\mathfrak h}_s^+,
\end{equation}
where $b_r$ is a Blaschke--Potapov product of degree $\kappa$ and
$s_r\in \cS^{\ptq}(\Omega_+)$ satisfies the condition
\begin{equation}\label{KRcanon}
\ker s_r(\lam)\cap \ker b_r(\lam)=\{0\},\quad\textrm{for}\ \lam\in\Omega_+.
\end{equation}
This condition can be rewritten in the equivalent form
\begin{equation}\label{KRcanon2}
\rank \left[\begin{array}{cc}
 b_r(\lam)^* & s_r(\lam)^*
\end{array}\right]  =q\quad\textrm{for}\  \lam\in\Omega_+.
\end{equation}
Under assumption~\eqref{KRcanon} the mvf $b_r$ is uniquely
defined up to a right unitary factor $V'\in\dC^{\qtq}$.

\begin{lem}\label{Corona}
A mvf $s_\ell\in\cS^{p\times q}$ and a finite Blaschke-Potapov
product $b_\ell\in \cS_{in}^{p\times p}$ meet the rank condition
(\ref{KLcanon2}), if and only if there exists a pair of mvf's $c\in
H^{p\times p}_\infty$ and $d\in H^{q\times p}_\infty$ such that
\begin{equation}\label{CorFormula}
b_\ell(\lam)c(\lam)+s_\ell(\lam)d(\lam)=I_p\quad \text{for}\ \lam\in
\Omega_+.
\end{equation}
\end{lem}

Lemma \ref{Corona} is a matrix version of the Carleson Corona Theorem.
A proof, which is adapted from Fuhrmann \cite{Fuhr68} who treated the case
$p=q$,  is furnished in \cite{DerDym08}.

A dual statement for Lemma~\ref{Corona} is obtained by applying
Lemma~\ref{Corona} to transposed vvf's.
\begin{lem}\label{Corona1}
A mvf $s_r\in\cS^{p\times q}$ and a finite Blaschke-Potapov product
$b_r\in \cS_{in}^{q\times q}$ meet the rank condition
(\ref{KRcanon2}), if and only if there exists a pair of mvf's $c\in
H^{q\times q}_\infty$ and $d\in H^{q\times p}_\infty$ such that
\begin{equation}\label{CorFormula1}
c(\lam)b_r(\lam)+d(\lam)s_r(\lam)=I_q\quad \text{for}\ \lam\in
\Omega_+.
\end{equation}
\end{lem}
The factorization~\eqref{KLleft} is called a {\it left coprime} factorization of $s$ if $s_\ell$ and
$b_\ell$ satisfy~\eqref{CorFormula}.
Similarly,  the factorization~\eqref{KLright} is called a {\it right coprime} factorization of $s$
if $s_r$ and $b_r$ satisfy~\eqref{CorFormula1}.

Every vvf $h(\lambda)$ from $H_2^p$ $({H_2^p}^\perp)$  has
nontangential limits $h(\mu)$ $(\mu\in\Omega_0)$ a.e. on the
boundary $\Omega_0$. These nontangential limits identify the vvf
$h$ uniquely. In what follows we often identify a vvf $h\in H_2^p
({H_2^p}^\perp)$ with its boundary values $h(\mu)$.

Let
$P_+$ and $P_-$ denote the orthogonal projections from $L_2^k$ onto $H_2^k$
and ${H_2^k}^\perp$, respectively, where $k$ is a positive integer that will
be understood from the context. The Hilbert spaces
\begin{equation}
\label{eq:feb10a8}
\cH(b_r)=H_2^q\ominus b_rH_2^q, \quad \cH_*(b_\ell):=(H_2^p)^\perp\ominus
b_\ell^*(H_2^p)^\perp
\end{equation}
and the operators
\begin{equation}
\label{eq:feb10b8}
X_r: h\in\cH(b_r)\mapsto P_-sh\quad \quad X_\ell: h\in\cH_*(b_\ell)
\mapsto P_+s^*h
\end{equation}
based on $s\in \cS_{\kappa}^{\ptq}$ will play an important role.

\begin{lem}\label{Ker*}{\rm(cf.~\cite{Der01})}
If $s\in \cS_{\kappa}^{\ptq}$, then:
\begin{enumerate}
\item[(i)]  The operator $X_\ell$ maps $\cH_*(b_\ell)$ injectively
onto $\cH(b_r)$.
\vspace{2mm}
\item[(ii)]
 The operator $X_r$ maps $\cH(b_r)$ injectively onto $\cH_*(b_\ell)$.
\vspace{2mm}
\item[(iii)] $X_\ell=X_r^*$.
\end{enumerate}
\end{lem}
\begin{proof} (i) For every  $h\in \cH_*(b_\ell)$ and $f=b_rh_+\in b_r H_2^q$,
it  is readily checked that
\[
\begin{split}
\left\langle P_+s^*h,f\right\rangle_{nst} &=\left\langle
s^*h,b_rh_+\right\rangle_{nst}\\
&=\left\langle h,s_rh_+\right\rangle_{nst}=0\quad \forall h_+\in
H_2^q\,,
\end{split}
\]
i.e., $X_\ell$ maps $\cH_*(b_\ell)$ into $\cH(b_r)$.
Therefore, since $\cH_*(b_\ell)$ and $\cH(b_r)$ are finite
dimensional spaces of the same dimension, and
\[
\dim\cH_*(b_\ell)=\dim\ker X_\ell+\dim\mbox{range } X_\ell,
\]
it suffices to show that
$\textup{ker}X_\ell=\{0\}$. But, if $h\in
\cH_*(b_\ell)$ and  $P_+s^* h=0$, then $b_\ell h\in H_2^p$ and, in view of
Lemma~\ref{Corona},
\begin{eqnarray*}
  \xi^* (b_\ell h)(\omega) &=&\langle b_\ell h,\frac{\xi}{\rho_\omega}\rangle_{nst}
   = \langle b_\ell h,(b_\ell c+s_\ell d)\frac{\xi}{\rho_\omega}\rangle _{nst}\\
  & =& \langle b_\ell h,s_\ell d\frac{\xi}{\rho_\omega}\rangle _{nst}
  =\langle P_+s^* h, d\frac{\xi}{\rho_\omega}\rangle_{nst}=0
\end{eqnarray*}
for every choice of $\omega\in \Omega_+$, and $\xi\in \dC^p$. Since
$b_\ell(\omega)\not\equiv 0$, this implies that $h(\omega)\equiv 0$.

 Statement (ii) can be obtained by similar calculations, and (iii) is easy.
 \end{proof}


\begin{definition}\label{GammaS}  Let
    \begin{equation}\label{GammaS2}
    \Gamma_\ell : f\in L_2^q \to  X_\ell^{-1}P_{\cH(b_r)}f\in\cH_*(b_\ell)\quad
\text{and}\quad
    \Gamma_r: g\in L_2^p\to X_r^{-1}P_{\cH_*(b_\ell)}g\in \cH(b_r),
\end{equation}
where $X_\ell$ and $X_r$ are defined in formula (\ref{eq:feb10b8}) and $P_{\mathcal{X}}$
denotes the orthogonal projection of $L_2^k$ onto a closed subspace $\mathcal{X}$ of
$L_2^k$ (where $k$ will always be understood from the context).
\end{definition}
It is readily checked that
\begin{equation}\label{GfGg}
 P_+s^*\Gamma_\ell f=P_{\cH(b_r)}f\quad\textrm{and}\quad
P_-s\Gamma_rg=P_{\cH_*(b_\ell)}g\quad\textrm{for}
 \ f\in L_2^q\ \textrm{and}\ g\in L_2^p.
\end{equation}

    \begin{lem}\label{Ker*2}
The operators  $\Gamma_\ell$ and $\Gamma_r$ satisfy the equalities:
\begin{enumerate}
\item[(i)]
$\Gamma_\ell^*=\Gamma_r$; \vspace{2mm}
\item[(ii)]
  $\langle
(\overline{\psi}\Gamma_\ell-\Gamma_\ell\overline{\psi})f, g\rangle_{nst}=0$ for
$f\in\cH(b_r)$, $g\in\cH_*(b_\ell)$ and $\psi$ a scalar inner function.
\end{enumerate}
\end{lem}
\begin{proof} (i)
Let $f\in L_2^q$ and $g\in L_2^p$. Then, since $X_\ell$ maps $\cH_*(b_\ell)$ onto
$\cH(b_r)$ and  $X_\ell^*=X_r$,
$$
\left<\Gamma_\ell f,g\right>_{nst}= \left<X_\ell
^{-1}P_{\cH(b_r)}f,P_{\cH_*(b_\ell)}g\right>_{nst}
    =\left<P_{\cH(b_r)}f, X_r^{-1}P_{\cH(b_\ell )}g\right>_{nst}
    =\left<f, \Gamma_r g\right>_{nst}.
$$
(ii) If $f\in\cH(b_r)$ and $g\in\cH_*(b_\ell )$, then
\[
\begin{split}
\langle (\overline{\psi}\Gamma_\ell -\Gamma_\ell \overline{\psi})f,
g\rangle_{nst}&=
\langle \overline{\psi}\Gamma_\ell f,P_-s\Gamma_rg\rangle_{nst} -\langle
P_+s^*\Gamma_\ell f,
\psi\Gamma_rg\rangle_{nst}\\
&= \langle \overline{\psi}\Gamma_\ell f,s\Gamma_rg \rangle_{nst}
-\langle s^*\Gamma_\ell
f, \psi\Gamma_rg\rangle_{nst}=0.
\end{split}
\]
\end{proof}

\subsection{Reproducing kernel Pontryagin spaces}
In this subsection we review some facts and notation from \cite{AI,Bo} on
the theory of
indefinite inner product spaces for the convenience of the reader.
A linear space $\cK$ equipped with a
sesquilinear form $\left<\cdot,\cdot\right>_\cK$ on $\cK\times\cK$
is called an indefinite inner product space. A subspace $\sF$ of $\cK$ is
called positive (negative) if $\left<f,f\right>_\cK>0\,(<0)$ for
all $f\in\cF$, $f\ne 0$.
If the full space $\cK$ is positive and complete with respect to the norm
$\|f\|=\left<f,f\right>_\cK^{1/2}$ then it is
a Hilbert space.

An indefinite inner product space $(\cK,\left<\cdot,\cdot\right>_\cK)$ is
called a Pontryagin space, if it can be decomposed as the
orthogonal sum
\begin{equation}\label{Pspace}
    \cK=\cK_+\oplus\cK_-
\end{equation}
of a positive subspace $\cK_+$ which is a Hilbert space and a
negative subspace $\cK_-$ of finite dimension. The number
$\mbox{ind}_-\cK :=\dim\cK_-$ is referred to as the negative index
of $\cK$. The convergence in a Pontryagin space
$(\cK,\left<\cdot,\cdot\right>_\cK)$ is meant with respect to the
Hilbert space norm
\begin{equation}\label{Pnorm}
    \|h\|^2=\left<h_+,h_+\right>_\cK-\left<h_-,h_-\right>_\cK\quad
\textrm{when}\    h=h_++h_-\quad\textrm{with}\ h_\pm\in\cK_\pm.
\end{equation}
It is easily seen that the convergence does not depend on a choice of
the decomposition~\eqref{Pspace}.

A Pontryagin space $(\cK,\left<\cdot,\cdot\right>_\cK)$ of
$\dC^m$-valued functions defined on a subset $\Omega$ of $\dC$ is called
a reproducing kernel Pontryagin space associated with the
Hermitian kernel ${\mathsf
K}_\omega(\lambda):\Omega\times\Omega\to\dC^{m\times m}$ if:
\begin{enumerate}
\item for every $\omega\in\Omega$  and  $u\in\dC^m$ the vvf ${\mathsf K}_\omega
(\lambda)u$ belongs to $\cK$;
\vspace{2mm}
\item for every $h\in\cK$,  $\omega\in\Omega$ and $u\in\dC^m$ the
following identity holds
\begin{equation}\label{RKprop}
    \left<h,{\mathsf K}_\omega u\right>_\cK=u^*h(\omega).
\end{equation}
\end{enumerate}

It is known (see~\cite{Sch}) that for every Hermitian
kernel ${\mathsf K}_\omega (\lambda):\Omega\times\Omega\to\dC^{m\times m}$
with a finite number of negative squares on $\Omega\times\Omega$ there is a
unique Pontryagin space $\cK$ with
reproducing kernel ${\mathsf K}_\omega (\lambda)$, and that
$\mbox{ind}_-\cK =\mbox{sq}_-{\mathsf K}=\kappa$. In the case
$\kappa=0$ this fact is due to Aronszajn~\cite{Ar50}.

\subsection{The class $\cU_\kappa(j_{pq})$ and the space ${\cK}(W)$}
Recall (see~\cite{AD86}) that the lower right hand $q\times q$ corner
$w_{22}(\lambda)$
of every $m\times m$ mvf $W\in\cU_\kappa(j_{pq})$ is invertible for all
$\lambda\in\gh_W^+$ except for at most $\kappa$ points. The
Potapov-Ginzburg transform
of $W(\lambda)$, which is  defined on $\gh_W^+$ by the formulas
\begin{equation}\label{PGtrans}
\begin{split}
    S(\lambda)=PG(W)&:=\left[\begin{array}{cc}
      w_{11}(\lambda) & w_{12}(\lambda) \\
      0 & I_q
    \end{array}      \right]
    \left[\begin{array}{cc}
            I_p &       0\\
      w_{21}(\lambda) & w_{22}(\lambda)
    \end{array}      \right]^{-1}\\
    &=\left[\begin{array}{cc}
      I_p & -w_{12}(\lambda) \\
      0   & -w_{22}(\lambda)
    \end{array}      \right]^{-1}
    \left[\begin{array}{cc}
      w_{11}(\lambda) &       0\\
      w_{21}(\lambda) & -I_q
    \end{array}      \right]\,,
\end{split}
\end{equation}
belongs to the class $\cS_\kappa^{\mtm}$. Since
\begin{equation}
\label{eq:jul4a7}
W=PG(S)\Longrightarrow S=PG(W)\,,
\end{equation}
the mvf $W$ is of bounded type. Thus, the nontangential limits $W(\mu)$ exist
a.e. on $\Omega_0$ and  assumption (2) makes sense.
It guarantees that $W(\lambda)$ is
invertible in $\Omega_+$ except for an isolated set of points.
Define  $W$ in $\Omega_-$
by the formula
\begin{equation}\label{eq:1.78}
    W(\lambda)=j_{pq}W^{\#} (\lambda)^{-1} j_{pq}= j_{pq}W(\lambda^\circ)^{-*}
    j_{pq}\quad \text{for}\
    \lambda \in \Omega_-.
\end{equation}
Since the nontangential
limits
\[
W_\pm(t)=\angle\lim_{\lam\to t}\{W(\lam):
  \lam\in\Omega_\pm \}
\]
coincide a.e. in $\Omega_0$,
$W$ in $\Omega_-$ is a pseudo-meromorphic extension of $W$ in $\Omega_+$.
If $W(\lambda)$ is rational this extension is
meromorphic on ${\dC}$. Formula~\eqref{eq:1.78} implies that $W(\lam)$ is
holomorphic and invertible in $\Omega_W:=\gh_W\cap\gh_{W^\#}$.

Let $W\in  \cU_\kappa(j_{pq})$ and let ${\cK}(W)$ be the reproducing
kernel Pontryagin space associated with the kernel ${\mathsf
K}^W_\omega(\lambda)$. The kernel ${\mathsf K}^W_\omega (\lambda)$
extended to ${\Omega_+}\cup{\Omega_-}$ by the
equality~\eqref{eq:1.78} has the same number $\kappa$ of negative
squares. This fact is due to a generalization of the Ginzburg
inequality (see~\cite[Theorem 2.5.2]{ADRS}.

The following example of a ${\cK}(W)$-space will play an important role in
Section~4.
\begin{example}\label{Mspace} Let $C\in \dC^{m\times n}$, let $M,N\in
\dC^{n\times n}$
and assume that
$P\in \dC^{n\times n}$ is an invertible Hermitian matrix, that
$\rho({M},{N})\ne\emptyset$ and that the observability condition (\ref{eq:A4})
is in force. Then the linear space of vvf's
\begin{equation}
\label{eq:3.26}
\cM=\{F(\lam)u:u\in \dC^n\}
\end{equation}
based on the mvf
\begin{equation}\label{eq:1.79}
    F(\lam)=C({M}-\lam {N})^{-1}\quad\textrm{for}\  \lam\in\rho({M},{N})
\end{equation}
and endowed with the inner product
\begin{equation}\label{eq:1.81}
   \langle Fu,\ Fv\rangle_{\cM}=v^*Pu
   \end{equation}
is an RKPS
(reproducing kernel Pontryagin space) with RK (reproducing kernel)
\begin{equation}\label{eq:1.82}
   {\mathsf K}_\omega(\lam)=F(\lam)P^{-1} F(\omega)^*
\end{equation}
and negative index $\textup{ind}_-\cM =\nu_-(P)$.
The assumption~\eqref{eq:A4} insures that the inner
product~\eqref{eq:1.81} is well defined.
\end{example}
We will need the following criterion for the space $\cM$ to be a
${\cK}(W)$ space:

\begin{thm}\label{HW}\cite{D8} Let $M, N, P\in\CC^{\ntn}$,
$C\in\CC^{m\times n}$ and assume that $P$ is invertible, $\rho(M,N)\ne\emptyset$ and that \eqref{eq:A4} holds. Then the RKPS $\cM$ considered in
Example  \ref{Mspace} with kernel ${\mathsf K}^W_\omega(\lambda)$
of the form (\ref{eq:1.82})  is a ${\cK}(W)$ space if and only if $P$
is a solution of
the
equation
\begin{equation}\label{eq:1.83}
   M^*PM-N^*PN=C^*JC\quad \textup{when} \ \Omega_+=\dD,
\end{equation}
\begin{equation}\label{eq:1.83B}
   M^*PN+N^*PM +C^*JC=0\quad \textup{when}\ \Omega_+=\Pi_+,
\end{equation}
where $J$ is a signature matrix ($J=J^*=J^{-1}$).
The mvf $W$ is uniquely defined
by the formula
\begin{equation}\label{eq:1.84}
   W(\lam)=I_m-\rho_{\mu}(\lam)F(\lam)P^{-1}F(\mu)^*J \quad
\text{where}\  \mu\in\Omega_0\cap\rho({M},{N}),
\end{equation}
up to a constant $J$-unitary factor on the right (that depends upon $\mu$).
\end{thm}

\subsection{Linear fractional transformations}

Let
\begin{equation}\label{eq:2.2}
  T_W[\varepsilon]:=(w_{11}(\lam)\varepsilon(\lam)+w_{12}(\lam))
  (w_{21}(\lam)\varepsilon(\lam)+w_{22}(\lam))^{-1}
\end{equation}
denote the linear fractional transformation of a mvf $\varepsilon\in
\cS^{\ptq}_{\kappa_2}$ $(\kappa_2\in \ZZ_+)$ based on the block
decomposition
\begin{equation}\label{eq:2.1}
W(\lam)=\left[\begin{array}{ll}
 w_{11}(\lam) & w_{12}(\lam) \\
  w_{21}(\lam) & w_{22}(\lam)
\end{array}\right]
\end{equation}
of a mvf $W\in\cU_{\kappa_1}(j_{pq})$ with blocks $w_{11}(\lam)$ and
$w_{22}(\lam)$ of sizes $p\times p$ and $q\times q$, respectively.
The transformation $T_W[\varepsilon]$ is well defined for
those points $\lam\in \gh_W\cap \gh_\varepsilon^+$ for which
$$
\det (w_{21}(\lam)\varepsilon(\lam)+w_{22}(\lam))\not= 0.
$$


\begin{lem}\label{pr:2}
If $W\in \cU_{\kappa_1}(j_{pq})$ and $\varepsilon\in
\cS_{\kappa_2}^{p\times q}$,  then:
\begin{enumerate}
\item[\rm(1)] 
$T_W[\varepsilon]$ admits the supplementary representation
\begin{equation}\label{eq:2.7}
    T_W[\varepsilon]=(w_{11}^\#(\lam)+\varepsilon(\lam)w^\#_{12}(\lam))^{-1}
    (w_{21}^\#(\lam)+\varepsilon(\lam)w^\#_{22}(\lam))\quad\textrm{for}\
    \lam\in\Omega_W\cap \gh_\varepsilon^+\cap \gh_s^+.
\end{equation}
\item[\rm(2)]
$s:=T_W[\varepsilon]\in \cS_{\kappa'}^{p\times q}$ with
$\kappa'\leq\kappa_2+\kappa_1$.
\end{enumerate}
\end{lem}

The set $T_W[\cS^{\ptq}_{\kappa-\kappa_1}]\cap \cS_\kappa^{\ptq}$ is
characterized by three conditions in Theorem 1.1 of of
\cite{DerDym08}, which is repeated immediately below for the
convenience of the reader; it is a generalization of Theorem 3.1 of
\cite{D8}, which treated the case $\kappa_1=\kappa_2=0$; see also
\cite{KKhYu} for more general formulations in a Hilbert space
setting. The interpolation conditions (C1)--(C3) coincide with the
three conditions in the next theorem.



The Krein-Langer representations for mvf's $s\in\cS_\kappa^{\ptq}$
insure that the nontangential limits $s(t)$ exist and are
contractions for almost all points $t\in\Omega_0$. Therefore,
\begin{equation}\label{DeltaS}
\Delta_s(t):=\left[\begin{array}{cc}
I_p & -s(t)\\
-s(t)^* & I_q   \end{array} \right]\ge 0  \quad\textrm{for almost all points}\ \ t\in \Omega_0.
\end{equation}

\begin{thm}
\label{thm:2.2} Let $\kappa,\kappa_1\in\dN\cup\{0\}$,
$(\kappa_1\leq\kappa)$, let $W\in \cU_{\kappa_1}(j_{pq})$ and let
$s\in \cS^{\ptq}_\kappa$ admit the Kre\u{\i}n-Langer factorizations
$s=b_\ell ^{-1}s_\ell =s_rb_r^{-1}$ and let $\Delta_s(\mu)$ and
$\Gamma_\ell$ be defined by~\eqref{DeltaS} and~\eqref{GammaS2}. Then
$s\in T_W[\cS^{\ptq}_{\kappa-\kappa_1}]$ if and only if the
following conditions hold:
\begin{enumerate}
    \item [\rm(1)]
    $\begin{bmatrix}
      b_\ell  & -s_\ell
    \end{bmatrix}f\in H_2^p$ for every $f\in
    {\cK}(W)$;
    \item [\rm(2)]
    $\begin{bmatrix}
      -s_r^* & b_r^*
    \end{bmatrix} f\in (H_2^q)^\perp$ for every $f\in
    {\cK}\rm(W)$;
    \item[\rm(3)]
    $
   {\displaystyle \langle \{\Delta_s
     +\Delta_s
     \left[\begin{array}{cc}
       0 & \Gamma_\ell  \\
       \Gamma_\ell ^* & 0 \\
    \end{array}\right]\Delta_s\}f,f
    \rangle_{nst}\le\langle f,f\rangle_{{\cK}(W)}}
    $ for every $f\in
    {\cK}(W)$.
    \end{enumerate}
\end{thm}

\section{Bitangential interpolation in the unit disc.}
\subsection{Main assumptions.}
In this section the basic theorem  (Theorem \ref{thm:2.2}) will be
used to obtain a linear fractional parametrization  of the set of
solutions $\widehat{\cS}_\kappa(A_1, A_2,C,P,X)$ of the
interpolation problem that is discussed in the Introduction. Note
that if $M$ and $N$ are as in~\eqref{eq:0.MN}, then
\[
\rho({M},{N})=\{\lam\in\dC:\,\det(\lam I_{n_1}-A_1)\ne 0\
\textrm{and}\ \ \det(I_{n_2}-\lam A_2)\ne 0\}.
\]

Let $P_s\in \dC^{n\times n}$ be defined by the
    formula
\begin{equation}\label{eq:3.7}
    v^*P_su=\left\langle \left\{\Delta_s
     +\Delta_s
     \left[\begin{array}{cc}
       0 & \Gamma_\ell  \\
       \Gamma_\ell ^* & 0 \\
    \end{array}\right]\Delta_s\right\}Fu,Fv
    \right\rangle_{nst},\quad\textrm{for}\ u,\, v\in\dC^{n},
    \end{equation}
where  the mvf $\Delta_s$ is given by~\eqref{DeltaS}, $\Gamma_\ell $ is the operator
from ${L_2^q}$ onto $\cH_*(b_\ell )$ defined by~\eqref{GammaS2} and $\Gamma_\ell ^*$ is
the adjoint of $\Gamma_\ell $ with respect to the standard inner product.


\begin{remark} \label{rem:3.1}
The mvf $s\in\widehat{\cS}_\kappa(A_1, A_2, C,P)$ if and only if
$\Delta_s Fu$ belongs to the de Branges-Rovnyak space ${\mathfrak
D}(s)$ and
\begin{equation}\label{eq:3.8}
  \langle\Delta_s Fu,\Delta_s Fu\rangle_{{\sD}(s)}\le u^*Pu,\quad
\text{for every}\
   u\in\dC^n.
\end{equation}
Therefore, since the left hand side of (\ref{eq:3.8}) is equal to $u^*P_su$
and  (C3) implies that $P_s\le P$,
the inequality
\begin{equation}\label{nuP}
    \nu_-(P)\leqslant \kappa
\end{equation}
is necessary for the problem (C1)--(C3) to be solvable in the class
$\cS^{\ptq}_\kappa$. The space $\cD(s)$ is considered in detail in
\cite{DerDym08}; see also \cite{ADRS} and the references cited
therein.

\end{remark}
\begin{remark}\label{Ajordan} Although the condition (B3) is not necessary
for the problem (C1)-(C3) to be solvable, it can be shown that for every data set
$(A_1,A_2,C,P)$ satisfying (B1), (B2) there exists a data set $(\wt A_1,\wt A_2,\wt
C,\wt P)$ satisfying (B1)--(B3),  such that
\[
\cS_\kappa(A_1,A_2,C,P)=\cS_\kappa(\wt A_1,\wt A_2,\wt C,\wt P).
\]
The same statement for the sets of rational mvf's in
$\cS_\kappa(A_1,A_2,C,P)$ and $\cS_\kappa(\wt A_1,\wt A_2,\wt C,\wt
P)$ was proved in~\cite{BGR}. It will be shown below that the
assumptions (B1)-(B4) combined with \eqref{nuP} are sufficient for
the problem (C1)-(C3) to be solvable in the class
$\cS_\kappa^{p\times q}$. The condition (B4) will be discussed in
Subsection~\ref{sub:4.4}.
\end{remark}

\begin{thm}\label{thm:3.1}
If (B1) is in force and $s\in \cS_\kappa^{p\times q}$ satisfies (C1)
and (C2), then $P_s$ is a solution of the Lyapunov-Stein
equation~\eqref{eq:0.LS}.
\end{thm}
\begin{proof}\quad For every $u\in \dC^n$ let us set
\[
h_1{(\lam)}=C({M}-\lam {N})^{-1}M u,\quad h_2{(\lam)}=C({M}-\lam
{N})^{-1}N u
\]
and
\begin{equation}\label{g1g2}
    g_1=P_+\begin{bmatrix}-s^* & I_q\end{bmatrix}h_1,\quad
    g_2=P_-\begin{bmatrix}I_p & -s\end{bmatrix}h_2.
\end{equation}
Since $s$ satisfies the assumptions (C1), (C2), one has
$
g_1\in\cH(b_r),\quad g_2\in\cH_*(b_\ell ).
$
It follows from the identity
\[
({M}-\lam {N})^{-1}{M}=I_n+\lam ({M}-\lam {N})^{-1}{N}
\]
that
\begin{equation}\label{h12Ident}
    h_1(\lam)=Cu+\lam h_2(\lam).
\end{equation}

Using the formulas (\ref{eq:3.7}), 
\eqref{h12Ident} one obtains
\begin{equation}\label{eq:3.12}
\begin{split}
u^*{M}^*P_s{M} u&-u^*{N}^*P_s{N} u=\langle
\Delta_s h_1,h_1\rangle _{nst}-\langle \Delta_s  h_2, h_2\rangle _{nst}\\
    &+\left\langle \left[\begin{array}{cc}
     0 & \Gamma_\ell  \\
        \Gamma_\ell ^* & 0
    \end{array}\right]\Delta_sh_1,\Delta_sh_1\right\rangle _{nst}
    -\left\langle \left[\begin{array}{cc}
     0 & \Gamma_\ell  \\
        \Gamma_\ell ^* & 0
    \end{array}\right]\Delta_sh_2,\Delta_sh_2\right\rangle _{nst}\\
    &=\langle
\Delta_s h_1,h_1\rangle _{nst}-\langle \Delta_s
(h_1-Cu),(h_1-Cu)\rangle _{nst}\\
    &+2{\gR}\left\langle \left[\begin{array}{cc}
     0 & \Gamma_\ell  \\
        0 & 0
    \end{array}\right]\Delta_sh_1,\Delta_sh_1\right\rangle _{nst}
-2{\gR}\left\langle \left[\begin{array}{cc}
     0 & 0\\
        \Gamma_r & 0
    \end{array}\right]\Delta_sh_2,\Delta_sh_2\right\rangle _{nst}\\
    &=2{\gR}\langle
\Delta_s h_1,Cu\rangle _{nst}
+2{\gR}X_1-2{\gR}X_2-\langle \Delta_s Cu,Cu\rangle _{nst},
\end{split}
\end{equation}
where $X_1$ and $X_2$ are given by
\begin{equation}\label{X1}
 \begin{split}
    X_1&=\langle\Gamma_\ell [-s^*\quad I_q] h_1,[I_p\quad -s]h_1\rangle _{nst}\\
    &=\langle\Gamma_\ell [-s^*\quad I_q] h_1,[I_p\quad -s](Cu+\lam h_2)\rangle
    _{nst},
    \end{split}
\end{equation}
\begin{equation}\label{X2}
\begin{split}
    X_2&=\langle \Gamma_r[I_p\quad -s]h_2,[-s^*\quad I_q] h_2\rangle _{nst}\\
    &=\langle\Gamma_r[I_p\quad -s] h_2,\bar\lam[-s^*\quad I_q] (h_1-Cu)\rangle
    _{nst}.
\end{split}
\end{equation}
Decomposing the matrix $C$ as in~\eqref{eq:C}
and using the definitions~\eqref{g1g2} of $g_1$, $g_2$ one obtains
\begin{equation}\label{X11}
    X_1=\langle \Gamma_\ell g_1,\lam g_2\rangle_{nst}-\langle
    s^*\Gamma_\ell g_1,{C}_2u\rangle_{nst},
\end{equation}
\begin{equation}\label{X12}
    X_2=\langle \Gamma_rg_2,\bar\lam g_1\rangle_{nst}+\langle s\Gamma_rg_2,\bar\lam {C}_1u\rangle_{nst}.
\end{equation}
The relations~\eqref{g1g2},~\eqref{X11} and~\eqref{X12} imply
\begin{equation}\label{eq:3.18}
\begin{split}
\langle\Delta_s h_1,Cu\rangle _{nst}+X_1&-X_2=\langle
\begin{bmatrix}I_p & -s\end{bmatrix}h_1,{C}_1u\rangle _{nst}+
\langle \begin{bmatrix}-s^* & I_q\end{bmatrix}h_1,{C}_2u\rangle _{nst}\\&+\langle
\bar\lam\Gamma_\ell g_1, g_2\rangle
    _{nst}-\langle s^*\Gamma_\ell g_1,{C}_2u\rangle
    _{nst}-\langle g_2,\Gamma_\ell \bar\lam g_1\rangle_{nst}-\langle s\Gamma_rg_2,\bar\lam {C}_1u\rangle
    _{nst}\\
    &=\langle
\begin{bmatrix}I_p & -s\end{bmatrix}(Cu+\lam h_2),{C}_1u\rangle _{nst}+
\langle g_1,{C}_2u\rangle _{nst}\\&+\langle \bar\lam\Gamma_\ell g_1, g_2\rangle
    _{nst}-\langle s^*\Gamma_\ell g_1,{C}_2u\rangle
    _{nst}-\langle g_2,\Gamma_\ell \bar\lam g_1\rangle_{nst}-\langle s\Gamma_rg_2,\bar\lam {C}_1u\rangle
    _{nst}\\
&=\langle g_1-s^*\Gamma_\ell g_1,{C}_2u\rangle _{nst}+\langle
g_2-s\Gamma_rg_2,\bar\lam  {C}_1u\rangle _{nst}\\
&+\langle \bar\lam\Gamma_\ell g_1, g_2\rangle
    _{nst}-\langle g_2,\Gamma_\ell \bar\lam g_1\rangle_{nst}+\langle ({C}_1-s{C}_2)u, {C}_1u\rangle
    _{nst}
 \end{split}
\end{equation}
Since $g_1\in\cH(b_r)$, $ g_2\in\cH_*(b_\ell )$ the equations~\eqref{GfGg}
characterizing the vectors $\Gamma_\ell g_1\in\cH_*(b_\ell )$ and $\Gamma_r
g_2\in\cH(b_r)$ can be rewritten as
\begin{equation}\label{Gg1Gg2}
 P_+(g_1-s^*\Gamma_\ell g_1)=0,\quad  P_-(g_2-s\Gamma_rg_2)=0,
\end{equation}
and hence \eqref{eq:3.18} takes the form
\begin{equation}\label{DsX12}
\langle\Delta_s h_1,Cu\rangle _{nst}+X_1-X_2= \langle \bar\lam\Gamma_\ell g_1, g_2\rangle
    _{nst}-\langle g_2,\Gamma_\ell \bar\lam g_1\rangle_{nst}+\langle ({C}_1-s{C}_2)u, {C}_1u\rangle
    _{nst} .
\end{equation}
Substituting \eqref{DsX12} in~\eqref{eq:3.12}  one obtains
\begin{equation}\label{eq:3.20}
\begin{split}
u^*{M}^*P_s{M} u&-u^*{N}^*P_s{N} u=2{\gR}\langle ({C}_1-s{C}_2)u,
{C}_1u\rangle_{nst}-\langle \Delta_s Cu,Cu\rangle _{nst}\\
&+\langle (\bar\lam\Gamma_\ell -\Gamma_\ell \bar\lam)g_1, g_2\rangle_{nst}+\langle
g_2,(\bar\lam\Gamma_\ell -\Gamma_\ell \bar\lam)\lam
    g_1\rangle_{nst}\\
&=\langle {C}_1u, {C}_1u\rangle_{nst}-\langle {C}_2u, {C}_2u\rangle_{nst}+2\mbox{Re
}\langle (\bar\lam\Gamma_\ell -\Gamma_\ell \bar\lam)g_1, g_2\rangle_{nst}.
 \end{split}
\end{equation}
In view of Lemma~\ref{Ker*2} the last term in~\eqref{eq:3.20} is
equal 0 and, hence, \eqref{eq:3.20} can be rewritten as
\[
M^*P_sM-N^*P_sN={C}_1^*{C}_1-{C}_2^*{C}_2.
\]
\end{proof}
Theorem~\ref{thm:3.1} implies that the problem (C1)-(C3) possesses the Parseval
identity property in the sense of~\cite{KhYu94}.

\begin{corollary}\label{cor:3.4}
If (B1) and (B2) are in force and $s\in
\widehat{\cS}_\kappa(A_1,A_2,C,P)$, then
\begin{equation}\label{PsP}
P_s=P,
\end{equation}
and hence
\begin{equation}\label{WSS}
\widehat{\cS}_\kappa(A_1,A_2,C,P)={\cS}_\kappa(A_1,A_2,C,P)\,.
\end{equation}
\end{corollary}
\begin{proof}
 Let $P$
be decomposed conformally with $M$ and $N$:
\begin{equation}\label{eq:3.110}
P =\left[
   \begin{array}{cc}
   P_{11}& P_{12}\\
   P_{21}& P_{22}
   \end{array}\right],\,P_{ij}\in\dC^{n_i\times n_j}\,,
    (i,j=1,2).
\end{equation}
Then, as $P$ is Hermitian,  equation (\ref{eq:0.LS}) is equivalent to the
system of
equations
\begin{equation}\label{eq:3.1011}
    A_1^*P_{11}A_1-P_{11}=C_{11}^*C_{11}-C^*_{21}C_{21},
\end{equation}
\begin{equation}\label{eq:3.1012}
    P_{22}-A_2^*P_{22}A_2=C_{12}^*C_{12}-C^*_{22}C_{22},
\end{equation}
\begin{equation}\label{eq:3.1013}
     P_{21}A_1-A_2^*P_{21}=C_{12}^*C_{11}-C^*_{22}C_{21}.
\end{equation}
Equations \eqref{eq:3.1011} and \eqref{eq:3.1012}
have unique solutions $P_{11}$ and $P_{22}$, respectively, since
$\sigma(A_1)\cup\sigma(A_2^*)\subset\DD$. Therefore, since
$P_s$ is also a solution
    of the Lyapunov-Stein equation~\eqref{eq:0.LS}, it must be of the form
  \begin{equation}
P_s =\left[
   \begin{array}{cc}
   P_{11}& \wt P_{12}\\
   \wt P_{21}& P_{22}
   \end{array}\right],\,\wt P_{ij}\in\dC^{n_i\times n_j},\,\,
    (i,j=1,2).
\end{equation}
Thus, the inequality $P_s\le P$ implies~\eqref{PsP} and so too ~\eqref{WSS}.
\end{proof}
Let $C$ be decomposed as $C=[{\bf C}_1,{\bf C_2}]$ (${\bf C}_j\in
\dC^{m\times n_j},\,j=1,2$) and define the mvf's
\begin{equation}\label{F12}
    F_1(\lam)={\bf C}_1(A_1-\lam I_{n_1})^{-1}\quad\textrm{and}\quad
        F_2(\lam)={\bf C}_2(I_{n_2}-\lam A_2)^{-1}
        \end{equation}
which belong to $\cR\cap (H_2^{m\times n_1})^\perp$ and
$\cR\cap H_2^{m\times n_2}$, respectively. Our first objective is to
check that the condition~\eqref{eq:A4} is satisfied.
\begin{proposition}\label{ChatJ}
If (B1)and (B3) are in force, then (\ref{eq:A4}) holds
and
\begin{equation}
\label{eq:nov4a8}
\ker\left[
      \begin{array}{c}
    M- \lambda N \\ C\\
      \end{array}
    \right]=\{0\}
\quad\text{and}\quad \ker\left[
      \begin{array}{c}
        \lambda M-N \\ C\\
      \end{array}
    \right]=\{0\}
    \quad
\text{for every} \ \lam \in\CC.
\end{equation}
\end{proposition}
\begin{proof}
If $\textup{col}\,(u_1, u_2)
   \in \ker C(M-\lam N)^{-1}$
   for all
   $\lam\in\rho(M,N)$ and some $u_1\in\dC^{n_1}$, $u_2\in\dC^{n_2}$,
then
\begin{equation}\label{ObservIII}
F_1(\lam)u_1+F_2(\lam)u_2\equiv 0,\quad\textrm{for}\ \lam\in\rho(M,N).
\end{equation}
Since  $F_1u_1\in (H_2^m)^\perp$ and $F_2u_2\in
H_2^m$, it follows that
\[
F_1(\lam)u_1={\bf C}_{1}(A_1-\lam I_{n_1})^{-1}u_1\equiv 0
\]
and
\[
F_2(\lam)u_2={\bf C}_{2}(I_{n_2}-\lam A_2)^{-1}u_2\equiv 0 \quad\text{for}\
\lam\in\rho(M,N),
\]
and hence that $u_1=0$, $u_2=0$ in view of (B3).

To complete the proof it suffices to verify the rank condition implicit in
(\ref{eq:nov4a8}) for all points $\lam\in\DD$. But, if $\lam\in\DD$, then,
in view of
(B3),
$$
n\ge\rank\begin{bmatrix}M-\lam N\\C\end{bmatrix} \ge\rank\begin{bmatrix}M
-\lam N\\C_2\end{bmatrix}=n.
$$
This proves the first equality in~\eqref{eq:nov4a8}.
The proof of the remaining assertion in (\ref{eq:nov4a8}) is similar.
\end{proof}

\subsection{Regular case}
We now parametrize the set of solutions
${\cS}_\kappa(A_1,A_2,C,P)$ of the problem (C1)--(C3) assuming that
the matrices $A_1$, $A_2$, $C$, $P$ satisfy the constraints
(B1)--(B3) and  $P$ is invertible.

\begin{thm}\label{thm:3.2}
Let (B1)--(B3) be in force, let $P$ be invertible and let
$\kappa_1:=\nu_-(P)\le\kappa$. Then
\begin{equation}\label{eq:3.25}
{\cS}_\kappa(A_1,A_2,C,P)=\cS^{p\times q}_\kappa\cap
T_W[\cS^{p\times q}_{\kappa-\kappa_1}],
\end{equation}
where the mvf $W(\lam)$ is given by formula (\ref{eq:1.84}).
\end{thm}

\begin{proof}  In view of Example~\ref{Mspace}, the
linear space $\cM=\{F(\lam)u:u\in \CC^n\}$ with inner product (\ref{eq:1.81})
(which is well defined, thanks to Proposition~\ref{ChatJ}) is a
reproducing kernel Pontryagin space with reproducing kernel
${\mathsf K}_\omega(\lam)$ given by (\ref{eq:1.82}) and, in view of
Theorem~\ref{HW}, also by (\ref{eq:0.4}), i.e.,
$$
{\mathsf K}_\omega(\lam)=F(\lam)P^{-1}F(\omega)^*
=\frac{j_{pq}-W(\lam)j_{pq}
   W(\omega)^*}{\rho_\omega(\lam)},
$$
where the mvf $W(\lam)$ is given by~\eqref{eq:1.84}. The rest of the theorem
follows directly from Theorem~\ref{thm:2.2}, since
the conditions (C1)-(C3) for
$s\in S^{p\times q}_\kappa$ are equivalent to the
conditions (1)--(3) of that theorem.
\end{proof}


\subsection{Examples}
We next present a few examples to
illustrate the interpolation problem (C1)--(C3).
\begin{example}\label{ex:6.1} ({\it Interpolation with multiplicities one}).
Let
\begin{equation}\label{eq:6.1}
  C=\left[\begin{array}{lll}
    \xi_1&\cdots & \xi_n \\
    \eta_1&\ \cdots & \eta_n
  \end{array}\right],\ (\xi_j\in\CC^p,\ \eta_j\in \CC^q,\ 1\leq
  j\leq n)
  \end{equation}
  and let the matrices ${M}$ and ${N}$ be as in~\eqref{eq:0.MN} with
$$
A_1=\textup{diag}\, (\alpha_1,\cdots,\alpha_{n_1}),\quad A_2=\textup{diag}\,
(\overline{\alpha}_{n_1+1},\cdots,\overline\alpha_{n})
 $$
 and $\alpha_j\in \DD$, $(j=1,\ldots,n)$.
 Then the interpolation condition (C1) is met if and only if
 \[
 \begin{split}
 \frac{b_\ell (\lam) \xi_j-s_\ell (\lam)\eta_j}{\alpha_j-\lam}&\in H^p_2   \quad
\textrm{for}\quad 1\leq j\leq n_1,
   \\
\textup{and}\\
   \frac{b_\ell {(\lam)} \xi_{j}-s_\ell {(\lam)}\eta_{j}}{1-\lam\overline{\alpha}_j}&\in H^p_2\quad
   \textrm{for}\quad n_1+1\leq j\leq n\,.
 \end{split}
 \]
The second  set of conditions is automatically fulfilled, and the first set
can be rewritten as
 \begin{equation}\label{eq:6.3}
   b_\ell (\alpha_j)\xi_j=s_\ell (\alpha_j)\eta_j\quad \textrm{for}\quad 1\leq j\leq
   n_1\,.
\end{equation}
If the $b_\ell (\alpha_j)$ are invertible, then the constraints in (\ref{eq:6.3}) take
the form
\begin{equation}\label{6.3}
   s(\alpha_j)\eta_j=\xi_j\quad\textrm{for}\ 1\leq j\leq
   n_1.
\end{equation}
Similarly, condition (C2) holds if and only if
$$
\frac{s^\#_r(\lam)\xi_{j}-b^\#_r(\lam)\eta_{j}}{1-\lam\overline{\alpha}_j}\in
(H^q_2)^\perp\quad \textrm{for}\quad n_1+1\leq j\leq n\,,
$$
or, equivalently, if and only if
\begin{equation}\label{eq:6.5}
\xi_j^*s_r(\alpha_j)=\eta_j^*b_r(\alpha_j)
\quad \textrm{for}\quad n_1+1\leq j\leq n.
\end{equation}
If the $b_r(\alpha_j)$ are invertible, then the constraints in (\ref{eq:6.5})
take the form
\begin{equation}\label{6.6}
 \xi_{j}^*s(\alpha_j)=\eta_{j}^*\quad (n_1+1\leq j\leq n).
\end{equation}

If $\sigma(A_1)\cap\sigma(A_2^*)=\emptyset$, then the Lyapunov-Stein equation
(\ref{eq:0.LS}) has exactly one solution. Therefore, $P_s=P$,
since both of these matrices are solutions of (\ref{eq:0.LS}) and hence
the condition (C3) is automatically met.

If $\sigma(A_1)\cap\sigma(A_2^*)\ne\emptyset$, then, since $P_s$ and $P$ are
both Hermitian, (\ref{eq:0.LS}) is equivalent to the system
(\ref{eq:3.1011})-(\ref{eq:3.1013}) for the blocks $P_{11}$, $P_{22}$ and
$P_{21}$ of $P$, respectively. Since the first two equations are uniquely
solvable,
$$
P-P_s=\begin{bmatrix}0&X^*\\X&0\end{bmatrix}\quad\textrm{with}\
X=P_{21}-(P_s)_{21}.
$$
Therefore, (C3) holds if and only if $s\in\cS_\kappa^{\ptq}$ is such that
$P_{21}=(P_s)_{21}$. This imposes an extra interpolation condition on $s$
(see e.g., p.368 of \cite{D7} for an explicit example when $\kappa=0$).
Because of this, (C1)--(C2) was called the basic interpolation
problem and  (C1)--(C3) was called the augmented interpolation problem in
\cite{MSRI}.
\end{example}

\begin{example} ({\it Nevanlinna-Pick matrix interpolation problem}).
Let $n_2=0,\ n_1=tq$, and let the matrices $A_1$ and $C$ have the following
block form
\begin{equation}\label{eq:6.7}
   A_1=\textup{diag}\,(\alpha_1I_q,\ldots,\alpha_t I_q),\,\, 
   C=\left[
   \begin{array}{lll}
     s_1 & \ldots&s_t \\
     I_q &\ldots& I_q
   \end{array}\right]
\end{equation}
where $\alpha_j$ are distinct points in $\DD$ and $s_j\in \CC^{\ptq}$ for
$j=1,\ldots,t$. Then the interpolation condition (C1)  reduces to
\begin{equation}\label{eq:6.8}
   s_\ell (\alpha_j)=b_\ell (\alpha_j)s_j\quad \textrm{for}\ j=1,\ldots,t.
\end{equation}
If $s(\lam)$ is holomorphic at $\alpha_j$ for $j=1,\ldots, t$, then
 these conditions take the form
\begin{equation}\label{eq:6.9}
   s(\alpha_j)=s_j\quad \textrm{for}\ j=1,\ldots,t.
   \end{equation}
 The data $\{A_1,C\}$ specified in
(\ref{eq:6.7}
   satisfies the conditions (B1)--(B3) and the corresponding Lyapunov-Stein
equation (\ref{eq:0.LS}) has a unique solution that may be written in
 block form as
   \begin{equation}\label{eq:6.10}
    P=\left[\frac{I_q-s^*_js_k}{1-
\overline{\alpha}_j\alpha_k}\right]^t_{j,k=1}.
\end{equation}
Let $\kappa_1=\nu_-(P)$. If $P$ is invertible and $\kappa_1\leq \kappa$,
then, by  I of Theorem \ref{thm:0.1}, the set of
solutions to the problem (\ref{eq:6.8}) is described by the formula
$$
{\cS}_\kappa(A_1,C,P)=T_W[\cS_{\kappa-\kappa_1}]\cap
\cS^{\ptq}_{\kappa},
$$
where $W$ is given by (\ref{eq:1.84}).

The Nevanlinna-Pick problem in the class of mvf's $s\in\cS_\kappa^{\ptq}$
that are holomorphic at the interpolation points was investigated by
Golinskii in \cite{Gol83}.
\bigskip
\end{example}

\begin{example}
\label{ex:nov3a8}
Suppose now that
$$
A_1=\begin{bmatrix}
  \alpha &1  & \\
  & \ddots & \ddots\\
  &&& 1\\
 && &\alpha
\end{bmatrix}=\alpha I_{n_1}+N_1 \quad\textrm{and} \quad A_2= \begin{bmatrix}
  \overline{\beta} &1  & \\
  & \ddots & \ddots\\
  &&& 1\\
 && &\overline{\beta}
\end{bmatrix}=\overline{\beta}I_{n_2}+N_2
$$
are Jordan cells of size $n_1\times n_1$ and $n_2\times n_2$,
respectively, and that $\alpha, \beta\in\DD$. Then (C1) holds if and
only if
$$
\begin{bmatrix}b_\ell(\lambda)&-s_\ell(\lambda)\end{bmatrix}\begin{bmatrix}
C_{11}\\C_{21}\end{bmatrix}
(\lambda I_{n_1}-A_1)^{-1}
$$
is holomorphic in $\DD$. A necessary condition for this is that the contour integral
around the unit circle $\dT$ 
$$
\frac{1}{2\pi i}\int_\dT \begin{bmatrix}b_\ell(\zeta)&-s_\ell(\zeta)\end{bmatrix}
\begin{bmatrix}
C_{11}\\C_{21}\end{bmatrix}
(\zeta I_{n_1}-A_1)^{-1}d\zeta =0.
$$
But upon substituting the formula
$$
(\zeta I_{n_1}-A_1)^{-1}
=\sum_{j=0}^{n_1-1}
\frac{N_1^j}{(\zeta-\alpha)^{j+1}}
$$
into the last integral and invoking Cauchy's formula this is readily seen to
reduce to the constraint
$$
\sum_{j=0}^{n_1-1}\frac{b_\ell^{(j)}(\alpha)}{j!}C_{11}N_1^j=
\sum_{j=0}^{n_1-1}\frac{s_\ell^{(j)}(\alpha)}{j!}C_{21}N_1^j,
$$
or, equivalently, in terms of the columns $\xi_1,\ldots,\xi_{n_1}$ of $C_{11}$ and the columns $\eta_1,\ldots,\eta_{n_1}$ of $C_{21}$,
$$
\begin{bmatrix}
a_0&0&\cdots&0\\
a_1&a_0&\cdots&0\\
\vdots & & &\vdots\\
a_{n_1-1}&a_{n_1-2}&\cdots&a_0
\end{bmatrix}
\begin{bmatrix}\xi_1\\ \\ \vdots  \\ \xi_{n_1}\end{bmatrix}=
\begin{bmatrix}
b_0&0&\cdots&0\\
b_1&b_0&\cdots&0\\
\vdots & & &\vdots\\
b_{n_1-1}&b_{n_1-2}&\cdots&b_0
\end{bmatrix}
\begin{bmatrix}\eta_1\\  \\ \vdots  \\ \eta_{n_1}\end{bmatrix},
$$
where
$$
a_j=\frac{b_\ell^{(j)}(\alpha)}{j!}\quad\textrm{and}\quad
b_j=\frac{s_\ell^{(j)}(\alpha)}{j!}.
$$
It is readily checked that this condition is sufficient as well as necessary,
and that it is equivalent to the asymptotic condition
\begin{equation}
\label{eq:nov3a8}
\begin{split}
{}&b_\ell(\lam)(\xi_1+(\lam-\alpha)\xi_2+\cdots+
(\lam-\alpha)^{n_1-1}\xi_{n_1})\\ &=
s_\ell(\lam)(\eta_1+(\lam-\alpha)\eta_2+\cdots+(\lam-\alpha)^{n_1-1}\eta_{n_1})
+O(\lam-\alpha)^{n_1}\quad\textrm{as}\ \lam\rightarrow \alpha.
\end{split}
\end{equation}

Similarly, the condition (C2) is met if and only if the columns of the mvf
$$
\begin{bmatrix}-s_r^\# & b_r^\#\end{bmatrix}\begin{bmatrix}C_{21}
\\C_{22}\end{bmatrix}(I_{n_2}-\lambda A_2)^{-1}
$$
belong to $(H_2^q)^\perp$, or, equivalently, if and only if
$$
(\lambda I-A_2^*)^{-1}\begin{bmatrix}C_{21}^*& C_{22}^*\end{bmatrix}
\begin{bmatrix}-s_r \\  b_r\end{bmatrix}
$$
is holomorphic in $\DD$, i.e., if and only if
$$
\sum_{j=0}^{n_2-1}(N_2^*)^jC_{21}^*\frac{s_r^{(j)}(\beta)}{j!}=
\sum_{j=0}^{n_2-1}(N_2^*)^jC_{22}^*\frac{b_r^{(j)}(\beta)}{j!}.
$$
The last constraint can be rewritten in terms of the columns
$\xi_{n_1+1},\ldots,\xi_n$ of $C_{21}$ and the columns
$\eta_{n_1+1},\ldots,\eta_n$ of $C_{22}$ as
$$
\begin{bmatrix}\xi_{n_1+1}^*&\cdots&\xi_n^*\end{bmatrix}
\begin{bmatrix}c_0&c_1&\cdots&c_{n_2-1}\\0&c_0&\cdots&c_{n_2-2}\\
\vdots& & &\vdots\\
0&0&\cdots&c_0\end{bmatrix}=
\begin{bmatrix}\eta_{n_1+1}^*&\cdots&\eta_n^*\end{bmatrix}
\begin{bmatrix}d_0&d_1&\cdots&d_{n_2-1}\\0&d_0&\cdots&d_{n_2-2}\\
\vdots& & &\vdots\\
0&0&\cdots&d_0\end{bmatrix},
$$
where
$$
c_j=\frac{s_r^{(j)}(\beta)}{j!}\quad\textrm{and}\quad
d_j=\frac{b_r^{(j)}(\beta)}{j!}\quad\textrm{for}\ j=0,\dots,n_2-1.
$$
This condition is equivalent to
\begin{equation}
\label{eq:nov3b8}
\begin{split}
{}&(\xi_{n_1+1}^*+(\lam-\beta)\xi_{n_1+2}^*+\cdots+(\lam-\beta)^{n_2-1}\xi_n^*)
s_r(\lam)\\ &=
(\eta_{n_1+1}^*+(\lam-\beta)\eta_{n_1+2}^*+\cdots+(\lam-\beta)^{n_2-1}
\eta_n^*)b_r(\lam)+O((\lam-\beta)^{n_2})\quad\textrm{as}\ \lam\rightarrow
\beta.
\end{split}
\end{equation}
\end{example}

Interpolation problems in which  $A_1$ and/or $A_2$ are made up of several
Jordan blocks lead to similar sets of formulas.

%
If $s$ is holomorphic on
$\sigma(A_1)\cup\overline{\sigma(A_2)}$ then the interpolation
conditions~\eqref{eq:nov3a8} and
\eqref{eq:nov3b8} take the form
\begin{equation}\label{eq:6.13}
\begin{array}{l}
s(\lam)(\eta_{1}+\eta_2(\lam-\alpha)+\cdots+
(\lam-\alpha)^{n_1-1}\eta_{n_1})\\
=\xi_{1}+(\lam-\alpha)\xi_2+\cdots+(\lam-\alpha)^{n_1-1}\xi_{n_1})
+O((\lam-\alpha)^{n_1})
   \end{array}
\end{equation}
\begin{equation}\label{eq:6.130a}
\begin{array}{l} (\xi_{n_1+1}^*+(\lam-\beta)\xi_{n_1+2}^*
+(\lam-\beta)^{n_2-1})\cdots+\xi_n^*
s(\lam)\\
 =  \eta_{n_1+1}^*+(\lam-\beta)\eta_{n_1+2}^*+\cdots
+(\lam-\beta)^{n_2-1}\eta_{n}^* +O((\lam-\beta)^{n_2})
    \end{array}
\end{equation}

One sided versions of the problem considered in Example \ref{ex:nov3a8}
($\ell_2=0$) of the problem
have been considered by T.~Takagi~\cite{T24} and
A. A.~Nudelman~\cite{Nud81}, respectively.
Rational solutions of the two-sided Takagi-Nudelman
problem~\eqref{eq:6.13}--\eqref{eq:6.130a} were described in~\cite{BGR}.

\begin{remark}
\label{rem:jun5a9}
If (C4) is in force, i.e., if $b_\ell(\lambda)$ and $b_r(\lambda)$ are
invertible at every point $\lam\in\sigma(A_1)\cup\overline{\sigma(A_2)}$,
then the conditions (C1)--(C3) are equivalent to the residue conditions
$$
\sum_{j=1}^{n_1}\textup{res}\{s(\lambda)C_{21}(\lambda I_{n_1}-A_1)^{-1}\}
=C_{11},
$$
$$
\sum_{j=n_1+1}^{n_1+n_2}\textup{res}\{(\lambda I_{n_2}-A_2^*)^{-1}C_{12}^*
s(\lambda)\}=C_{22}^*
$$
and
$$
\sum_{j=1}^{n_1+n_2}\textup{res}\{(\lambda I_{n_2}-A_2^*)^{-1}C_{12}^*
s(\lambda)C_{21}(\lambda I_{n_1}-A_1)^{-1}\}=P_{21}
$$
in \cite{BGR}, respectively.
\end{remark}


\subsection{Resolvent matrix in the singular case.\label{sub:4.4}}
If $P$ is not invertible, then the construction of the model space
$\cM$ depends upon assumption (B4); cf.~\cite{D8}, \cite{BD}. Since $\ran X$
is presumed to be invariant under $M$ and $N$, there exist matrices
$M_0\in\CC^{\ntn}$
and $N_0\in\CC^{\ntn}$ such that
\begin{equation}
\label{eq:dec16a8}
MX=XM_0\quad \text{and}\quad NX=XN_0.
\end{equation}
Moreover, the matrices $M_0$ and $N_0$ may be chosen so that
\begin{equation}
\label{eq:dec16b8}
\sigma(M_0)\subseteq\sigma(M)\quad\text{and}\quad
\sigma(N_0)\subseteq\sigma(N);
\end{equation}
see e.g., \cite{linalg}.

\begin{lem}\label{Odecom}
If (B4) is in force and $M$ and $N$ are as in \eqref{eq:0.MN}, then
\begin{equation}\label{eq:XP}
\nu_-(X)=\nu_-(P),\quad \textup{rank }X=\textup{rank }P
\end{equation}
and
\begin{equation}\label{ranXdecom}
    \ran X=\left(\ran X\bigcap
\begin{bmatrix}
 \dC^{n_1} \\
  0 \\
\end{bmatrix}
\right) \dot{+} \left(\ran X\bigcap
\begin{bmatrix}
  0 \\
   \dC^{n_2} \\
\end{bmatrix}
\right).
\end{equation}
\end{lem}

\begin{proof}
The equality $XPX=X$ implies that $\nu_-(X)\le\nu_-(P)$ and
$\mbox{rank }X\le\mbox{rank }P$; whereas  the equality $PXP=P$ implies
the opposite inequalities. Therefore, (\ref{eq:XP}) holds.

Next, since
$$
\CC^n=\begin{bmatrix}
 \dC^{n_1} \\
  0 \\
\end{bmatrix}
 \dot{+}\begin{bmatrix}
  0 \\
   \dC^{n_2} \\
\end{bmatrix},
$$
to verify (\ref{ranXdecom}), it suffices to show that
$$
\text{if}\ u\in\textup{rng} X,\quad\text{then}\quad
\begin{bmatrix}0&0\\0&I_{n_2}\end{bmatrix}u\in\textup{rng} X.
$$
The invariance assumption implies that $MX=XM_0$ for some $M_0\in\CC^{\ntn}$
with $\sigma(M_0)\subseteq\sigma(M)$. Therefore, $M^kX=XM_0^k$, and since
$\sigma(A_1)\subset\DD$, the entries in  $M_0^k$ are bounded and so a
subsequence $M_0^{k_j}$ tends to a limit $H\in\CC^{\ntn}$ as $k_j$ tends to
infinity. Thus,  if $u=Xv$, then
$$
\begin{bmatrix}0&0\\0&I_{n_2}\end{bmatrix}u=\lim_{j\uparrow\infty}M^{k_j}u=
\lim_{j\uparrow\infty}XM_0^{k_j}v=XHv\,.
$$
\end{proof}

\begin{lem}\label{lem:3.5}
If (B1), (B2) and (B4) are in force, then $X$ is a solution of the
Riccati equation
\begin{equation}
\label{eq:dec16c8}
XM^*PMX-XN^*PNX=XC^*j_{pq}CX,
\end{equation}
or, equivalently (in view of (\ref{eq:dec16a8})),
\begin{equation}\label{eq:3.47}
   M_0^*XM_0-N^*_0XN_0=XC^*j_{pq}CX.
\end{equation}
\end{lem}
\begin{proof} This follows from (\ref{eq:0.LS}) upon multiplying it on the
left and on the right by $X$ and using (\ref{eq:dec16a8}).
\end{proof}
\bigskip

 Let $\cM$ be the space of rational vvf's
 \begin{equation}\label{eq:3.48}
   \cM=\{F(\lam)Xu:u\in \CC^n\}
\end{equation}
endowed with the inner product
\begin{equation}\label{eq:3.49}
    \left\langle FXu,\ FXv\right\rangle_{\cM}=v^*Xu.
\end{equation}
If the condition (B3)
holds, then the inner product in $\cM$ is well defined, since by
Proposition~\ref{ChatJ} the identity $F(\lam)Xu\equiv 0$
implies $Xu=0$ and hence $\left\langle FXu,\ FXv\right\rangle_{\cM}=0$ for all
$v\in \CC^n$. $\cM$ is a reproducing kernel space with kernel
\begin{equation}\label{eq:3.50}
    {\mathsf K}_\omega(\lam)=F(\lam)XF(\omega)^*\quad\textrm{for}\ \lam,\omega\in\rho(M,N).
\end{equation}

\begin{lem}\label{lem:3.6}
Let (B1)--(B4)  be in force, let  the reproducing kernel space $\cM$
be given by (\ref{eq:3.48}), (\ref{eq:3.49}) and let the mvf
$W(\lambda)$ be given by~\eqref{eq:0.5}
Then $\cM$  is a finite dimensional de
Branges-Krein
space ${\cK}(W)$ with reproducing kernel
\begin{equation}\label{eq:3.52}
   {\mathsf K}_\omega(\lam)=F(\lam)XF(\omega)^*={\mathsf K}^W_\omega(\lam):
   =\frac{j_{pq}-W(\lam)j_{pq}W(\omega)^*}{\rho_{\omega}(\lam)}
\end{equation}
and negative index $\nu_-(P)$.
\end{lem}

\begin{proof} Upon substituting (\ref{eq:0.5}) into (\ref{eq:0.4}), direct
calculations show that
\begin{equation}\label{eq:3.53}
  {\mathsf K}^W_\omega(\lam)=F(\lam)X(M^*-\overline{\mu}N^*)^{-1}
\frac{Y}{\rho_\omega(\lam)}
  (M-\mu N)^{-1} XF(\omega)^*,
\end{equation}
where, with the help of the Lyapunov-Stein equation (\ref{eq:0.LS}) and the
formulas
$$
X(M^*-\overline{\mu}N^*)^{-1}=X(M^*-\overline{\mu}N^*)^{-1}PX\quad\textrm{and}
\quad (M-\mu N)^{-1} X=XP(M-\mu N)^{-1} X
$$
based on the invariance of the $\ran X$ under $M$ and $N$, the central term
\begin{eqnarray}\label{eq:3.54}
  Y&=&\rho_\mu(\lam)(M^*-\overline{\omega}N^*)P(M-\mu
  N)\nonumber\\&+& \rho_\omega(\mu)(M^*-\overline{\mu}N^*)
  P(M-\lam  N)- \rho_\mu(\lam)\rho_\omega(\mu)C^*j_{pq}C\\
&=&\rho_\omega(\lam)(M^*-\overline{\mu}N^*)P(M-\mu N). \nonumber
\end{eqnarray}
The verification of the formula for $Y$ rests on the evaluations
\begin{eqnarray*}
F(\lambda)XF(\mu)^*&=&F(\lambda)X(M^*-\overline{\mu}N^*)^{-1}(M^*
-\overline{\omega}N^*)F(\omega)^*\\
&=&F(\lambda)X(M^*-\overline{\mu}N^*)^{-1}(M^*
-\overline{\omega}N^*)PXF(\omega)^*\\
&=&F(\lambda)X(M^*-\overline{\mu}N^*)^{-1}(M^*
-\overline{\omega}N^*)P(M-\mu N)(M-\mu N)^{-1}XF(\omega)^*
\end{eqnarray*}
and
$$
F(\mu)XF(\omega)^*=F(\lambda)X(M^*-\overline{\mu}N^*)^{-1}
(M^*-\overline{\mu}N^*)P(M-\lambda N)(M-\mu N)^{-1}XF(\omega)^*.
$$
Now (\ref{eq:3.53}) and (\ref{eq:3.54}) yield (\ref{eq:3.52}). \end{proof}
\begin{lem}\label{lem:3.7}
If the Hermitian matrices $X,\,P\in \CC^{n\times n}$ satisfy $XPX=X$ and
$\textup{rank}\,X=\textup{rank}\,P$, then
\begin{equation}\label{eq:3.29}
\CC^n=\textup{\ ker} P\dot{+}\ran\ X
\end{equation}
\end{lem}
\begin{proof}
(Cf. \cite{D8}). If $u=Xv$ and $Pu=0$, then $0=XPXv=Xv=u$. Thus,  the
indicated sum is direct. The rest follows from the fact that
\[
n=\dim\ran P+\dim\ker P,
\]
and, by assumption, $\dim\ran P=\dim\ran X$.
\end{proof}

The assumptions of Lemma~\ref{lem:3.7} are
fulfilled if $X$ and $P$ satisfy (i) and (ii) in (B4).
\begin{lem}\label{lem:6.1}
If $\lam\in\rho(M,N)\cap \rho (N^*,M^*)$, $W$ is as in
\eqref{eq:0.5} and
\begin{equation}\label{eq:9.3}
    \wt C=F(\mu)X(\mu M^*-N^*),
\end{equation}
then
\begin{equation}\label{eq:9.1}
    j_{pq}W^\#(\lam)j_{pq}F(\lam)X=\wt C(\lam M^*-N^*)^{-1}
\end{equation}
and
\begin{equation}\label{eq:9.2}
    j_{pq}W^\#(\lam)j_{pq}F(\lam)v=F(\mu)(I_n-X P)(M-\mu N)(M-\lam N)^{-1}v
\quad\text{for} \ v\in\ker P.
\end{equation}
\end{lem}

\begin{proof}
Formula \eqref{eq:0.5} implies that
\[
\begin{split}
j_{pq}W^\#&(\lam)j_{pq}F(\lam)=C[I_n-(\lam-\mu)(M-\mu N)^{-1}X(\lam
M^*-N^*)^{-1}C^*j_{pq}C](M-\lam N)^{-1}\\
&=C(M-\mu N)^{-1}[M-\mu N-(\lam-\mu)X(\lam
M^*-N^*)^{-1}C^*j_{pq}C](M-\lam N)^{-1}.
\end{split}
\]
Therefore, since
\[
C^*j_{pq}C=M^*PM-N^*PN=M^*P(M-\lam N)+(\lam M^*-N^*)PN,
\]
it follows that
\begin{equation}\label{eq:9.4}
\begin{split}
j_{pq}W^\#(\lam)j_{pq}F(\lam) &=F(\mu)[-(\lam-\mu)X(\lam
M^*-N^*)^{-1}M^*P]\\
&+F(\mu)[M-\mu N-(\lam-\mu)XPN](M-\lam N)^{-1} \\
=&F(\mu)[I_n-(\lam-\mu)X(\lam
M^*-N^*)^{-1}M^*P]\\
&+(\lam-\mu)F(\mu)(I_n-XP)N(M-\lam N)^{-1}.
\end{split}
\end{equation}
Thus, as
$$
(I_n-XP)N(M-\lambda N)^{-1}X=0\quad\textrm{and}\quad X(\lam M^*-N^*)^{-1}M^*(I_n-PX)=0,
$$
it is readily seen that
$$
X(\lambda M^*-N^*)^{-1}M^*PX
=X(\lambda M^*-N^*)^{-1}M^*,
$$
and hence that formula (\ref{eq:9.4}) implies that
\[
\begin{split}
j_{pq}W^\#(\lam)j_{pq}F(\lam)X &=F(\mu)X[I_n-(\lam-\mu)(\lam
M^*-N^*)^{-1}M^*]\\
&=F(\mu)X[\mu M^*- N^*](\lam M^*-N^*)^{-1}.
\end{split}
\]
This proves~\eqref{eq:9.1}.

Next, if $v\in\ker P$, then formula (\ref{eq:9.4}) implies that
\begin{equation}\label{eq:9.5}
\begin{split}
j_{pq}W^\#(\lam)j_{pq}F(\lam)v &=F(\mu)[I_n+(\lam-\mu)
(I_n-XP)N(M-\lam N)^{-1}]v \\
&=F(\mu)[I_n+(I_n-XP)\{(\lam N-M)+(M-\mu N)\}(M-\lam
N)^{-1}]v\\
&=F(\mu)XPv+ F(\mu)(I_n-XP)(M-\mu N)(M-\lam N)^{-1}v.
\end{split}
\end{equation}
This proves~\eqref{eq:9.2} since $F(\mu)XPv=0$.
\end{proof}

\begin{remark}
Lemma~\ref{lem:6.1} is similar to  \cite[Theorem 5.1 ]{Dym01} (see
also~\cite[Lemma 2.1]{Bol04}). If $P$ is invertible, then $P^{-1}$ satisfies
the  Lyapunov-Stein equation
\[
MP^{-1}M^*-NP^{-1}N^*=\wt C^*j_{pq}\wt C.
\]
\end{remark}
\begin{lem}\label{lem:7.2}
Let $W$ and $\wt{C}$ be given by~\eqref{eq:0.5} and~\eqref{eq:9.3},
respectively. Then
\begin{equation}\label{eq:7.09}
  W(\lambda)=D+C\begin{bmatrix}
   (A_1-\lambda I_{n_1})^{-1} & 0 \\
    0 & \lambda(I_{n_2}-\lambda A_2)^{-1} \
  \end{bmatrix}\widetilde{C}^* j_{pq},\quad\lam\in\rho(M,N)
\end{equation}
and
\begin{equation}\label{eq:7.09A}
\begin{split}
  W^{-1}(\lambda)
  &=j_{pq}D^* j_{pq}-\wt{C}\begin{bmatrix}
    \lambda(I_{n_1}-\lambda A_1^*)^{-1} & 0 \\
    0 & (A_2^*-\lambda I_{n_2})^{-1}\
  \end{bmatrix}C^* j_{pq},\quad\lam\in\rho(N^*,M^*),
\end{split}
\end{equation}
where
\[
{D}=I_m-C\begin{bmatrix}
  \bar{\mu}I_{n_1} & 0 \\
  0 & I_{n_2}
\end{bmatrix}XF(\mu)^*j_{pq}.
\]
\end{lem}
\begin{proof}
This is a tedious but straightforward computation based on the identities
\[
(1-\bar{\mu}\lambda)(A_1-\lambda I_{n_1})^{-1}=\bar{\mu} I_{n_1}+(A_1-
\lambda I_{n_1})^{-1} (I_{n_1}-\bar{\mu}A_1)\,,
\]
\[
(1-\bar{\mu}\lambda)(I_{n_2}-\lambda A_2)^{-1}=I_{n_2}+\lambda
(I_{n_2}-\lambda A_2)^{-1} (A_2-\bar{\mu}I_{n_2})
\]
and
$
 W^{-1}(\lambda)=j_{pq}W^\#(\lambda)j_{pq}.
$
\end{proof}
\begin{lem}\label{lem:7.1}
Let $(A_1, A_2, C, P)$ satisfy $(B1)-(B4)$, let $X$ be defined in
$(B4)$, and let
\begin{equation}\label{eq:7.1}
\wt{\bC}_1=\wt{C}\begin{bmatrix}I_{n_1}\\0\end{bmatrix}\quad\text{and}\quad
\wt{\bC}_2= \wt{C}\begin{bmatrix}0\\I_{n_2}\end{bmatrix}.
\end{equation}
Then
\begin{equation}\label{eq:7.4}
\cL_1:=\ker\begin{bmatrix}
  \wtilde{\bC}_1\\
  \wtilde{\bC}_1 A_1^* \\
    \vdots\\
  \wtilde{\bC}_1 {A_1^*}^{n_1 -1}
\end{bmatrix}=\ker\begin{bmatrix}
   X_{11} \\
   X_{21}
\end{bmatrix}\quad\textrm{and}\quad
\cL_2:=\ker\begin{bmatrix}
  \wtilde{\bC}_2 \\
  \wtilde{\bC}_2 A_2^* \\
    \vdots\\
  \wtilde{\bC}_2 {A_2^*}^{n_2 -1}
\end{bmatrix}=\ker\begin{bmatrix}
   X_{12} \\
   X_{22}
\end{bmatrix}\,.
\end{equation}
\end{lem}

\begin{proof}
Since the subspace $\cL_1$ is invariant under $A_1 ^*$, it can be
decomposed into the sum of algebraic subspaces of $A_1
^*|_{\cL_1}$. Let $u_1, u_2, \dots, u_k$ be a maximal chain of
generalized eigenvectors of $A_1 ^*|_{\cL}$ corresponding to an
eigenvalue $\alpha$ and let $U=\begin{bmatrix}u_1& \cdots &
u_k\end{bmatrix}$. Then
\begin{equation}\label{eq:7.51}
\wt{\bC}_1(A_1^*)^jU=0 \quad\text{for}\ j=0, 1,\ldots \quad
\textrm{and}\quad A_1^*U=U(\alpha I_k+Z),
\end{equation}
 where $Z$ is an upper triangular  $k\times k$ Jordan cell with $0$ on the diagonal.
Therefore, since $\sigma(A_1)\subset\DD$,
$$
-C(M-\mu N)^{-1}\begin{bmatrix}X_{11}\\X_{21}\end{bmatrix}U=\wt{\bC}_1
(I_{n_1}-\mu A_1^*)^{-1}U
=\sum_{j=0}^\infty\wt{\bC}_1\mu^j(A_1^*)^jU=0.
$$
Thus, if
$$
\wt{u}_j=\begin{bmatrix}
  u_j \\
  0
\end{bmatrix}
\quad\textrm{for $j=1, \ldots, k$ and correspondingly}\quad
\wt{U}=\begin{bmatrix}U\\0\end{bmatrix}\,,
$$
then
\begin{equation}
\label{eq:dec20a7}
C(M-\mu N)^{-1}X\wt{U}=0\,.
\end{equation}
The Lyapunov-Stein equation implies that
\begin{equation}\label{eq:7.6}
  M^* P(M- \mu N)+(\mu M^*-N^*)PN=C^* j_{pq} C
\end{equation}
and hence that
\begin{equation}\label{eq:7.7}
\begin{split}
X(\mu M^*-N^*)^{-1}M^* PX &+X P N(M-\mu N)^{-1}X \\
&=X(\mu M^*-N^*)^{-1}C^* j_{pq} C(M-\mu N)^{-1}X\,.
\end{split}
\end{equation}
In view of the presumed invariance of the range of $X$ under
multiplication by $M$ and $N$,
\[
X(\mu M^*-N^*)^{-1}M^* P X=X(\mu M^*-N^*)^{-1}M^*
\]
and
\[
X P N(M-\mu N)^{-1}X=N(M-\mu N)^{-1}X.
\]
Therefore, equations \eqref{eq:7.7} and (\ref{eq:dec20a7}) imply that
\begin{equation}\label{eq:7.10}
X(\mu M^*-N^*)^{-1}M^*\wt{U} +N(M-\mu N)^{-1}X \wt{U}=0.
\end{equation}
However, since
\begin{eqnarray*}
X(\mu M^*-N^*)^{-1}M^*\wt{U}&=&\begin{bmatrix}X_{11}\\X_{21}\end{bmatrix}
(\mu A_1^*-I_{n_1})^{-1}A_1^*U\\
&=&\begin{bmatrix}X_{11}\\X_{21}\end{bmatrix}U(\beta I_k+\mu Z)^{-1}
(\alpha I_k+\mu Z),
\end{eqnarray*}
where $\beta=\mu\alpha -1\ne 0$, the matrix
$$
E=\begin{bmatrix}X_{11}\\X_{12}\end{bmatrix}U
$$
is a solution of the equation
$$
E(\beta I_k+\mu Z)^{-1}(\alpha I_k+Z)+N(M-\mu N)^{-1}E=0.
$$
Moreover, since
$$
(\beta I_k+\mu Z)^{-1}(\alpha I_k+Z)=\gamma I_k+(1-\mu\gamma)
(\beta I_k+\mu Z)^{-1}Z
\quad\text{with} \quad \gamma=\alpha/\beta
$$
and
$$
\gamma I_n+N(M-\mu N)^{-1}=\beta^{-1}(\alpha M-N)(M-\mu N)^{-1},
$$
it is readily seen that
$$
(\alpha M-N)(M-\mu N)^{-1}E=(\mu\gamma-1)\beta E(\beta I_k+\mu Z)^{-1}Z.
$$
Let $e_j$, $j=1,\ldots,k$, denote the $j$'th column of $I_k$. Then, since
$$
C(M-\mu N)^{-1}E=0,\quad Ze_1=0\quad \textrm{and}\quad Ze_{j+1}=e_j \quad
\textrm{for}\ j=1,\ldots,k-1,
$$
it follows that
$$
\begin{bmatrix}\alpha M-N\\C\end{bmatrix}(M-\mu N)^{-1}Ee_1=0
$$
and hence that $Ee_1=0$. Next, proceeding inductively, suppose that
$Ee_1=\cdots=Ee_j=0$. Then
\begin{eqnarray*}
(\alpha M-N)(M-\mu N)^{-1}Ee_{j+1}&=&
(\mu\gamma-1)\beta E(\beta I_k+\mu Z)^{-1}Ze_{j+1}\\
&=&
(\mu\gamma-1)\beta E(\beta I_k+\mu Z)^{-1}e_j\\
&=&
(\mu\gamma-1)\{Ee_j-\mu E(\beta I_k+\mu Z)^{-1}Ze_j\\
&=& (\mu/\beta)(\alpha M-N)(M-\mu N)^{-1}Ee_j=0.
\end{eqnarray*}
Therefore,
$$
\begin{bmatrix}\alpha M-N\\C\end{bmatrix}(M-\mu N)^{-1}Ee_{j+1}=0
$$
too.
This proves that the columns of $U$ are in the kernel of
$\textup{col}(X_{11}, X_{21})$. The opposite inclusion follows from the
presumed invariance:
\[
\wt{\bC}_1(A_1^*)^j=F(\mu)X(\mu M^*-N^*)^{-1}(M^*)^j
\begin{bmatrix}I_{n_1}\\0\end{bmatrix}
=F(\mu)X(\mu M^*-N^*)^{-1}(M^*)^jPX
\begin{bmatrix}I_{n_1}\\0\end{bmatrix}.
\]
This completes the proof of the first assertion. The proof of the
second is similar.
\end{proof}

\begin{corollary}\label{cor:7.1}
Let the assumptions of Lemma~\ref{lem:7.1} hold. Then
\begin{equation}\label{eq:7.2A}
\ran\begin{bmatrix}
  \wt{\bC}_1^* & A_1\wt{\bC}_1^* & \dots & A_1^{n_1-1}\wt{\bC}_1^*
  \end{bmatrix}
  =\ran\begin{bmatrix}
    X_{11} & X_{12} \
  \end{bmatrix},
\end{equation}
\begin{equation}\label{eq:7.3A}
\ran\begin{bmatrix}
  \wt{\bC}_2^* & A_2\wt{\bC}_2^* & \dots & A_2^{n_2-1}\wt{\bC}_2^*
  \end{bmatrix}
  =\ran\begin{bmatrix}
    X_{21} & X_{22} \
  \end{bmatrix}.
\end{equation}
\end{corollary}
\subsection{Interpolation with singular $P$}
 Given a set of matrices $(A_1, A_2, C, P, X)$ that satisfy (B1)--(B4),
let $\widehat{\cS}_\kappa(A_1, A_2,C,P,X)$ denote the set of mvf's
$s\in \cS^{p\times q}_\kappa$ which satisfy the following
conditions:
\begin{enumerate}
    \item [(X1)]$\begin{bmatrix}
      b_\ell  & -s_\ell  \
    \end{bmatrix}FXu\in H^p_2$\ for every $u\in \CC^n$.
    \vspace{2mm}
    \item [(X2)]$\begin{bmatrix}
      -s^\#_r & b^\#_r \
    \end{bmatrix}FXu\in (H^q_2)^\perp$\ for every $u\in \CC^n$.
    \vspace{2mm}
    \item [(X3)]$XP_sX\leq XPX\,(=X)$.
\end{enumerate}

\begin{thm}\label{thm:3.9}
If (B1)--(B4) are in force and $\kappa\ge \kappa_1=\nu_-(P)$, then
$$\widehat{\cS}_\kappa(A_1, A_2,C,P,X)=
T_W[\cS_{\kappa-\kappa_1}^{\ptq}]\cap \cS^{p\times q}_\kappa,
$$
where the mvf $W$ is given by (\ref{eq:0.5}).
\end{thm}
\begin{proof}
1) Consider the linear space $\cM$ defined by (\ref{eq:3.48}) with
 inner product (\ref{eq:3.49}). By Lemma~\ref{lem:3.6} $\cM$
is a finite dimensional de Branges-Krein space ${\cK}(W)$ with $W$
given by (\ref{eq:0.5}).

Let $s\in \widehat{\cS}_\kappa(A_1, A_2,C,P,X)$. Then it follows
from (X3) that $XP_sX\le X$ and for every $u\in \CC^n$
\begin{equation}\label{eq:3.103}
\left\langle\left\{ \Delta_s
     +\Delta_s
     \left[\begin{array}{cc}
       0 & \Gamma_\ell  \\
       \Gamma_\ell ^* & 0 \\
    \end{array}\right]\Delta_s\right\}FXu,FXu\right
    \rangle_{nst}\le\langle FXu,FXu\rangle_{{\cK}(W)}.
   \end{equation}
Due to Theorem~\ref{thm:2.2} this implies that $s\in
T_W[\cS^{p\times q}_{\kappa-\kappa_1}]$.
 This
proves the inclusion
\begin{equation}\label{eq:3.104}
\widehat{\cS}_\kappa(A_1, A_2,C,P,X)\subseteq \cS^{p\times
q}_{\kappa}\cap T_W[\cS^{p\times q}_{\kappa-\kappa_1}].
\end{equation}

2) Conversely, let  $s\in \cS^{p\times q}_{\kappa}\cap
T_W[\cS^{p\times q}_{\kappa-\kappa_1}]$. Then, by
Theorem~\ref{thm:2.2} the conditions (X1), (X2) and the
inequality~\eqref{eq:3.103} hold. In view of~\eqref{eq:3.7}
and~\eqref{eq:3.49} this implies $XP_sX\le X$.
 Therefore $s\in
    \widehat{\cS}_\kappa(A_1, A_2,C,P,X)$ which serves to prove the opposite
    inclusion in (\ref{eq:3.104}).
\end{proof}

Clearly,
$$
{\cS}_\kappa(A_1, A_2,C,P)\subseteq \widehat{\cS}_\kappa(A_1,
A_2,C,P,X),
$$
since
$$
s\in {\cS}_\kappa(A_1, A_2,C,P)\Longrightarrow P_s= P\Longrightarrow
XP_sX= XPX=X.
$$
In fact $XP_sX=X$ for every solution $s\in\widehat{\cS}_\kappa(A_1,
A_2,C,P,X)$, but we will not verify that in this paper.  Every mvf
$s\in {\cS}_\kappa(A_1, A_2,C,P)$ satisfies also an extra condition
that will be presented below in Lemma~\ref{lem:3.90}. It is
convenient, however, to first establish a preliminary fact from
linear algebra.

\begin{lem}
\label{lem:may31a9}
If $A\in\CC^{\ntn}$ is positive semidefinite, $B\in\CC^{\ntn}$ is Hermitian
and $\nu_-(A+B)=\nu_-(B)$, then $\textup{ker}(A+B)\subseteq\textup{ker}\,A$.
\end{lem}

\begin{proof} Let $k=\nu_-(A+B)=\nu_-(B)$. Since the conclusion of the lemma
is clear if $k=0$, suppose that $k>0$ and let $\lambda_1\le\cdots\le\lambda_n$
denote the eigenvalues of $A+B$, let $x_1,\ldots,x_n$ be a corresponding
set of eigenvectors and suppose further that $\textup{ker}(A+B)\ne\emptyset$.
Then, $\lambda_1\le\cdots\le\lambda_k<0$,
$\lambda_{k+1}=0$.  and  $\langle Ax_{k+1}, x_{k+1}\rangle\ge 0$. However, if
$\langle Ax_{k+1}, x_{k+1}\rangle> 0$, then
$$
\langle Bx_{k+1}, x_{k+1}\rangle =-\langle Ax_{k+1}, x_{k+1}\rangle<0.
$$
Therefore, if $X\in\CC^{n\times k}$ denotes the matrix the matrix with columns
$x_1,\ldots,x_{k+1}$, then
$$
k+1=\nu_-(X^*BX)\le \nu_-(B)=k,
$$
which is impossible. Therefore $Ax_{k+1}=0$, as claimed.
\end{proof}

   \begin{lem}\label{lem:3.90}
   Let (B1)--(B4) be in force, let $\nu_-(P)=\kappa$
and let $s\in {\cS}_\kappa(A_1, A_2,C,P)$. Then
   \begin{equation}\label{eq:3.108}
  \Delta_sFv=0\quad\text{a.e. on}\ \Omega_0\ \text{for every}\ v\in ker\ P.
\end{equation}
   \end{lem}

\begin{proof} Let $s\in{\cS}_\kappa(A_1, A_2,C,P)$ and let $P_s$ be
defined by (\ref{eq:3.7}). Then, since
$P_s=P$ by Corollary~\ref{cor:3.4}, the result follows from Lemma
\ref{lem:may31a9}, upon letting $u_1,\ldots,u_n$ be any orthonormal basis of
$\CC^n$ and setting $A$ and $B$ be the matrices with entries
$$
a_{ij}=
\langle \Delta_s
    FXu_j,FXu_i
   \rangle_{nst}
\quad\textrm{and}\quad
b_{ij}=\left\langle\left\{
    \Delta_s
    \left[\begin{array}{cc}
      0 & \Gamma_\ell  \\
      \Gamma_\ell ^* & 0 \\
   \end{array}\right]\Delta_s\right\}FXu_j,FXu_i\right
   \rangle_{nst},
$$
respectively.
\end{proof}

    \begin{lem}\label{lem:3.9}
    Let (B1)-(B4) be in force, and let $\nu_-(P)=\kappa$.
 Then $s\in{\cS}_\kappa(A_1, A_2,C,P)$ if and only if $s\in
    {\cS}_\kappa(A_1, A_2,C,P,X)$ and~\eqref{eq:3.108} holds.
    \end{lem}
\begin{proof} If $s\in{\cS}_\kappa(A_1, A_2,C,P)$, then clearly, $s\in
\widehat{\cS}_\kappa(A_1, A_2,C,P,X)$ and (\ref{eq:3.108}) is
implied by Lemma~\ref{lem:3.90}.

Conversely, if $s\in \widehat{\cS}_\kappa(A_1, A_2,C,P,X)$ and
\eqref{eq:3.108} holds, then, in view of Lemma \ref{lem:3.7}, (C1)
and (C2) are satisfied. Moreover, it follows from \eqref{eq:3.108}
and formula (\ref{eq:3.7}) that
\[
u^*P_sv=0\quad\textrm{for every }v\in\ker P\ \textrm{and}\ u\in\dC^n,
\]
i.e., $\textup{ker} P\subseteq \textup{ker} P_s$.
Thus, if $u=u_1+Xu_2$ with $u_1\in\ker P$ and  $u_2\in\CC^n$, then
\[
\begin{split}
u^*P_su&=(u_1+Xu_2)^*P_s(u_1+Xu_2)=u_2^*XP_sXu_2\\
&\le u_2^*XPXu_2=(u_1+Xu_2)^*P(u_1+Xu_2)=u^*Pu.
\end{split}
\]
This proves (C3) and hence, in view of Corollary~\ref{cor:3.4}, that
$P_s=P$. Therefore, $s\in{\cS}_\kappa(A_1, A_2,C,P)$.
\end{proof}

\begin{lem}\label{lem:3.10}
If (B1)--(B4) are in force, $W$ is defined by (\ref{eq:0.5}),
$\varepsilon\in \cS^{p\times q}$ and $s=T_W[\varepsilon]$, then
condition (\ref{eq:3.108}) holds if and only if
\begin{equation}\label{eq:3.133}
  \begin{bmatrix}I_p &\ -\varepsilon(\lam)\end{bmatrix}F(\mu)v=0
\end{equation}
for every  $\lam\in \sH_\varepsilon^+$ and every $v\in\ker\,P$,
where $\mu\in\Omega_0$ is the point chosen to normalize $W$ in
(\ref{eq:0.5}) (i.e., $W(\mu)=I_m$).
\end{lem}
\begin{proof} Let $s=T_W[\varepsilon]$ for $\varepsilon\in \cS^{p\times q}$ and
\begin{equation}
\label{eq:nov5a8}
R(\lam)=(w_{11}^\#(\lam)+\varepsilon(\lam)w_{12}^\#(\lam))^{-1}
\quad\textrm{for}\ \lam\in\Omega_W\cap\gh_\varepsilon^+.
\end{equation}
Then it follows from (\ref{eq:2.7}) and (\ref{eq:nov5a8}) that
\begin{equation}\label{eq:3.134}
  \begin{split}
 \begin{bmatrix}I_p &\ -s(\lam)\end{bmatrix}
&=(w_{11}^\#(\lam)+\varepsilon(\lam)w_{12}^\#(\lam))^{-1}
 \begin{bmatrix}I_p &\ -\varepsilon(\lam)\end{bmatrix}
j_{pq}W^\#(\lam)j_{pq}\\
  &=R(\lam) \begin{bmatrix}I_p &\ -\varepsilon(\lam)\end{bmatrix}
j_{pq}W^\#(\lam)j_{pq}
  \end{split}
\end{equation}
for $\lam\in\Omega_W\cap\sH_\varepsilon^+\cap\sH_s^+\subset \Omega_+$.
Using the formula
(\ref{eq:9.2}) one obtains
\begin{equation}\label{eq:3.135}
    \begin{split}
 \begin{bmatrix}I_p &\ -s(\lam)\end{bmatrix}
F(\lam)v&=R(\lam)
\begin{bmatrix}I_p &\ -\varepsilon(\lam)\end{bmatrix}j_{pq}W^\#(\lam)j_{pq}
F(\lam)v\\
&= R(\lam) \begin{bmatrix}I_p &\
-\varepsilon(\lam)\end{bmatrix}F(\mu)(I_n-XP)(M-\mu N)(M-\lam
N)^{-1}v.
\end{split}
\end{equation}  Since
\[
(I_n-XP)(M-\mu N)(M-\lam N)^{-1}X= (I_n-XP)X(M_0-\mu N_0)(M_0-\lam N_0)^{-1}=0,
\]
it follows that
$  \begin{bmatrix}I_p &\ -s(\lam)\end{bmatrix} F(\lam)v\equiv 0 $
for all $v\in\ker\,P$
if and only
$$
 \begin{bmatrix}I_p &\ -\varepsilon(\lam)\end{bmatrix}
F(\mu)(I_n-XP)\equiv 0,
$$
which is equivalent to~\eqref{eq:3.133}.
\end{proof}

\begin{lem}\label{lem:3.11}
If (B1) and (B2) are in force and $\mu\in\Omega_0$, then the subspace
$F(\mu)\ker\, P$ is a $j_{pq}$-neutral subspace of $\CC^n$.
\end{lem}
\begin{proof} If $v, w\in \ker\, P$, then  the
Lyapunov-Stein equation~\eqref{eq:0.LS} implies that
\[
\begin{split}
w^*F(\mu)^*j_{pq}F(\mu)v&=
w^*(M^*-\overline{\mu}N^*)^{-1}\{M^*P(M-\mu N)+(\mu M^*-N^*)PN\}
(M-\mu N)^{-1}v\\
&=w^*\{(M^*-\overline{\mu}N^*)^{-1}M^*P+\mu PN(M-\mu N)^{-1}\}v=0.
\end{split}
\]
This proves the statement.
\end{proof}

\begin{thm}\label{thm:3.13}
If (B1)--(B4) are in force and $\nu_-(P)=\kappa$, then there are
unitary matrices  $U\in \dC^{p\times p}$, $V\in \dC^{q\times q}$
such that
\begin{equation}\label{Sdescr}
\cS_\kappa(A_1,A_2,C,P)=\left\{T_W\left[U\left[\begin{array}{cc}
                                  \wt \varepsilon & 0 \\
                                        0 & I_\nu \\
                                      \end{array}  \right]V^*\right]:
                    \, \wt \varepsilon \in
                    \cS^{(p-\nu)\times(q-\nu)}\right\}
                    \cap \cS_\kappa^{\ptq},
\end{equation}
where $W$ and $\nu$ are given by \eqref{eq:0.5} and (\ref{nu}),
respectively.
\end{thm}

\begin{proof}
Let $s\in \cS_\kappa(A_1,A_2,C,P)$. Then, in view of
Lemma~\ref{lem:3.9}, $s$ is a solution of (X1)--(X3) such that
\eqref{eq:3.108} holds. By Theorem~\ref{thm:3.9} $s$ admits the
representation $s=T_{W}[\varepsilon]$, where $W$ is given
by~\eqref{eq:0.5}, $\varepsilon\in \cS_0^{p\times q}$ and, by
Lemma~\ref{lem:3.10},
\begin{equation}\label{eq:3.148}
\begin{bmatrix}I_p &\ -\varepsilon(\lam)\end{bmatrix}
F(\mu)u=0 \quad \text{for}\ u\in\ker P\ \text{and}\ \mu\in\Omega_0.
\end{equation}
Thus, if $\textup{col}\,(x,y)$ is a nonzero vector in the subspace
$\cK_\mu=F(\mu)\ker P$ for some choice of $x\in\CC^p$, $y\in\CC^q$ and
$\mu\in\Omega_0$, then, since
$\cK_\mu$ is $j_{pq}$-neutral,
$$
x=\varepsilon(\lambda)y\quad\text{and}\quad x^*x=y^*y.
$$
Therefore, (see \cite[Lemma 0.13]{D89})
$\varepsilon(\lambda)$ admits the representation
\[
   \varepsilon(\lambda)=U\left[\begin{array}{cc}
           \widetilde{\varepsilon}(\lambda) & 0 \\
          0 & I_\nu \\
         \end{array}\right] V^\ast,
\]
where $U\in \dC^{p\times p}$ and $V\in \dC^{q\times q}$ are unitary
matrices,  $\nu=\dim \cK_\mu$ and $\widetilde{\varepsilon}\in
\cS^{(p-\nu)\times(q-\nu)}$.

Conversely, if $\varepsilon$ is of the form~\eqref{epsilon}, then
\eqref{eq:3.148} holds; and $s\in
 \cS_\kappa(A_1,A_2,C,P)$ by Lemma~\ref{lem:3.11}.

Next, to verify (\ref{nu}), we first observe that if $\mu\in\Omega_0$,
then $M-\mu N$ defines an invertible map of
\begin{equation}\label{eq:3.500}
\cL=\ker PM\cap\ker PN\cap \ker C\quad\text{onto}\quad
\cL_\mu=\ker F(\mu)\cap \ker P.
\end{equation}

Let $u\in\CC^n$ and $v=(M -\mu N)^{-1}u$. Then it is readily checked that
$v\in\cL\Longrightarrow u\in\cL_\mu$, i.e., $(M -\mu N)\cL\subseteq \cL_\mu$.
Conversely, if $u\in\cL_\mu$, then clearly $P(M -\mu N)v=0$ and $Cv=0$.
Therefore, the Lyapunov-Stein equation~\eqref{eq:0.LS} implies that
\[
0=C^* j_{pq}Cv=M^*PMv-N^* PNv=M^*P(M-\mu N)v+(\mu M^*-N^*)PNv
\]
and hence that $PNv=0$. Thus, $(M -\mu N)\cL=\cL_\mu$ and
\begin{equation}\label{eq:3.154}
\dim\cL_\mu=\dim\cL=\dim\ker(M^*P^2M+N^*P^2N+C^*C).
\end{equation}
Therefore,
\[
\nu=\dim\cK_\mu=\dim\ker P -\dim\cL_\mu
= \dim\ker P-\dim\ker(M^*P^2M+N^*P^2N+C^*C),
\]
which is equivalent to formula (\ref{nu}).
\end{proof}

\section{Resolvent matrices}
\subsection{Pole and zero multiplicities.}
Let  $G(\lambda)$ be a $p\times q$ mvf that is meromorphic on $\Omega_+$ with
a Laurent expansion
\begin{equation}\label{Gexp}
    G(\lambda)=(\lambda-\lambda_0)^{-k}G_{-k}+\cdots+(\lambda-\lambda_0)^{-1}G_{-1}+G_0+\cdots
\end{equation}
 in a neighborhood of a pole $\lambda_0\in\Omega_+$. The pole
multiplicity $M_{\pi}(G, \lambda_0)$ is defined by (see~\cite{KL})
\begin{equation}\label{MTexp}
M_{\pi}(G, \lambda_0)=\rank T(G,\lam_0), \quad\textrm{where}\quad
T(G,\lam_0)=\begin{bmatrix}
  G_{-k} &  & {\bf 0} \\
  \vdots & \ddots &  \\
  G_{-1} &\hdots & G_{-k}
\end{bmatrix}.
\end{equation}
The pole multiplicity of $G$ over $\Omega_+$ is given by
\begin{equation}\label{eq:7.11p}
  M_{\pi}(G, \Omega_+)=\sum_{\lambda\in\Omega_+}M_{\pi}(G, \lambda).
\end{equation}
The zero multiplicity of a square  mvf $G$ over $\Omega_+$ is
defined by
\begin{equation}\label{eq:7.11z}
  M_{\zeta}(G, \Omega_+)=M_{\pi}(G^{-1}, \Omega_+)
  \quad \mbox{if} \quad \det G(\lambda)\not\equiv 0.
\end{equation}
Note that definitions~\eqref{eq:7.11p},~\eqref{eq:7.11z} of pole and
zero multiplicities coincide with those based on Smith-McMillan
representations of $G$ (see \cite{BGR}); and the degree of a Blaschke-Potapov
product $b$ of the form (\ref{BPprod}) coincides with
$M_{\zeta}(b, \Omega_+)$.

Let $H_{\kappa,\infty}^{p\times q}(\Omega_+)$ denote the class of $p\times q$
mvf's $G$ of the form
$G=F+B$, where $B$ is a  rational $p\times q$ mvf of pole
multiplicity $M_{\pi}(B, \Omega_+)\le\kappa$ and $F\in
H_{\infty}^{p\times q}(\Omega_+)$.

It follows from the results of \cite{KL}, \cite{Der03} (see also
\cite[Theorem 5.2]{Ki}) that every  mvf $G\in
H_{\kappa,\infty}^{p\times q}(\Omega_+)$ admits coprime
factorizations
\begin{equation}\label{eq:7.110}
  G(\lambda)=b_\ell(\lambda)^{-1}\varphi_\ell(\lambda)=\varphi_r(\lambda)b_r(\lambda)^{-1},
\end{equation}
where $b_\ell\in \cS_{in}^{p\times p}$, $b_r\in
{\cS}_{in}^{q\times q}$ are Blaschke-Potapov factors of degree
$M_{\pi} (G, \Omega_+)$ and $\varphi_\ell,\varphi_r\in
H_{\infty}^{p\times q}$.

Left and
right coprime factorizations may be characterized in terms of the pole and zero
multiplicities of its factors:

\begin{proposition}\label{Prop:5.1} {\rm (\cite[Proposition~3.4]{DerDym08} and \cite{BGR} in the rational case)}
Let $H_\ell, H_r\in H_{\infty}^{p\times q}$, $G_\ell, G_r\in
H_{\infty}^{p\times p}$ be a pair of mvf's such that $ G_\ell^{-1},
G_r^{-1}\in H_{\kappa,\infty}^{p\times p}$ for some
$\kappa\in\dN\cup\{0\}$. Then:
\begin{enumerate}
\item[(i)] The pair $G_\ell$, $H_\ell$ is left coprime
$\Longleftrightarrow$ $M_{\pi} (G_\ell^{-1}H_\ell, \Omega_+)=M_{\pi} (G_\ell^{-1},
\Omega_+)$. \item[(ii)] The pair $G_r$, $H_r$ is right coprime $\Longleftrightarrow$
$M_{\pi} (H_rG_r^{-1}, \Omega_+)=M_{\pi} (G_r^{-1}, \Omega_+)$.
\end{enumerate}
\end{proposition}
\subsection{Associated pairs.}
Let $b_1\in{\cS}_{in}^{\ptp}$ and $b_2\in \cS_{in}^{q\times  q}$ be
the  Blaschke-Potapov factors that are uniquely defined (up to right
and left constant unitary factors, respectively) by the coprime
factorizations
\begin{equation}\label{eq:7.12}
  \lambda C_{12}(I_{n_2}-\lambda A_2)^{-1}\begin{bmatrix}
   X_{21} & X_{22} \end{bmatrix}=b_1(\lambda)\varphi_1(\lambda),\quad
   \lambda\in\Omega_-\setminus\sigma(A_2)^{-1},\,
\varphi_1\in {\cR}\cap H_2^{p\times n}(\Omega_-),
\end{equation}
and
\begin{equation}\label{eq:7.11}
  C_{21}(A_1 -\lambda I_{n_1})^{-1}\begin{bmatrix}
  X_{11} & X_{12} \\
\end{bmatrix}%
=b_2(\lambda)^{-1}\varphi_2(\lambda), \quad
\lambda\in\Omega_+\setminus\sigma(A_1),\, \varphi_2\in
{\cR}\cap H_2^{q\times  n}(\Omega_+).
\end{equation}

\begin{remark}
\label{rem:nov13a8}
Since $(C_{21}, A_1)$ and $(C_{12}, A_2)$ are observable pairs and
$\sigma(A_1)\subset\DD$ and $\sigma(A_2)\subset\DD$, Theorem \ref{HW}
can be exploited to obtain Blaschke-Potapov products $\wt{b_1}(\lam)$
and $\wt{b_2}(\lam)$ such that
$$
\wt{b_1}(\lam)^{-1}\lambda C_{12}(I_{n_2}-\lambda A_2)^{-1}\in
H_\infty^{p\times n_2}(\Omega_-)\quad\textrm{and}\quad
\wt{b_2}(\lam)C_{21}(A_1 -\lambda I_{n_1})^{-1}\in
H_\infty^{q\times n_1}(\Omega_+).
$$
In particular, the Stein equations $Q_1-A_1^*Q_1A_1=C_{21}^*C_{21}$ and
$Q_2-A_2^*Q_2A_2=-C_{12}^*C_{12}$
have unique solutions
$$
Q_1=\sum_{j=0}^\infty (A_1^*)^jC_{21}^*C_{21}A_1^j
=\frac{1}{2\pi}\int_0^{2\pi}
(A_1^*-e^{-i\theta}I_{n_1})^{-1}C_{21}^*C_{21}
(A_1-e^{i\theta}I_{n_1})^{-1}d\theta
$$
and
$$
\quad Q_2=-\sum_{j=0}^\infty (A_2^*)^jC_{12}^*C_{12}A_2^j=
-\frac{1}{2\pi}\int_0^{2\pi}(I_{n_2}-e^{-i\theta}A_2^*)^{-1}C_{12}^*C_{12}
(I_{n_2}-e^{i\theta}A_2)^{-1}d\theta,
$$
respectively. Moreover, $Q_1$ is positive definite and $Q_2$ is negative
definite. Therefore, Theorem \ref{HW} (with
$C=C_{21}$, $M=A_1$, $N=I_{n_1}$ and $J=-I_q$) implies that
\begin{equation}
\label{eq:nov13a8}
\wt{b}_2(\lam)=I_q+(\lam-\mu)C_{21}(A_1-\mu I_{n_1})^{-1}Q_1^{-1}
(\lam A_1^*-I_{n_1})^{-1}C_{21}^*;
\end{equation}
see also Theorem 5.4 in \cite{Dym01} for additional discussion, if need be.
Similarly,
\begin{equation}
\label{eq:nov13b8}
\wt{b}_1(\lam)=I_p+(1-\overline{\mu}\lambda)C_{12}(I_{n_2}-\lam A_2)^{-1}
Q_2^{-1}(I_{n_2}-\overline{\mu}A_2^*)^{-1}C_{12}^*.
\end{equation}
Furthermore, if $F_{12}(\lam)=C_{12}(I_{n_2}-\lambda A_2)^{-1}$,
$$
\widetilde{\cM}_2=\{F_{12}u: u\in \CC^{n_2}\}\quad\text{and}\quad
{\cM}_2=\{F_{12}[X_{21}\quad X_{22}]u: u\in \CC^n\},
$$
then ${\cM}_2$ is a subspace of $\widetilde{\cM}_2$ that is invariant under
the backwards shift operator
$$
R_\alpha\, :\, f\rightarrow \frac{f(\lambda)-f(\alpha)}{\lambda-\alpha} \quad
\text{for}\ \alpha\in\DD.
$$
In particular,
\begin{eqnarray*}
R_0 F_{12}(\lambda)[X_{21}\quad X_{22}]
&=&F_{12}(\lambda)A_2[X_{21}\quad X_{22}]=
F_{12}(\lambda)[0\quad I_{n_2}]NX\\
&=&F_{12}(\lambda)[0\quad I_{n_2}]XPNX=F_{12}(\lambda)[X_{21}\quad X_{22}]PNX.
\end{eqnarray*}
Thus, upon endowing these two spaces with the normalized
standard inner product,
it follows from the either the arguments cited earlier in this remark or the
Beurling-Lax theorem that
$$
{\cM}_2=\cH(b_1)\subseteq \cH(\wt{b}_1)=\widetilde{\cM}_2.
$$
Similar considerations imply that if
$F_{21}(\lam)=C_{21}(A_1-\lambda I_{n_1})^{-1}$,
$$
\widetilde{\cM}_1=\{F_{21}u: u\in \CC^{n_1}\}\quad\text{and}\quad
{\cM}_1=\{F_{21}[X_{11}\quad X_{12}]u: u\in \CC^n\},
$$
then
$$
{\cM}_1=\cH_*(b_2)\subseteq \cH_*(\widetilde{b}_2)=\widetilde{\cM}_1.
$$
These conclusions also follow from Lemmas
\ref{lem:7.3} and \ref{lem:7.4}, below, which are based on state space
calculations. Moreover,
$$
b_1(\lam)=\wt{b}_1(\lam)\quad\textrm{if}\
\rank [X_{21}\quad X_{22}]=n_2\quad\textrm{and}\quad
b_2(\lam)=\wt{b}_2(\lam)\quad\textrm{if}\ \rank [X_{11}\quad X_{12}]
=n_1.
$$
\end{remark}

\begin{lem}\label{lem:7.5} Let the data set $(A_1, A_2, C, P)$ satisfy the
assumptions (B1), (B2) and (B4). Then:
\begin{equation}\label{eq:7.13}
  b_2^{-1}H_2^q=\{C_{21}(A_1-\lambda I_{n_1})^{-1}\begin{bmatrix}
  X_{11} & X_{12} \\
\end{bmatrix}%
u+h_+: u\in\dC^n, h_+\in H_2^q \},
\end{equation}
\begin{equation}\label{eq:7.14}
  b_1(H_2^p)^{\bot}=\{C_{12}(I_{n_2}-\lambda A_2)^{-1}\begin{bmatrix}
   X_{21} & X_{22} \end{bmatrix}u+h_-: u\in\dC^{n}, h_-\in(H_2^p)^{\bot}\}.
\end{equation}
\end{lem}
\begin{proof}
The first step in the verification of (\ref{eq:7.13}) is to
observe that
if $g\in H_2^{n_1}$ and $v\in\CC^{n_1}$, then
\begin{equation}
\label{eq:mar17c8}
(A_1-\lambda I_{n_1})^{-1}(g(\lambda)-v)\in H_2^{n_1}\Longleftrightarrow
v=\frac{1}{2\pi i}\int_{\dT}(\zeta I_{n_1}-A_1)^{-1}g(\zeta)d\zeta\,.
\end{equation}
Next, recall that since $b_2$ and $\varphi_2$ are left coprime, there exist
a pair of mvf's $c\in H_\infty^{q\times q}$ and $d\in H_\infty^{n\times q}$
such that $b_2c+\varphi_2 d=I_q$ in $\Omega_+$. Let $f\in H_2^q$ and let
$$
g=[X_{11}\quad X_{12}]df\quad\textrm{and}\quad
v=\frac{1}{2\pi i}\int_{\TT}(\zeta I_{n_1}-A_1)^{-1}g(\zeta)d\zeta.
$$
Then
$$
b_2^{-1}f= b_2^{-1}(b_2c+\varphi_2 d)f=cf +C_{21}(A_1-\lambda I_{n_1})^{-1}
(g-v)+C_{21}(A_1-\lambda I_{n_1})^{-1}v,
$$
and since $g\in H_2^{n_1}$, the first two terms on the far right belong to
$H_2^q$. Moreover, since
\begin{eqnarray*}
(A_1-\zeta I_{n_1})^{-1}[X_{11}\quad X_{12}]&=&[I_{n_1}\quad 0]
(M-\zeta N)^{-1}X=[I_{n_1}\quad 0]XP(M-\zeta N)^{-1}X\\
&=&[X_{11}\quad X_{12}]P(M-\zeta N)^{-1}X,
\end{eqnarray*}
it follows that $v=[X_{11}\quad X_{12}]w$ for some $w\in\CC^n$ and hence that
$$
C_{21}(A_1-\lambda I_{n_1})^{-1}v=b_2^{-1}\varphi_2w=C_{21}(A_1-\lambda I_{n_1})^{-1}[X_{11}\quad X_{12}]w.
$$
Thus, the left-hand side of (\ref{eq:7.13}) is
a subset of the right-hand side. The opposite inclusion is easy to check and
is left to the reader, as is the verification of (\ref{eq:7.14}).\
\end{proof}

\begin{lem}\label{lem:7.6}
Let the assumptions of Lemma~\ref{lem:7.5} be in force and let $\wt{b}_1$
and $\wt{b}_2$ be defined by (\ref{eq:nov13b8}) and (\ref{eq:nov13a8}),
respectively.
\begin{enumerate}
\item[\rm(1)]
If the pair $(C_{12}, A_2)$ is observable, then
\begin{equation}
\label{eq:mar17a8}
\{P_-(\lambda I_{n_2}-A_2^*)^{-1}h:\,h\in H_2^{n_2}\}
=\{P_-(\lambda I_{n_2}-A_2^*)^{-1}C_{12}^*g:\,
g\in{\cH}(\wt{b_1})\}.
\end{equation}
\item[\rm(2)]
If the pair $(C_{21}, A_1)$ is observable, then
\begin{equation}
\label{eq:mar17b8}
\{P_+(A_1^*-\lambda I_{n_1})^{-1}h:\, h\in (H_2^{n_1})^\perp\}=
\{P_+(A_1^*-\lambda I_{n_1})^{-1}C_{21}^*g:\, g\in{\cH}_*(\wt{b_2})\}.
\end{equation}
\end{enumerate}
\end{lem}

\begin{proof}
If $h\in H_2^{n_2}$ and
\[
x=\frac{1}{2\pi i}\int_{\dT}(\zeta I_{n_2}-A_2^*)^{-1}h(\zeta)d\zeta, \quad
\text{then}\quad
P_-(\lambda I_{n_2}-A_2^*)^{-1}\{h(\lambda)-x\}=0.
\]
Therefore,
\begin{equation}
\label{eq:nov21a8}
\{P_-(\lambda I_{n_2}-A_2^*)^{-1}h:\,h\in H_2^{n_2}\}\subseteq
\{(\lambda I_{n_2}-A_2^*)^{-1}x:x\in\CC^{n_2}\}.
\end{equation}
On the other hand, if $g=C_{12}(I_{n_2}-\lam A_2)^{-1}u$ for some
$u\in\CC^{n_2}$, then
$$
P_-(\lambda I_{n_2}-A_2^*)^{-1}C_{12}^*g=(\lambda I_{n_2}-A_2^*)^{-1}v,
$$
where
\begin{eqnarray*}
v&=&\frac{1}{2\pi i}\int_{\TT}(\zeta I_{n_2}-A_2^*)^{-1}C_{12}^*C_{12}
(I_{n_2}-\zeta A_2)^{-1}d\zeta u\\
&=&\frac{1}{2\pi}\int_0^{2\pi}(I_{n_2}-e^{-i\theta}A_2^*)^{-1}C_{12}^*C_{12}
(I_{n_2}-e^{i\theta}A_2)^{-1}d\theta u =-Q_2u,
\end{eqnarray*}
where $Q_2$ is the negative definite matrix introduced in Remark
\ref{rem:nov13a8}. Thus, as $Q_2$ is invertible, and $g\in\cH(\wt
b_1)$
\begin{equation}
\label{eq:nov21b8}
\{(\lambda I_{n_2}-A_2^*)^{-1}x:x\in\CC^{n_2}\}\subseteq
\{P_-(\lambda I_{n_2}-A_2^*)^{-1}C_{12}^*g:\,
g\in{\cH}(\wt{b_1})\}.
\end{equation}
This completes the proof of (1). The proof of (2) is similar.
\end{proof}

\begin{lem}\label{lem:7.3}
Let the data set $(A_1, A_2, C, P)$ satisfy the
assumptions (B1)--(B4) and let $W$
and $F_1$ be given by~\eqref{eq:0.5} and  \eqref{F12}, respectively.
Then:
\begin{enumerate}
\item[(i)]
 $P_-WH_2^m=\{F_1\begin{bmatrix}
  X_{11} & X_{12} \\
\end{bmatrix}%
u:
u\in\dC^n \}$.
\item[(ii)]
$WH_2^m=\{F_1 \begin{bmatrix}
  X_{11} & X_{12} \\
\end{bmatrix}%
u+h: u\in\dC^n, h\in H_2^m\,\&\,P_-W^{-1}(F_1 \begin{bmatrix}
  X_{11} & X_{12} \\
\end{bmatrix}%
u+h)=0 \}$.
\item[(iii)]
$f\in H_2^m$, $[0\quad I_q]Wf\in H_2^q$ $\Longrightarrow$ $Wf\in H^m_2$.
\item[(iv)]
$[0\quad I_q]WH_2^m=b_2^{-1}H_2^q$.
\item[(v)]
$[0\quad I_q]W(\cR\cap H_2^m)=b_2^{-1}(\cR\cap H_2^q)$.
\end{enumerate}
\end{lem}
\begin{proof}
(i) Since
\begin{equation}
\label{eq:nov5b8}
\wtilde{C}^*=(\overline{\mu}M-N)XF(\mu)^*=XP(\overline{\mu}M-N)XF(\mu)^*,
\end{equation}
it is readily checked with the help of
formulas (\ref{eq:7.09}) and  (\ref{eq:mar17c8}) that
\[
 P_-WH_2^m\subseteq P_-F_1[X_{11}\quad X_{12}]H_2^n\subseteq
\{{\bC}_1(A_1-\lambda I_{n_1})^{-1}\begin{bmatrix}
  X_{11} & X_{12} \\
\end{bmatrix}%
u: u\in\dC^n \}.
\]
More precisely, if $h\in H_2^m$, then, in view of (\ref{eq:7.09}) and
(\ref{eq:mar17c8}),
$$
P_-Wh=P_-F_1\wt{\bC}_1^*j_{pq}h={\bC}_1(A_1-\lam I_{n_1})^{-1}x,
$$
where
$$
x=\frac{1}{2\pi i}\int_{\TT}(\zeta I_{n_1}-A_1)^{-1}\wt{\bC}_1^*
j_{pq}h(\zeta)d\zeta.
$$
Thus, if $h=j_{pq}\wtilde{\bC}_1(I_{n_1}-\zeta A_1^*)^{-1}u$,
then
$$
x=Qu,\quad\textrm{where}\quad Q=-\frac{1}{2\pi i}\int_{\TT}
(\zeta I_{n_1}-A_1)^{-1}\wtilde{\bC}_1^*
\wtilde{\bC}_1(I_{n_1}-\zeta A_1^*)^{-1}d\zeta.
$$
But this serves to complete the proof of (i), since
$\textup{rng}\,Q=\textup{rng}[X_{11}\quad X_{22}]$ by Corollary \ref{cor:7.1}.

(ii) For every $f\in H_2^m$ there is a vvf $h\in H_2^m$  such that
\[
Wf=P_-(Wf)+h.
\]
Thus, in view of (i), there is a vector $u\in\dC^n$ such that
\[
Wf={\bC}_1(A_1-\lam I_{n_1})^{-1}\begin{bmatrix}
  X_{11} & X_{12} \\
\end{bmatrix}
u+h=g\quad\text{and}\quad P_-W^{-1}g=0.
\]

(iii) Since the pair $(C_{21}, A_1)$ is observable the third statement is
immediate from (i) and (ii).

(iv) The inclusion $[0\quad I_q]WH_2^m\subseteq b_2^{-1}H_2^q$ follows
readily from (ii) and formula (\ref{eq:7.11}). To verify the opposite
inclusion, let $f\in H_2^q$. Then, by Lemma \ref{lem:7.5},
\begin{equation}
\label{eq:nov22a8}
b_2^{-1}f=C_{21}(A_1-\lambda I_{n_1})^{-1}[X_{11}\quad X_{12}]u+h_2
\end{equation}
for some choice of $u\in\CC^n$ and $h_2\in H_2^q$. In view of Lemma
\ref{Odecom}, we may
assume that $[X_{21}\quad X_{22}]u=0$ and hence, upon writing
$u=\textup{col}(u_1,u_2)$ with $u_j\in\CC^{n_j}$ for $j=1,2$, that
\begin{equation}
\label{eq:nov21c8}
P_-W^{-1}\left(F_1[X_{11}\quad X_{12}]u+\begin{bmatrix}h_1\\h_2\end{bmatrix}
\right)=P_-\widetilde{\bC}_2(\lambda
I_{n_2}-A_2^*)^{-1}\{u_2+(C^*_{12}h_1-C^*_{22}h_2)\}
\end{equation}
for every choice of $h_1\in H_2^p$, since
\[
P_-W^{-1} {\bC}_1(A_1 -\lambda I_{n_1})^{-1}\begin{bmatrix}
  X_{11} & X_{12}
\end{bmatrix}u=P_-\widetilde{\bC}_2(\lambda
I_{n_2}-A_2^*)^{-1}u_2,
\]
by formula (\ref{eq:9.1}) and
$$
P_-W^{-1}\begin{bmatrix}h_1\\h_2\end{bmatrix}
=P_-\wt{\bC}_2(A_2^*-\lam I_{n_2})^{-1}(C^*_{12}h_1-C^*_{22}h_2).
$$
by formula (\ref{eq:7.09A}). Now, upon choosing $h_1\in
\cH(\wt{b_1})$ so that the right hand side of (\ref{eq:nov21c8}) is
equal to zero, it follows from (ii) that
\begin{equation}
\label{eq:nov21d8}
F_1[X_{11}\quad X_{12}]u+\begin{bmatrix}h_1\\h_2\end{bmatrix}\in WH_2^m
\end{equation}
and hence that
$$
b_2^{-1}f\in[0\quad I_q]WH_2^m,
$$
as needed.
\end{proof}

\begin{corollary}\label{thm:7.8}
If $W$ is defined by (\ref{eq:0.5}) and $b_2$ is defined
by~\eqref{eq:7.11}, then there exist a pair of rational mvf's
 $g_1\in H_\infty^{p\times q}$ and $g_2\in H_\infty^{q\times q}$ such that
\begin{equation}\label{g12}
w_{21}g_1+w_{22}g_2=b_2^{-1}.
\end{equation}
Moreover, if $(g_1, g_2)$
is any pair of functions in
$H_\infty^{p\times q}\times H_\infty^{q\times q}$
that satisfies equation
(\ref{g12}), then the mvf
\begin{equation}\label{K}
    K=(w_{11}g_1+w_{12}g_2)(w_{21}g_1+w_{22}g_2)^{-1}=(w_{11}g_1+w_{12}g_2)b_2
\end{equation}
is holomorphic in $\Omega_+$.
\end{corollary}

\begin{proof}
By Lemma~\ref{lem:7.3} (v) for each $f\in \cR\cap H_2^q$ there
exists $f_1\in \cR\cap H_2^p$ and $f_2\in \cR\cap H_2^q$ such that
\[
w_{21}f_1+w_{22}f_2=b_2^{-1}f
\]
which leads easily to~\eqref{g12} by successively choosing columns
of $I_q$.

The second assertion is implied by (iii) of Lemma~\ref{lem:7.3}.
\end{proof}

There is an analogue of Lemma \ref{lem:7.3} that focuses on $(H_2^m)^\perp$
that we state without proof.

\begin{lem}\label{lem:7.4}
Let the data set $(A_1, A_2, C, P)$ satisfy the
assumptions (B1)--(B4)
and let $W$
and $F_2$ be given by~\eqref{eq:0.5} and  \eqref{F12}, respectively.
Then:
\begin{enumerate}
\item[(i)]
 $P_+W(H_2^m)^\perp=\{F_2(\lambda )\begin{bmatrix}
  X_{21} & X_{22} \\
\end{bmatrix}%
u: u\in\dC^n \}$.
\item[(ii)]
$W(H_2^m)^\perp=\{F_2[X_{21} \ X_{22}]u
+h: u\in\dC^n, h\in (H_2^m)^\perp\,
\&\,P_+W^{-1}(F_2[X_{21} \ X_{22}]u+h)=0 \}$.
\item[(iii)] If $f\in (H_2^m)^\perp$ and  $[I_p\quad 0]Wf\in (H_2^p)^\perp$, then
$Wf\in (H_2^m)^\perp$.
\item[(iv)]
$[I_p\quad 0]W(H_2^m)^\perp=b_1(H_2^p)^\perp$.
\end{enumerate}
\end{lem}

\begin{lem}\label{lem:7.6a}
Let $W\in\cU_\kappa(j_{pq})$ and let the conditions (iii) and (iv)
of Lemma~\ref{lem:7.3} be in force. Then $s_{21}\in
\cS_\kappa^{\ptq}$ and hence $W\in\cU_\kappa^\circ(j_{pq})$.
\end{lem}

\begin{proof}
  Our first objective is to show that
\begin{equation}\label{eq:7.114}
    \dim P_-\begin{bmatrix}
  s_{11} \\ s_{21}
\end{bmatrix}H^p_2= \dim P_-s_{21}H_2^p.
\end{equation}
To prove the inequality $\le$ assume that $ P_-s_{21}g_1=0 $ for some $g_1\in H_2^p$. Then
$g_2=s_{21}g_1\in H_2^q$ and, hence,
\begin{equation}\label{eq:7.115}
    \begin{bmatrix}
  s_{11}g_1 \\ 0
\end{bmatrix}=\begin{bmatrix}
  w_{11}g_1+w_{12}g_2 \\ w_{21}g_1+w_{22}g_2
\end{bmatrix}=W\begin{bmatrix}
  g_1 \\ g_2
\end{bmatrix}\in W H_2^m.
\end{equation}
In view of Lemma~\ref{lem:7.3} (iii),  $W\begin{bmatrix}
  g_1 \\ g_2
\end{bmatrix}\in H_2^m$.
Therefore $s_{11}g_1\in H_2^p$ and hence the asserted inequality is justified.
This completes the proof of~\eqref{eq:7.114}, since the opposite
inequality is self-evident.

Next, to prove that
\begin{equation}\label{eq:s21col}
    P_-\begin{bmatrix}
  s_{12} \\ s_{22}
\end{bmatrix}H^p_2\subseteq P_-\begin{bmatrix}
  s_{11} \\ s_{21}
\end{bmatrix}H^p_2,
\end{equation}
let $h$ be an arbitrary vvf from $H_2^q$. Then,
by Lemma~\ref{lem:7.3} (iv), there exists a pair of vvf's $g_1\in H_2^p$
and $g_2\in H_2^q$ such that
\begin{equation}\label{eq:7.120}
    h=w_{21}g_1+w_{22}g_2.
\end{equation}
Then 
\[
\begin{bmatrix}
  w_{11}g_1+w_{12}g_2 \\w_{21}g_1+w_{22}g_2
\end{bmatrix}=\begin{bmatrix}
  w_{11}g_1+w_{12}g_2 \\ h
\end{bmatrix}\in H_2^m
\]
due to Lemma~\ref{lem:7.3} (iii).  Therefore,
$g_2=-w_{22}^{-1}w_{21}+w_{22}^{-1}h$, and hence
$$
\begin{bmatrix}
  s_{11}g_1 \\ s_{21}g_1
\end{bmatrix}+
\begin{bmatrix}
  s_{12}h \\ s_{22}h
\end{bmatrix}=
\begin{bmatrix}
  w_{11}g_1+w_{12}g_2 \\ g_2
\end{bmatrix}\in H_2^m.
$$
Thus, for any $h\in H_2^q$ there exists a vvf $g_1\in H_2^p$ such that
\[
P_-\begin{bmatrix}
  s_{11}g_1 \\ s_{21}g_1
\end{bmatrix}=-
P_-\begin{bmatrix}
  s_{12}h \\ s_{22}h
\end{bmatrix},
\]
which justifies~\eqref{eq:s21col}.

Formulas~\eqref{eq:7.114} and \eqref{eq:s21col} imply
\[
\dim P_-s_{21}H_2^p=\dim P_-SH_2^m=\kappa.
\]
and, hence, since $s_{21}$ is contractive on $\Omega_0$ than
$s_{21}\in \cS_\kappa^{\qtp}$.

The last statement now follows from the definition of the class
$\cU_\kappa^\circ(j_{pq})$ (which is given in the Introduction).
\end{proof}

\begin{remark}\label{IndS12}
If the data set $(A_1, A_2, C, P)$ satisfies  assumptions
(B1)--(B4) and the mvf $W$ given by (\ref{eq:0.5}) belongs to
$\cU_{\kappa_1}(j_{pq})$, then
conditions (iii) and (iv) of Lemma~\ref{lem:7.3}
are in force and, hence, $s_{21}=-w_{22}^{-1}w_{21}\in
\cS_{\kappa_1}^{\qtp}$.  This fact was proved  in~\cite{BGR} for
invertible $P$ when (B1)-(B3) hold.
\end{remark}

\begin{thm}\label{factofW}
Let the data set $(A_1, A_2, C, P)$ satisfy assumptions
(B1)--(B4), let $W$ be given by~\eqref{eq:0.5}, and let the pair
$\{b_1,b_2\}$ and the mvf $K\in H_\infty^{\ptq}$ be defined
by~\eqref{eq:7.12}, \eqref{eq:7.11}
and \eqref{K}. Then:
\begin{enumerate}
\item[\rm(1)]
 The mvf $W$ admits the  factorization
\begin{equation}\label{eq:7.41}
W=\begin{bmatrix}b_1& 0\\
                           0 & b_2^{-1}\end{bmatrix}
\begin{bmatrix}\wt\varphi_{11}&\wt\varphi_{12}\\
                \varphi_{21}&\varphi_{22}\end{bmatrix}\quad \textrm{a.e. in}\
                \Omega_0
\end{equation}
where the pair $(\wt\varphi_{11}, \wt\varphi_{12}) \in
\cR\cap  H_\infty^{p\times p}(\Omega_-)\times \cR\cap  H_\infty^{p\times q}(\Omega_-)$ is left coprime over $\Omega_-$ and the pair $(\varphi_{21}, \varphi_{22})\in \cR\cap  H_\infty^{q\times p}(\Omega_+) \times\cR\cap  H_\infty^{q\times q}(\Omega_+)$ is left coprime over $\Omega_+$.
\vspace{2mm}
\item[\rm(2)]
The pair $\{b_1,b_2\}$ is an associated pair for $W$.
\vspace{2mm}
\item[\rm(3)]
$W$ admits the factorizations
\begin{equation}
\label{eq:0.6a}
W=\Theta\,\Phi\quad\text{and}\quad W=\wt{\Theta}\,\wt{\Phi},
\end{equation}
over $\Omega_+$ and $\Omega_-$, respectively, with
\begin{equation}
\label{eq:0.7a}
\Theta=\begin{bmatrix}b_1&Kb_2^{-1}\\0&b_2^{-1}\end{bmatrix},\quad
\wt{\Theta}=\begin{bmatrix}b_1&0\\K^\# b_1&b_2^{-1}\end{bmatrix}
\end{equation}
\begin{equation}
\label{eq:0.8a}
\Phi= \begin{bmatrix}\varphi_{11}& \varphi_{12}\\
\varphi_{21}&\varphi_{22}\end{bmatrix}\in\cR\cap H_\infty^{m\times m}
(\Omega_+)\quad\text{and}\quad
\wt{\Phi}= \begin{bmatrix}\wt\varphi_{11}& \wt\varphi_{12}\\
\wt\varphi_{21}&\wt\varphi_{22}\end{bmatrix}\in\cR\cap H_\infty^{m\times m}
(\Omega_-).
\end{equation}
Moreover, $\Phi$ is outer in $H_\infty^{m\times m}(\Omega_+)$,  $\wt{\Phi}$
is outer in $H_\infty^{m\times m}(\Omega_-)$ and
\begin{equation}
\label{eq:sep22a8}
\wt{\Theta}^\#j_{pq}\Theta=j_{pq}.
\end{equation}
\end{enumerate}
\end{thm}
\begin{proof}
Formulas~\eqref{eq:7.09} and (\ref{eq:nov5b8}) yield the representation
\[
\begin{bmatrix}w_{21}(\lam)& w_{22}(\lam)\end{bmatrix}=C_{21}(A_1-\lambda I_{n_1})^{-1}[X_{11}\quad X_{12}]U+v(\lambda),
\]
for some choice of $U\in\CC^{n\times m}$ and $v\in \cR\cap H_\infty^{q\times m}(\Omega_+)$.
Lemma~\ref{lem:7.5} then guarantees the existence of a factorization of the
form
\begin{equation}\label{eq:7.41a}
 \begin{bmatrix}w_{21}& w_{22}\end{bmatrix}=
b_2^{-1}\begin{bmatrix}\varphi_{21}&\varphi_{22}\end{bmatrix}
\end{equation}
with $\varphi_{21}\in \cR\cap H_\infty^{q\times p}(\Omega_+)$ and
$\varphi_{22}\in \cR\cap H_\infty^{q\times q}(\Omega_+)$.
By Lemma~\ref{lem:7.3} there are mvf's $g_1\in
\cR\cap H_\infty^{p\times q}$ and  $g_2\in \cR\cap H_\infty^{q\times q}$ such that
\[
\varphi_{21}g_1+\varphi_{22}g_2=b_2(w_{21}g_1+w_{22}g_2)=I_q.
\]
Thus, Lemma~\ref{Corona} implies the rank condition
\[
\mbox{rank }\begin{bmatrix}\varphi_{21}(\lambda)
&\varphi_{22}(\lambda)\end{bmatrix}=q\quad\textrm{for}\ \lambda\in\Omega_+.
\]
Therefore, the factorization $\varphi_{22}^{-1}\varphi_{21}$ is left coprime over $\Omega_+$ and, as
$$
\textup{rank }[b_2(\lam)\quad \varphi_{21}(\lam)\quad
\varphi_{22}(\lam)]=q\quad\textrm{for}\ \lam\in\Omega_+,
$$
the factorization \eqref{eq:7.41a} is also left coprime over $\Omega_+$.

Similarly, the coprimeness of the factorizations
\begin{equation}\label{eq:7.43}
\wt{\varphi}_{11}^{-1}\wt{\varphi}_{12}\quad\textrm{and}\quad
 \begin{bmatrix}w_{11} & w_{12}\end{bmatrix}
 =b_1\begin{bmatrix}\wt\varphi_{11}& \wt\varphi_{12}\end{bmatrix}=
(b_1^\#)^{-1}\begin{bmatrix}\wt\varphi_{11}& \wt\varphi_{12}\end{bmatrix}
\quad\textrm{over}\
\Omega_-
\end{equation}
follows from (\ref{eq:7.09}),~\eqref{eq:7.14} and Lemma~\ref{lem:7.4}. This
completes the proof of (1).

Next, (2) follows by comparing the coprime
factorizations~\eqref{eq:7.41a} and \eqref{eq:7.43} with those
in~\cite[Corollary~4.15]{DerDym08}.

Finally, (3) follows from~\cite[Theorem~4.12]{DerDym08}.
\end{proof}

\begin{remark}
\label{rem:jun19a9}
If $W\in\cU_{\kappa_1}^\circ(j_{pq})$, $S=PG(W)=[s_{ij}]_{i,j=1}^2$
and the mvf's $c_\ell\in H_\infty^{\qtq}$, $d_\ell\in H_\infty^{\ptq}$ are
related to the factors in the Kre\u{\i}n-Langer factorizations~\eqref{eq:0.4a}
of $s_{21}$ by  the condition
\[
\gb_\ell c_\ell+\gs_\ell d_\ell=I_q,
\]
then (see e.g., \cite{DerDym08})  $W$ admits the
factorizations~\eqref{eq:0.6a}, \eqref{eq:0.7a}, where the mvf's
$b_1$, $b_2$ are defined by~\eqref{eq:0.8} and \eqref{eq:0.7}, and
\begin{equation}\label{eq:11.28}
    K=(-w_{11}d_\ell+w_{12}c_\ell)(-w_{21}d_\ell+w_{22}c_\ell)^{-1}.
\end{equation}
Moreover, if $\kappa_1=0$, then the second matrix on the right hand side of
(\ref{eq:7.41}) is a (right) $\gamma$-generating matrix (in the sense of
Arov); see Chapter 7 of \cite{ArovD08} for the definition and references.
\end{remark}

\section{Parametrization of solutions}

In this section we will give a parametrization of the set of
solutions of the problem (C1)--(C3) in terms of 
the linear fractional transformation $T_W$. The main result of this
section is based on Theorem~\ref{thm:3.13} augmented by the
factorization result of Theorem~\ref{factofW} and a special case of
the Kre\u{\i}n-Langer generalization of Rouche's Theorem, which is
formulated below.

\subsection{Counting zeros.}
The analysis in this subsection is valid for $\Omega_+=\DD$ and
$\Omega_+=\Pi_+$. We will need the notation
$$
\widetilde{L}_r^{\ptq}=\left\{\begin{array}{ll} L_r^{\ptq}
&\quad\text{if}\quad
\Omega_+=\DD;\\ \\
\{f:\,(1+\vert\mu\vert^2)^{-1/r}f\in
L_1^{\ptq}\}&\quad\text{if}\quad \Omega_+=\Pi_+.
\end{array}\right., \quad r=1,2.
$$

\begin{thm}\label{Rouche}
{\rm\cite{KL81}}
    Let $\varphi,\psi\in H_\infty^{q\times q}$, $\det(\varphi+\psi)\not\equiv 0$ in $\Omega_+$, $M_\zeta(\varphi,\Omega_+)<\infty$, and
\begin{equation}\label{eq:8.3}
    \|\varphi(t)^{-1}\psi(t)\|\le 1\quad \mbox{a.e. on }\Omega_+.
\end{equation}
Then $M_\zeta(\varphi+\psi,\Omega_+)\le M_\zeta(\varphi,\Omega_+)$
with equality if
\begin{equation}\label{eq:8.4}
    (\varphi +\psi)^{-1}\varphi|_{\Omega_0} \in \wt L_1^{q\times q}.
\end{equation}
\end{thm}
This theorem is used to estimate the zero multiplicity of the
denominator in linear-fractional transformation $T_W$ associated
with the data set $(A_1,A_2,C,P)$.
\begin{lem}\label{lem:5.3}
Let (B1)--(B4) be in force, 
let $\kappa_1=\nu_-(P)$, let the mvf's $W$, $b_1$, $b_2$ be defined
by~\eqref{eq:0.5}, \eqref{eq:7.11}, \eqref{eq:7.12}, let
$\varphi_{21}$, $\varphi_{22}$, $\wt\varphi_{11}$, $\wt\varphi_{12}$
be defined by~\eqref{eq:7.41a}, \eqref{eq:7.43} and let
$\varepsilon\in
\cS^{\ptq}_{\kappa_2}$ 
admit the
following Kre\u{\i}n-Langer factorizations
\begin{equation}\label{KLepsilon}
\varepsilon=\theta^{-1}_\ell  \varepsilon_\ell =\varepsilon_r\theta^{-1}_r.
\end{equation}
Then
\begin{equation}\label{eq:8.5}
M_\zeta(\varphi_{21}\varepsilon_r+\varphi_{22}\theta_r,\Omega_+)= M_\zeta(\theta_\ell
\wt\varphi_{11}^\#+\varepsilon_\ell \wt\varphi_{12}^\#,\Omega_+)=\kappa_1+\kappa_2.
\end{equation}
\end{lem}
\begin{proof}
1) By (1) of Theorem~\ref{factofW} and (3) of Lemma~\ref{lem:7.6a},
the factorization $\varphi_{22}^{-1}\varphi_{21}$ is left coprime
and $\varphi_{22}^{-1}\varphi_{21}\in \cS_{\kappa_1}^{p\times q}$.
Thus,  Proposition~\ref{Prop:5.1} guarantees that
 $M_\zeta(\varphi_{22},\Omega_+)=\kappa_1$. Moreover,  the identities
\begin{equation}\label{eq:7.52}
    w_{22}(\zeta)w_{22}(\zeta)^*-w_{21}(\zeta)w_{21}(\zeta)^*\equiv I_q \quad
    (\zeta\in\Omega_0),
\end{equation}
\[
w_{22}^{-1}(\zeta)w_{21}(\zeta)=\varphi_{22}^{-1}(\zeta)\varphi_{21}(\zeta) \quad
(\zeta\in\Omega_0)
\]
imply that
\[
\|\varphi_{22}(\zeta)^{-1}\varphi_{21}(\zeta)\|\le \rho<1\quad
(\zeta\in\Omega_0)
\]
for some $\rho<1$. Therefore,
\begin{equation}\label{eq:7.54}
\begin{split}
\|(\varphi_{21}(\zeta)\varepsilon(\zeta)+\varphi_{22}(\zeta))^{-1}\|\
&\le\|\varphi_{22}(\zeta)^{-1}\|
\|(I_q+\varphi_{22}^{-1}(\zeta)\varphi_{21}(\zeta)\varepsilon(\zeta))^{-1}\|\\
&\le\frac{\|\varphi_{22}(\zeta)^{-1}\|}{1-\rho}\le L\quad\text{for
some $L<\infty$}.
\end{split}
\end{equation}

2) Let
\[
\psi=\varphi_{21}\varepsilon_r,\quad \text{and}\quad
\varphi=\varphi_{22}\theta_r.
\]
Then $M_\zeta(\varphi,\Omega_+)=\kappa_1+\kappa_2$,
since $M_\zeta(\varphi_{22},\Omega_+)=\kappa_1$ and
 $M_\zeta(\theta_r,\Omega_+)=\kappa_2$.   Moreover,
\[
\|\varphi(t)^{-1}\psi(t)\|\le 1\quad (t\in\Omega_0)
\]
and, since $\theta_r$ is inner, it follows from~\eqref{eq:7.54} that
\[
\|\theta_r(t)(\varphi(t)+\psi(t))^{-1}\|=
\|(\varphi_{21}(t)\varepsilon(t)+\varphi_{22}(t))^{-1}\|\le L\quad
\text{for}\ t\in\Omega_0.
\]
Therefore, by Theorem~\ref{Rouche},
\[
M_\zeta(\varphi_{21}\varepsilon_r+\varphi_{22}\theta_r,\Omega_+)
=M_\zeta(\varphi_{22}\theta_r,\Omega_+)=\kappa_1+\kappa_2.
\]
\end{proof}

\begin{remark}
The equalities~\eqref{eq:8.5} were proved in~\cite{DerDym08} under the
less restrictive assumption that $W\in
\cU_{\kappa_1}^\circ(j_{pq})\cap \wt L_2^{\mtm}$.
\end{remark}
\begin{lem}\label{lem:5.3a}
Let $W\in \cU_{\kappa_1}^\circ(j_{pq})\cap \wt L_2^{\mtm}$, let the
mvf $\Phi$ be defined by~\eqref{eq:0.6a}, \eqref{eq:0.7a}, let
$\varepsilon\in \cS_{\kappa_2}^{p\times q}$ and let
\begin{equation}\label{G}
\ G=T_\Phi[\varepsilon]=(\varphi_{11}\varepsilon+\varphi_{12})
     (\varphi_{21}\varepsilon+\varphi_{22})^{-1}.
\end{equation}
Then
\begin{equation}\label{eq:8.8}
M_\pi(G,\Omega_+)= \kappa_1+\kappa_2.
\end{equation}
\end{lem}
\begin{proof}
Let $\varepsilon\in \cS^{\ptq}_{\kappa_2}$ admit the
Kre\u{\i}n-Langer factorizations \eqref{KLepsilon}. Then the
factorization (\ref{G}) of $G$ can be rewritten as
\begin{equation}\label{Gr}
\ G=T_\Phi[\varepsilon]=(\varphi_{11}\varepsilon_r+\varphi_{12}\theta_r)
     (\varphi_{21}\varepsilon_r+\varphi_{22}\theta_r)^{-1}.
\end{equation}
Since $\Phi$ is outer,
\[
\mbox{ker}\begin{bmatrix}
 \varphi_{11}(\lambda)\varepsilon_r(\lambda)+\varphi_{12}(\lambda)
\theta_r(\lambda) \\
  \varphi_{21}(\lambda)\varepsilon_r(\lambda)+\varphi_{22}(\lambda)
\theta_r(\lambda)
\end{bmatrix}
=\mbox{ker }\Phi(\lambda)\begin{bmatrix}
  \varepsilon_r(\lambda)\\
  \theta_r(\lambda)
\end{bmatrix}=\{0\}
\]
for every  $\lam\in\Omega_+$. By Lemma~\ref{Corona1} (see also~\cite[Lemma
3.3]{DerDym08}) the factorization (\ref{Gr}) of $G$ is right coprime over
$\Omega_+$.
Consequently, Proposition~\ref{Prop:5.1} and (\ref{eq:8.5}) imply that
\begin{equation}\label{eq:nov08}
    M_\pi(G,\Omega_+)=
M_\pi((\varphi_{21}\varepsilon_r+\varphi_{22}\theta_r)^{-1},\Omega_+)
=\kappa_1+\kappa_2.
\end{equation}
\end{proof}

We will also need the following general  noncancellation lemma
from~\cite{DerDym08}.
\begin{lem}\label{lem:5.4}
Let $G\in H_{\kappa,\infty}^{\ptq}$, $H_1\in H_{\infty}^{p\times p}$
and $H_2\in H_{\infty}^{q\times q}$. Then
\begin{equation}\label{b1G}
    M_\pi(H_1G,\Omega_+)=M_\pi(G,\Omega_+)\Longrightarrow
M_\pi(H_1GH_2,\Omega_+)=M_\pi(GH_2,\Omega_+),
\end{equation}
whereas
\begin{equation}\label{Gb2}
    M_\pi(GH_2,\Omega_+)=M_\pi(G,\Omega_+)\Longrightarrow
M_\pi(H_1GH_2,\Omega_+)=M_\pi(H_1G,\Omega_+).
\end{equation}
\end{lem}

\subsection{Proof of Theorem~\ref{thm:0.2}.}
\label{subsec:prfthm:0.2}
Since the first assertion (i.e., the verification of (\ref{epsilon}))
is covered by Theorem \ref{thm:3.13}, it remains only to
justify the second assertion. Towards this end, it is convenient to first
note that
if $\varepsilon\in\cS^{\ptq}$, then, by Lemma \ref{lem:5.3},
\begin{equation}
\label{eq:sep12a8}
M_\pi(\varphi_{21}\varepsilon +\varphi_{22})^{-1},\Omega_+)
=M_\pi((\wt\varphi_{11}^\#+\varepsilon\wt\varphi_{12}^\#)^{-1},\Omega_+)=
\nu_-(P)=\kappa.
\end{equation}
Thus, in view of the factorization $W=\Theta\Phi$ supplied in Theorem
\ref{factofW} and the properties of these factors, it follows that the mvf
$s=T_W[\varepsilon]$ is equal to
\begin{equation}\label{eq:10.3}
s=T_W[\varepsilon]=T_\Theta[G]=b_1Gb_2+K,
\end{equation}
where $G$ is defined by~\eqref{G}.
Therefore, since $K$ is holomorphic in $\Omega_+$,
\begin{equation}\label{MbGb}
M_\pi(b_1Gb_2,\Omega_+)=M_\pi(s,\Omega_+),
\end{equation}
 and  hence, if $s\in\cS_\kappa^{\ptq}$, then,  by (\ref{eq:sep12a8}),
\begin{equation}\label{eq:10.5}
\kappa=M_\pi(b_1Gb_2,\Omega_+)\le
M_\pi((\varphi_{21}\varepsilon+\varphi_{22})^{-1}b_2,\Omega_+) \le
M_\pi((\varphi_{21}\varepsilon+\varphi_{22})^{-1},\Omega_+)=\kappa.
\end{equation}
Thus, in view of  Proposition~\ref{Prop:5.1},
the factorization
$b_2^{-1}(\varphi_{21}\varepsilon+\varphi_{22})$ is coprime over
$\Omega_+$.

Similarly, since
$$
j_{pq}=W^\#j_{pq}W=\wtilde{\Phi}^\#\wtilde{\Theta}^\#j_{pq}\Theta\Phi=
\wtilde{\Phi}^\#j_{pq}\Phi,
$$
the mvf $G$ can be written as
\begin{equation}\label{eq:10.7}
    G=(\wt\varphi_{11}^\#+\varepsilon\wt\varphi_{12}^\#)^{-1}
    (\wt\varphi_{21}^\#+\varepsilon\wt\varphi_{22}^\#),
\end{equation}
and consequently the assumption $s\in\cS_\kappa^{\ptq}$
and~\eqref{MbGb} imply that
\[
M_\pi(b_1(\wt\varphi_{11}^\#+\varepsilon\wt\varphi_{12}^\#)^{-1},\Omega_+)
=M_\pi((\wt\varphi_{11}^\#+\varepsilon\wt\varphi_{12}^\#)^{-1},\Omega_+)=\kappa.
\]
Therefore,  Proposition~\ref{Prop:5.1} implies that the
factorization
$(\wt\varphi_{11}^\#+\varepsilon\wt\varphi_{12}^\#)b_1^{-1}$ is
coprime over $\Omega_+$. This completes the proof of the implication
$$
s\in\cS_\kappa^{\ptq}\Longrightarrow \textrm{the two factorizations in (a)
and (b) are coprime}.
$$

Suppose now that the two conditions (a) and (b) are met. Then by
assumption (b) and Proposition~\eqref{Prop:5.1}
\[
M_\pi((\varphi_{21}\varepsilon+\varphi_{22})^{-1}b_2,\Omega_+)=
M_\pi((\varphi_{21}\varepsilon+\varphi_{22})^{-1},\Omega_+)=\kappa.
\]
Since by  Lemma~\ref{lem:5.3a} $M_\pi(G,\Omega_+)=\kappa$ it follows from
Lemma~\ref{lem:5.4} that
\[
M_\pi(Gb_2,\Omega_+)=
M_\pi((\varphi_{21}\varepsilon+\varphi_{22})^{-1}b_2,\Omega_+)=\kappa.
\]
Similarly, (\ref{eq:sep12a8}), Lemma~\ref{lem:5.3a}, Lemma~\ref{lem:5.4}
and assumption (a) imply that
\[
M_\pi(b_1G,\Omega_+)=
M_\pi(b_1(\wt\varphi_{11}^\#+\varepsilon\wt\varphi_{12}^\#)^{-1},\Omega_+)
=\kappa.
\]
Therefore,
\[
M_\pi(b_1G,\Omega_+)=M_\pi(G,\Omega_+)=\kappa,
\]
which, with the help of Lemma~\ref{lem:5.4}, implies that
\[
M_\pi(b_1Gb_2,\Omega_+)=M_\pi(Gb_2,\Omega_+)=\kappa.
\]
Consequently,  $s\in \cS_{\kappa}^{p\times q}$, by~\eqref{eq:10.3}.
$\square$
\medskip

If $P$ is invertible, then  $\nu=0$ and Theorem  \ref{thm:0.2} takes the form
\begin{corollary}\label{cor:C13}
Let the data set $(A_1,A_2,C,P)$ satisfy the assumptions (B1)-(B3),
let $P$ be invertible, let $\kappa=\nu_-(P)$, and let the mvf's $W$,
$b_1$, $b_2$ be defined by~\eqref{eq:2.7}, \eqref{eq:7.11},
\eqref{eq:7.12}. Then:
\begin{enumerate}
\item[\rm(I)] $s\in
\cS_\kappa(A_1,A_2,C,P)$ if and only if $s$ belongs to $\cS_\kappa^{\ptq}$
and is of the form $s=T_W[\varepsilon]$ for some
$\varepsilon\in \cS^{p\times q}$.
\vspace{2mm}
    \item [\rm(II)] If $\varepsilon\in\cS^{\ptq}$, then $T_W[\varepsilon]\in
\cS_\kappa^{\ptq}$ if and only if
\begin{enumerate}
\item[\rm(a)]  the factorization
    $w^\#_{11}+\varepsilon w^\#_{12}=(\wt\varphi_{11}^\#
+\varepsilon\wt\varphi_{12}^\#)b_1^{-1}$ is
    coprime over $\Omega_+$ and
    \item [\rm(b)] the factorization
    $w_{21}\varepsilon+ w_{22}
=b_2^{-1}(\varphi_{21}\varepsilon+\varphi_{22})$ is coprime over $\Omega_+$.
\end{enumerate}
\end{enumerate}
\end{corollary}
\begin{remark}\label{rem:C13}
    The statement of Corollary~\ref{cor:C13} is a partial case of a general statement in~\cite[Theorem~1.3]{DerDym08}
    for arbitrary mvf $W\in\cU_{\kappa_1}^\circ(j_{pq})$. However the factorizations in (i), (ii)
    of Corollary~\ref{cor:C13} are simplified with respect to those in~\cite[Theorem~1.3]{DerDym08}
    since the outer mvf's $\varphi_1\in \cS_{out}^{\ptp}$ and $\varphi_2\in \cS_{out}^{\qtq}$
    in the problem (C1)--(C3) satisfy the conditions $\varphi_1^{-1}\in H_\infty^{\ptp}$,
    $\varphi_2^{-1}\in H_\infty^{\qtq}$.
\end{remark}

\subsection{Proof of Theorem~\ref{thm:0.1}.}
\label{subsec:prfthm:0.1}
Since assertion I is covered by
Theorem \ref{thm:3.2}, it remains only to justify II.

{\it Necessity.} Let $s\in \cS_\kappa(A_1,A_2,C,P)$. Then  by
Theorem \ref{thm:3.2} $s=T_W[\varepsilon]\in \cS^{\ptq}_\kappa$
 where  $\varepsilon\in
\cS_{\kappa-\kappa_1}^{\ptq}$. It follows from~\eqref{eq:10.3} that
\begin{equation}\label{G1}
   s-K=b_1(\varphi_{11}\varepsilon_r+\varphi_{12}\theta_r)
     (\varphi_{21}\varepsilon_r+\varphi_{22}\theta_r)^{-1}b_2\in \cS_\kappa^{p\times q
   }.
\end{equation}
Due to Lemma~\ref{lem:5.3} one has
\begin{equation}\label{MbGb2}
M_\pi((\varphi_{21}\varepsilon_r+\varphi_{22}\theta_r)^{-1},\Omega_+)=\kappa,
\end{equation}
 and, hence,
\begin{equation}\label{eq:10.52}
M_\pi((\varphi_{21}\varepsilon_r+\varphi_{22}\theta_r)^{-1}b_2,\Omega_+)=\kappa.
\end{equation}
It follows from~\eqref{eq:10.52}, \eqref{MbGb2} and Proposition~\ref{Prop:5.1} that the
factorization~\eqref{Reg2} is coprime over $\Omega_+$.

Similarly, it follows from \eqref{eq:10.3} and \eqref{eq:10.7} that
\begin{equation}\label{G2}
    s-K=b_1(\theta_\ell \wt\varphi_{11}^\#+\varepsilon_\ell\wt\varphi_{12}^\#)^{-1}
    (\theta_\ell \wt\varphi_{21}^\#+\varepsilon_\ell\wt\varphi_{22}^\#)b_2.
\end{equation}
Since
\begin{equation}\label{eq:10.72}
M_\pi((\theta_\ell
\wt\varphi_{11}^\#+\varepsilon_\ell\wt\varphi_{12}^\#)^{-1},\Omega_+)=\kappa,
\end{equation}
and $s\in \cS_\kappa^{p\times q}$ one obtains
\[
M_\pi(b_1(\theta_\ell
\wt\varphi_{11}^\#+\varepsilon_\ell\wt\varphi_{12}^\#)^{-1},\Omega_+)
=\kappa.
\]
This implies that the factorization~\eqref{Reg1} is coprime over
$\Omega_+$.

{\it Sufficiency.} The proof of sufficiency is similar to that in
Theorem~\ref{thm:0.2} except that now one should replace the
factorizations~\eqref{G} and~\eqref{eq:10.7} by ~\eqref{G1},
~\eqref{G2}.
$\Box$

\subsection{Takagi-Nudelman problem}
Let the data set $(A_1,A_2,C,P)$  satisfy the assumptions (B1)-(B4)
and let $\kappa\in\dN\cup\{0\}$. A mvf $s\in {\cS}_\kappa^{p\times
q}$ is said to be a solution of the Takagi-Nudelman problem  if $s$
satisfies the conditions $(C1)-(C3)$ and:
\begin{enumerate}
\item[(C4)] $s$ is holomorphic on $\sigma(A_1)\cup\overline{\sigma(A_2)}$.
\end{enumerate}

The set of all mvf's $s\in \cS_\kappa^{p\times q}$ satisfying
$(C1)-(C4)$ will be denoted  ${\cT\cN}_\kappa(A_1, A_2,C,P)$.

The following result gives a description of Takagi-Nudelman
interpolation problem. It contains, as a partial case, the Nudelman
description of the tangential indefinite interpolation problem
considered in \cite{Nud81}.
\begin{thm}\label{thm:8.5}
Let the data set $(A_1,A_2,C,P)$ satisfy the assumptions (B1)-(B4),
 let $\kappa=\nu_-(P)$, let the mvf $W(\lam)$ be defined by
(\ref{eq:0.5}) and let $\{b_1,b_2\}$ be an associated pair of the
mvf $W(\lam)$. Then $s\in {\cT\cN}_\kappa(A_1,A_2,C,P)$ if and only if
$s=T_W[\varepsilon]$, where $\varepsilon\in \cS^{p\times q}$ has the
form~\eqref{epsilon} and $\varepsilon$ satisfies one of the
equivalent conditions:
\begin{enumerate}
    \item [(i)] the mvf
    $(w^\#_{11}+\varepsilon w^\#_{12})b_1$ has no zeros on $\sigma(A_1)\cup\overline{\sigma(A_2)}$;
    \vspace{2mm}
    \item [(ii)] the mvf $b_2(w_{21}\varepsilon+w_{22})$ has no zeros
on $\sigma(A_1)\cup\overline{\sigma(A_2)}$.
\end{enumerate}
\end{thm}
\begin{proof}
{\it Necessity.} Let $s\in {\cT\cN}_\kappa(A_1,A_2,C,P)$. Then by Theorem
\ref{thm:3.13} $s=T_W[\varepsilon]\in \cS^{\ptq}_\kappa$, where the
mvf $\varepsilon\in \cS^{\ptq}$ has the form~\eqref{epsilon}. Since
\[
M_\pi(s,\Omega_+)=M_\pi((\varphi_{21}\varepsilon+\varphi_{22})^{-1},\Omega_+)
=\kappa,
\]
it follows from~\eqref{G1} that poles of $(\varphi_{21}\varepsilon
+\varphi_{22})^{-1}$
and hence poles of $(w_{21}\varepsilon+w_{22})^{-1}b_2^{-1}$ coincide
with poles of $s$.
Then (ii) is implied by (C4). Similarly the equalities~\eqref{G2} and
\[
M_\pi(s,\Omega_+)=M_\pi((\wt\varphi^\#_{11}+\varepsilon
\wt\varphi^\#_{12})^{-1},\Omega_+)=\kappa,
\]
imply that poles of
$(\wt\varphi^\#_{11}+\varepsilon \wt\varphi^\#_{12})^{-1}$ coincide
with poles of $s$, which serves to prove (i).

{\it Sufficiency.} Lemma \ref{lem:5.3a} guarantees that
$$
M_\pi(G,\Omega_+)=\kappa_1+\kappa_2=\kappa
$$
in the present setting for the mvf $G(\lambda)$ defined by~\eqref{G}.
Moreover, in view of (\ref{eq:7.41a}) and assumption (ii),
$G$ has no poles in $\sigma(A_1)\cup\overline{\sigma(A_2)}$. Therefore,
\[
M_\pi(G,\Omega_+\setminus(\sigma(A_1)\cup\overline{\sigma(A_2)})) =M_\pi(G,\Omega_+)
=\kappa.
\]
Thus, as $b_1(\lambda)$ and $b_2(\lambda)$ are invertible in
$\sigma(A_1)\cup\overline{\sigma(A_2)}$,
\[
M_\pi(b_1Gb_2,\Omega_+\setminus(\sigma(A_1)\cup\overline{\sigma(A_2)}))
=M_\pi(b_1Gb_2,\Omega_+) =\kappa
\]
and hence $s=K+b_1Gb_2\in \cS_{\kappa}^{p\times q}$.

Similarly if (i) is in force, then~\eqref{eq:nov08}, the second
equality in~\eqref{eq:sep12a8} and the dual
representation~\eqref{eq:2.7} of $s$ lead to the desired
conclusions.
\end{proof}

\begin{remark}
\label{rem:jun21a9}
Descriptions of the set ${\cT\cN}_\kappa(A_1,A_2,C,P)$ are available for
$\kappa\ge\nu_-(P)$ when $P$ is invertible;
see~\cite[Theorem 19.2.1]{BGR} for a description of rational solutions,
and, for the general case, \cite[Theorem~2]{Nud81} for the one sided problem
and~\cite[Proposition~4.7]{AmDer} for the two-sided problem
with $\sigma(A_1)\cap\overline{\sigma(A_2)}=\emptyset$.
\end{remark}

\begin{corollary}\label{cor:6.8}
Let (B1)-(B3) be in force, let $P$ be invertible,
$\kappa_1:=\nu_-(P)\le\kappa$, let the mvf $W(\lam)$ be defined by
(\ref{eq:1.84}), let $\{b_1,b_2\}$ be an associated pair of the mvf
$W(\lam)$, and let a mvf $\varepsilon\in
\cS^{\ptq}_{\kappa-\kappa_1}$ admits the Kre\u{\i}n-Langer
factorizations~\eqref{KLepsilon}.
Then $T_W[\varepsilon]\in {\cT\cN}_\kappa(A_1,A_2,C,P)$ if and only if one
of the equivalent conditions holds:
\begin{enumerate}
    \item [(i)] the mvf
    $(\theta_\ell w^\#_{11}+\varepsilon_\ell w^\#_{12})b_1$ has no zeros on $\sigma(A_1)\cup\overline{\sigma(A_2)}$;
    \item [(ii)] the mvf $b_2(w_{21}\varepsilon_r+w_{22}\theta_r)$ has no zeros on $\sigma(A_1)\cup\overline{\sigma(A_2)}$.
\end{enumerate}
\end{corollary}

 The next example shows that the inclusion
\[
{\cT\cN}_\kappa(A_1,A_2,C,P)\subseteq \cS_\kappa(A_1,A_2,C,P)
\]
may be strict.

\begin{example}\label{ex:6.8} Let $n_1=n=1$,
$M=A_1=0$, $N=1$, $C=\mbox{col}(2,0,0,1)$ and $\kappa=1$. Then the
unique solution of the Lyapunov-Stein equation (\ref{eq:0.LS}) is
$P=-3$, and, hence, $\kappa_1=1$. By formula (\ref{eq:1.84}) with
$\mu=1$, the resolvent matrix $W(\lam)$ can be written as
$$
W(\lam)=I_4+\frac{1-\lam}{3\lam} C\ C^* j_{22} $$ Direct
calculations show that
$$
W(\lam) =\left[
\begin{array}{ccccc}
\frac{4-\lam}{3\lam}& 0 &|&0&\frac{2(\lam-1)}{3\lam}\\
                 &   &|& & \\
               0 & 1 &|&0&0\\
                 &   &|& & \\
\hline\\
                0& 0 &|&1&0\\
                 &   &|& & \\
\frac{2(1-\lam)}{3\lam}&0&|&0&\frac{4\lam-1}{3\lam}
\end{array}\right],\quad
W^\#(\lam) =\left[
\begin{array}{ccccc}
\frac{4\lam-1}{3}& 0 &|&0&\frac{2(\lam-1)}{3}\\
                 &   &|& & \\
               0 & 1 &|&0&0\\
                 &   &|& & \\
\hline\\
                0& 0 &|&1&0\\
                 &   &|& & \\
\frac{2(1-\lam)}{3}&0&|&0&\frac{4-\lam}{3}
\end{array}\right]
$$
and hence, as follows from (\ref{eq:7.41}),
$$
b_1(\lam)=I_2\quad\textrm{and}\quad b_2(\lam)=\begin{bmatrix}1&0\\0&\lambda
\end{bmatrix}.
$$
The mvf
\[
\varepsilon(\lam)=\left[\begin{array}{cc}
  \frac{3}{4-\lam} & \frac{2(1-\lam)}{4-\lam} \\ \\
 \frac{2(1-\lam)}{4-\lam} & \frac{3\lam}{4-\lam}
\end{array}\right]\in \cS_0^{2\times 2},
\]
since $\varepsilon(\lam)$ is holomorphic in $\DD$ and
$$
I_2-\varepsilon(\lam)^*\varepsilon(\lam)=\frac{1-\vert\lam\vert^2}
{\vert 4-\lam\vert^2}\begin{bmatrix}3&6\\6&12\end{bmatrix}\ge
0\quad\textrm{for}\ \lam\in\DD.
$$
However, $\varepsilon(\lam)$
does not satisfy  condition (i) of Theorem~\ref{thm:8.5} since
the mvf
$$
(w^\#_{11}+\varepsilon w^\#_{12})^{-1}= \left[\begin{array}{cc}
  \frac{4-\lam}{3\lam} &0 \\ \\
 \frac{2(\lam-1)}{3} & 1
\end{array}\right]
$$
has a pole at 0. The corresponding linear fractional transform
$$
s(\lam)=T_W[\varepsilon]=\left[
\begin{array}{cc}
  \lam^{-1} & 0\\
  0 &\lam
\end{array}\right]=b^{-1}_\ell  s_\ell \quad \left(b_\ell =
\left[
\begin{array}{cc}
  \lam & 0\\
  0 &1
\end{array}\right],\,s_\ell =\left[
\begin{array}{cc}
  1 & 0\\
  0 &\lam
\end{array}\right]\right)
$$
also has a pole at $0$, and is not a solution of the Takagi-Nudelman
problem (C1)-(C4). However, $s\in {\cS}_1(A_1,C,P)$ by I of Theorem
\ref{thm:0.1}. It is reassuring to check that
$$
[b_\ell \quad -s_\ell ]C(A_1-\lam)^{-1}=\frac{b_\ell  \xi-s_\ell  \eta}{-\lam}=\
\begin{bmatrix}-2\\1\end{bmatrix}\in H_2^p\quad\textrm{with}\ p=2.
$$
\end{example}

\begin{thm}\label{thm:6.8a}
Let (B1)-(B3) be in force, let $P$ be invertible, and let $\kappa_1=\nu_-(P)\le\kappa$.
Then the problem (C1)-(C4) is solvable.
\end{thm}
\begin{proof}
Suppose first that $q\le p$ for the sake of definiteness and let
$\alpha_1,\dots,\alpha_l$ be all the points in
$\sigma(A_1)\cup\overline{\sigma(A_2)}$. In view of
Corollary~\ref{cor:6.8} one has to show that there is a mvf
$\varepsilon\in \cS^{\ptq}_{\kappa-\kappa_1}$ with the
Kre\u{\i}n-Langer factorizations~\eqref{KLepsilon} such that
\begin{equation}\label{Solv}
\det(\varphi_{21}(\alpha_j)\varepsilon_r(\alpha_j)+\varphi_{22}(\alpha_j)
\theta_r(\alpha_j))\ne 0 \quad\text{for}\ j=1,\dots,l.
\end{equation}
Let us choose the Blaschke-Potapov factor $\theta_r$ of degree
$\kappa-\kappa_1$ in such a way that $\theta_r(\alpha_j)$ is
invertible. If $\varphi_{21}(\alpha_j)=0$, then
$\varphi_{22}(\alpha_j)$ is invertible and, hence, \eqref{Solv} is
satisfied for every $\varepsilon_r\in \cS^{\ptq}$.

Assume now that $r:=\mbox{rank }\varphi_{21}(\alpha_j)>0$ and let us show that the
algebraic manifold
\begin{equation}\label{Man1}
    \cM_j=\{Y\in \dC^{\ptq}:
    \det(\varphi_{21}(\alpha_j)Y+\varphi_{22}(\alpha_j)\theta_r(\alpha_j))=0\}
\end{equation}
does not coincide with $\dC^{\ptq}$. Choosing invertible
matrices $V_1\in\dC^{q\times
q}$ and  $V_2\in\dC^{p\times p}$ such that
\[
V_1\varphi_{21}(\alpha_j)V_2=\begin{bmatrix} I_r & 0\\ 0
&0\end{bmatrix},
\]
and setting
\[
\wt Y=V_2^{-1}Y,\quad Z=\begin{bmatrix} Z_1 \\
Z_2\end{bmatrix}=V_1\varphi_{22}(\alpha_j)\theta_r(\alpha_j),\quad
Z_{1}\in\dC^{r\times q},\quad Z_{2}\in\dC^{(q-r)\times q},
\]
one can rewrite the equality in~\eqref{Man1} in the form
\begin{equation}\label{Man2}
    \det\left(\begin{bmatrix} I_r & 0\\ 0 &0\end{bmatrix}\wt Y
+\begin{bmatrix} Z_1 \\
Z_2\end{bmatrix}\right)=0.
\end{equation}
Since $\mbox{rank}\begin{bmatrix} \varphi_{21}(\alpha_j) &
\varphi_{22}(\alpha_j)\end{bmatrix}=q$ one obtains
\[
\mbox{rank}\,V_1\begin{bmatrix} \varphi_{21}(\alpha_j) &
\varphi_{22}(\alpha_j)\end{bmatrix}\begin{bmatrix} V_2 & 0\\ 0
&\theta_r(\alpha_j)\end{bmatrix}
=\mbox{rank}\begin{bmatrix} I_r & 0 & Z_1 \\
0 &0 & Z_2\end{bmatrix}=q
\]
and hence $\mbox{rank}\, Z_2=q-r$. Therefore, there is an invertible
matrix $V_3\in\dC^{q\times q}$ such that
$
Z_2V_3=\begin{bmatrix} 0 & I_{q-r}\end{bmatrix}.
$
Setting
\[
\begin{bmatrix} \wh Y_{11} & \wh Y_{12}\\ \wh Y_{21} & \wh Y_{22}\end{bmatrix}
=\wt Y V_3,\quad
\wh Y_{11}\in\dC^{r\times r}, \quad
\begin{bmatrix} Z_{11} & Z_{12}\end{bmatrix}=Z_1V_3,\quad Z_{11}\in\dC^{r\times r},
\quad Z_{12}\in\dC^{r\times (q-r)},
\]
one can rewrite the equality~\eqref{Man2} in the form
\begin{equation}\label{Man3}
    \det\left(\begin{bmatrix} I_r & 0\\ 0 &0\end{bmatrix}
    \begin{bmatrix} \wh Y_{11} & \wh Y_{12}\\ \wh Y_{21} & \wh Y_{22}\end{bmatrix}
    +\begin{bmatrix} Z_{11} & Z_{12} \\
0 & I_{q-r}\end{bmatrix}\right)=0
\end{equation}
or, equivalently,
\begin{equation}\label{Man4}
    \det( \wh Y_{11} + Z_{11} )=0.
\end{equation}
It follows from~\eqref{Man4} that the manifold $\cM_j$ has the complex dimension at most
$pq-1$ in $\dC^{\ptq}$. Therefore, the manifold $\cup_{j=1}^l\cM_j$ is nowhere dense in
$\dC^{\ptq}$ and hence there is a constant contractive matrix $\varepsilon_r\in
\dC^{\ptq}$ which satisfies~\eqref{Solv}.

In the case $p<q$ one has to use the first condition (i) in
Corollary~\ref{cor:6.8}.
\end{proof}

The conclusions of Theorem~\ref{thm:6.8a} may fail to hold if either
(B3) is not in force or $P$ is not invertible.
\bigskip

\begin{example} \label{ex:8} Let $n_2=0$, $M=A_1=O_{2\times 2}$,  $N=I_2$,
$C=I_2$, $P=-j_{11}$ and $\kappa=1$. Then, clearly, assumptions (B1)
and (B2) are in force (see Remark \ref{rem:oct29a8}), $P$ is
invertible  and $\nu_-(P)=1=\kappa$, while (B3) fails to hold, since
$C_2=[0\quad 1]$ and
\[ \bigcap_{j=0}^1\ker
C_2A_1^j=\ker C_2=\textup{span
}\left\{\begin{bmatrix}1\\0\end{bmatrix}\right\}.\] The
interpolation condition (C1) can be rewritten as
$$
b_\ell (0)\xi_j=s_\ell (0)\eta_j\quad (j=1,2).
$$
This implies that
$$
b_\ell (0)=0,\quad s_\ell (0)=0,
$$
which contradicts the noncancellation condition \eqref{KLcanon}.
Therefore the problem (C1)--(C3) 
and hence the problem (C1)--(C4) has no solution in $\cS_1$. We
remark that formula~\eqref{eq:3.25} is still valid: by formula
(\ref{eq:1.84}) with $\mu=1$,
$$
W(\lam)={\lam}^{-1} I_2
$$ and for every parameter $\varepsilon\in S_0$ one has
$ T_W[\varepsilon]=\varepsilon\not\in \cS_1. $ Therefore,
$T_W[S_0]\cap \cS_1=\emptyset$ and $\cS_1(A_1,C_1,P)=\emptyset$.
\bigskip
\end{example}

If $P$ is not invertible, then  the set $\cS_\kappa(A_1,A_2,C,P)$
may be empty. A criterion for the solvability of degenerate scalar
Nevanlinna-Pick problems  can be found in~\cite{Wor97}.

\subsection{Excluded parameters}
\begin{definition} Let the data $(A_1,A_2,C,P,\kappa)$
of the problem (C1)-(C4) satisfies the assumptions (B1)-(B4) and let
$\kappa_1:=\nu_-(P)(\leq \kappa)$. Then the parameter
$\varepsilon\in \cS^{\ptq}_{\kappa-\kappa_1}$ of the
form~\eqref{epsilon} is said to be excluded
for the problem (C1)-(C4), if
$T_W[\varepsilon]\not\in{\cT\cN}_\kappa(A_1,A_2,C,P)$.
\end{definition}

According to this definition, all the parameters $\varepsilon\in
S_0$ are excluded for the problem (C1)-(C4) in Example \ref{ex:8}.
Theorem~\ref{thm:8.5} leads to the following descriptions of
excluded parameters for the problem (C1)--(C4).
\begin{proposition}\label{ExcludedP}
Under the assumptions of Theorem~\ref{thm:8.5} the parameter
$\varepsilon\in \cS_{\kappa-\kappa_1}^{\ptq}$ of the
form~\eqref{epsilon} is excluded for the problems (C1)-(C4), if and
only if
\begin{equation}\label{Excl1}
\det(\varphi_{21}(\alpha_j)\varepsilon_r(\alpha_j)+\varphi_{22}(\alpha_j)
\theta_r(\alpha_j))=0.
\end{equation}
for at least one of the point
$\alpha_j\in\sigma(A_1)\cup\overline{\sigma(A_2)}$ $(j=1,\dots,l)$.
\end{proposition}
In accordance with the above statement the parameter $\varepsilon\in
\cS^{\ptq}_{\kappa-\kappa_1}$ will be said to be excluded at the
point $\alpha_j\in\sigma(A_1)\cup\overline{\sigma(A_2)}$ if
\eqref{Excl1} holds. In a special case when $\nu_-(P)=\kappa$ we can
give a criterion for the problem (C1)-(C4) to have no excluded
parameters.
\begin{proposition} Let (B1)-(B4) be in force, and let
$\nu_-(P)=\kappa$. Then the problem (C1)-(C4) has no excluded parameters at
$\alpha_j\in\sigma(A_1)\cup\overline{\sigma(A_2)}$ $(j=1,\dots,l)$ if and only if
\begin{equation}\label{Excl2}
\varphi_{21}(\alpha_j)\varphi_{21}(\alpha_j)^*-\varphi_{22}(\alpha_j)\varphi_{22}(\alpha_j)^*<0.
\end{equation}
\end{proposition}
\begin{proof}
A matrix $\varepsilon\in \cS^{\ptq}$ is an excluded parameter at
$\alpha_j$ if and only if there is a vector $v\in\dC^q$, $(v\ne 0)$
such that
\begin{equation}\label{Tan1}
    v^*\varphi_{21}(\alpha_j)\varepsilon(\alpha_j)+v^*\varphi_{22}(\alpha_j)=0.
\end{equation}
The tangential interpolation problem~\eqref{Tan1} has a solution
$\varepsilon\in \cS^{\ptq}$ if and only if the corresponding Pick
matrix
\begin{equation}\label{Pick1}
    (1-|\alpha_j|^2)^{-1}
    (v^*\varphi_{21}(\alpha_j)\varphi_{21}(\alpha_j)^*v-v^*\varphi_{22}(\alpha_j)\varphi_{22}(\alpha_j)^*v)
\end{equation}
is nonnegative. Therefore the problem~\eqref{Tan1} has no solutions for every
$v\in\dC^q$ if and only if the matrix in~\eqref{Excl2} is negative.
\end{proof}
\begin{corollary}
Let (B1) - (B4) be in force, let $\nu_-(P)=\kappa$ and let the
conditions~\eqref{Excl2} are satisfied for all
$\alpha_j\in\sigma(A_1)\cup\overline{\sigma(A_2)}$ $(j=1,\dots,l)$.
Then the problem (C1)-(C4) has a solution $s\in \cS_\kappa^{\ptq}$.
This solution is unique if and only if
\[
\rank(M^*P^2M+N^*P^2N+C^*C)-\rank P=\min\{p,q\}.
\]
\end{corollary}
The proof is immediate from Theorem~\ref{thm:8.5} and
Proposition~\ref{ExcludedP}.

\subsection{The Schur-Takagi interpolation problem} \label{TSproblem}
In this subsection $\Omega_+$ is either $\DD_+$ or $\Pi_+$. Let
$b_1\in \cS_{in}^{p\times p}$, $b_2\in \cS_{in}^{q\times q}$ be
inner mvf's,
let $K\in H_\infty^{p\times q}$ and let $\kappa\in\dN\cup\{0\}$. Consider the
following
problem.
\begin{definition}\label{def:GSIP}
    A ${p\times q}$ mvf $s$ is called a solution of the Schur-Takagi
interpolation problem in the class $\cS_{\kappa}^{p\times q}$, if
$\kappa\ge 1$, $s\in \cS_{\kappa}^{p\times q}$ and
\begin{equation}\label{eq:10.1}
    b_1^{-1}(s-K)b_2^{-1}\in H_{\kappa,\infty}^{p\times q}\setminus
H_{\kappa-1,\infty}^{p\times
    q}.
\end{equation}
The set of solutions of this problem is denoted
$\cS_\kappa(b_1,b_2,K)$.
\end{definition}

The analogue of this problem for $\kappa=0$ (in which the right hand side of
(\ref{eq:10.1}) is replaced by $H_\infty^{p\times q}$) has been extensively
studied by D. Z. Arov; see e.g., \cite{Arov93} and, for a more accessible
discussion and additional references, Chapter 7 of \cite{ArovD97}.

A ${p\times q}$ mvf $s$ is called a solution of the Takagi-Sarason
problem with the data set $(b_1,b_2,K)$ if  $s$  belongs to
$\cS_{\kappa'}^{p\times q}$ for some $\kappa'\le \kappa$ and
satisfies~\eqref{eq:10.1} (see~\cite{BGR}). The set of solutions of
the Takagi-Sarason problem is designated as ${\cT\cS}_\kappa(b_1,b_2,K)$.

\begin{thm}\label{thm:9.1}
Let $W\in \cU_{\kappa_1}^\circ(j_{pq})\cap \wt L_2^{\mtm}$,  and let
the mvf's $b_1$, $b_2$, $K$ be defined by~\eqref{eq:0.8},
\eqref{eq:0.7} and \eqref{eq:11.28}, respectively. Then:
\begin{enumerate}
\item[(i)] $s\in {\cT\cS}_\kappa(b_1,b_2,K)$ if and only if $s\in
T_W[\cS_{\kappa-\kappa_1}^{p\times q}]$;
\item[(ii)] $s\in \cS_\kappa(b_1,b_2,K)$ if and only if $s\in
T_W[\cS_{\kappa-\kappa_1}^{p\times q}]\cap \cS_\kappa^{\ptq}$.
\end{enumerate}
\end{thm}
\begin{proof}
(i) Let $s=T_W[\varepsilon]$, where $\varepsilon\in
\cS_{\kappa_2}^{p\times q}$ ($\kappa_2=\kappa-\kappa_1$) and let the
mvf $\varepsilon$ admit the Kre\u{\i}n-Langer
factorizations~\eqref{KLepsilon}. Then it follows from~\eqref{G1}
that
\begin{equation}\label{eq:9.2A}
    b_1^{-1}(s-K)b_2^{-1}=(\varphi_{11}\varepsilon_r+\varphi_{12}\theta_r)
     (\varphi_{21}\varepsilon_r+\varphi_{22}\theta_r)^{-1}
\end{equation}
and hence coincides with $G$ in~\eqref{Gr}. Then by Lemma~\ref{lem:5.3a}
\[
M_\pi(b_1^{-1}(s-K)b_2^{-1},\Omega_+)=\kappa_1+\kappa_2=\kappa.
\]
Since  $s\in \cS_{\kappa'}^{\ptq}$ for some $\kappa'\le\kappa$ by
Lemma~\ref{pr:2}, this proves that $s\in {\cT\cS}_\kappa(b_1,b_2,K)$.

Conversely, let $s\in {\cT\cS}_\kappa(b_1,b_2,K)$. Then $s\in
\cS_{\kappa'}^{\ptq}$ for some $\kappa'\le\kappa$. Consequently, the
mvf
\[
\varepsilon=T_{W^{-1}}[s].
\]
belongs to $\cS_{\kappa_2}^{\ptq}$, where
$\kappa'-\kappa_1\le\kappa_2$ by Lemma~\ref{pr:2} and, by an
analogous argument, $\kappa_2\le\kappa'+\nu_+(P)$.
In view of Lemma~\ref{lem:5.3a},
\[
M_\pi(b_1^{-1}(s-K)b_2^{-1},\Omega_+)=\kappa_1+\kappa_2.
\]
Now it follows from~\eqref{eq:10.1} that $\kappa_2=\kappa-\kappa_1$,
and hence that $s\in T_W[\cS_{\kappa-\kappa_1}^{\ptq}]$.

(ii) The second statement is immediate from (i) and
Definition~\ref{def:GSIP}.
\end{proof}

\begin{corollary}\label{cor:9.1}
Let the data set $(A_1,A_2,C,P)$ satisfy the assumptions (B1)-(B3),
let $P$ be invertible, let the mvf's $W$, $b_1$, $b_2$, $K$ be
defined by~\eqref{eq:0.5}, \eqref{eq:7.12}, \eqref{eq:7.11} and
\eqref{K}, respectively,  and assume that
$\kappa_1=\nu_-(P)\le\kappa$. Then $
\cS_\kappa(A_1,A_2,C,P)=\cS_\kappa(b_1,b_2,K)$.
\end{corollary}
\begin{proof}
The proof follows from the descriptions of the sets
$\cS_\kappa(A_1,A_2,C,P)$ and $\cS_\kappa(b_1,b_2,K)$ in
Theorem~\ref{thm:3.2} and Theorem~\ref{thm:9.1}, respectively.
\end{proof}
 It follows from Theorem~\ref{thm:9.1} and Corollary~\ref{cor:9.1}
 that
for every data set $(A_1,A_2,C,P)$ satisfying the assumptions
(B1)-(B3) with invertible $P$ there exists  a data set $(b_1,b_2,K)$
($b_1\in \cS_{in}^{p\times p}$, $b_2\in \cS_{in}^{q\times q}$,
$K\in H_\infty^{p\times q}$), such that
\[
\cS_\kappa(A_1,A_2,C,P)=\cS_\kappa(b_1,b_2,K)\subseteq
{\cT\cS}_\kappa(b_1,b_2,K).
\]

\section{Bitangential interpolation in the right half plane}
\label{sec:bipinrhp}
In this section we shall summarize the main
formulas that come into play when $\Omega_+$ (resp., $\Omega_-$) is the open
right (resp., left)  half plane $\Pi_+$ (resp., $\Pi_-$) and $\Omega_0=i\RR$.
In this setting, the bitangential interpolation problem under consideration
is formulated in terms of the matrices
\begin{equation}
\label{eq:dec10a8}
    M = A=\left[
   \begin{array}{cc}
   A_1&0\\
   0& A_2
   \end{array}\right],\quad
  N=I_n\quad \text{and},\quad J=j_{pq},\quad
C=\begin{bmatrix}
  C_{11} & C_{12} \\
  C_{21} & C_{22} \\
\end{bmatrix}:
\begin{bmatrix}
  \dC^{n_1} \\
  \dC^{n_2} \\
\end{bmatrix}\to
\begin{bmatrix}
  \dC^{p} \\
  \dC^{q} \\
\end{bmatrix},
\end{equation}
where $n_1+n_2=n$, $n_1>0$, $n_2>0$,
$A_1\in\dC^{n_1\times n_1}$ and  $A_2\in\dC^{n_2\times n_2}$. These matrices
are subject to the following constraints:
\begin{enumerate}
\item[(B1a)] $\sigma(A_1)\subset \Pi_+$ and $\sigma(A_2)\subset \Pi_-$.
\item[(B2a)] 
$P$ is  a Hermitian solution of the Lyapunov-Stein equation
\begin{equation}
\label{eq:dec10b8}
A^*P+PA+C^*JC=0
\end{equation}
(or, equivalently, in terms of $M$ and $N$, $M^*PN+N^*PM+C^*JC=0$).
\item[(B3a)]
The pairs $(C_{12}, A_2)$ and $(C_{21}, A_1)$ are observable.
\item [(B4a)]
$X\in\CC^{\ntn}$ is a Hermitian solution of the Riccati equation
\begin{equation}
\label{eq:dec10c8}
XA^*+AX+XC^*JCX=0
\end{equation}
such that:
\begin{enumerate}
    \item [(i)] $XPX=X$.
\vspace{2mm}
    \item [(ii)] $PXP=P$.
 \end{enumerate}
\end{enumerate}
Notice that because of the special form of (\ref{eq:dec10c8}),
$\textup{rng} X$ is automatically invariant under $M$ and $N$. If
$P$ is invertible, then (B4a) is superfluous, since it is
automatically satisfied by $X=P^{-1}$.

Let
\begin{equation}
\label{eq:dec10d8}
W(\lambda)=I_m+C(A-\lambda I_n)^{-1}XC^*J.
\end{equation}
Then
\begin{equation}
\label{eq:dec10e8}
W^\#(\lambda)=W(-\overline{\lambda})^*\quad\text{and}\quad
W^{-1}(\lambda)=JW^\#(\lambda)J=I_m+JCX(A^*+\lambda I_n)^{-1}C^*.
\end{equation}

\begin{example}\label{Mspace1}
If  assumptions (B1a)-(B4a) are in force, then the linear space of vvf's
$$
\cM=\{F(\lam)Xu:u\in \dC^n\}
$$
based on the mvf $F(\lambda)=C(A-\lambda I_n)^{-1}$  for
$\lam\not\in\sigma(A)$ and endowed with the inner product
\begin{equation}\label{eq:dec10f8}
   \langle FXu,\ FXv\rangle_{\cM}=v^*Xu
   \end{equation}
is an RKPS with RK
\begin{equation}\label{eq:1.82X}
   {\mathsf K}_\omega(\lam)=F(\lam)X F(\omega)^*
=\frac{J-W(\lambda)JW(\omega)^*}{\lambda+\overline{\omega}}
\end{equation}
and $\nu_-(X)$ negative squares.
\end{example}

\begin{lem}
\label{lem:dec10a8}
If  assumptions (B1a)-(B4a) are in force, and if
$F(\lambda)=C(A-\lambda I_n)^{-1}$ and $W(\lambda)$ is
given by (\ref{eq:dec10d8}), then
\begin{equation}
\label{eq:dec10g8}
W(\lambda)^{-1}F(\lambda)X=-CX(A^*+\lambda I_n)^{-1}
\end{equation}
and
\begin{equation}
\label{eq:dec10h8}
W(\lambda)^{-1}F(\lambda)v=C(I_n-XP)(A-\lambda I_n)^{-1}v\quad\text{for}
\quad v\in\ker P.
\end{equation}
\end{lem}

\begin{lem}
\label{lem:dec10b8}
If assumptions (B1a)-(B4a) are in force and if
$$
X_1=X\begin{bmatrix}I_{n_1}\\ 0\end{bmatrix}\quad\textrm{and}\quad
X_2=X\begin{bmatrix}0\\ I_{n_2}\end{bmatrix},
$$
then
$$
\ker\begin{bmatrix}CX_1\\ CX_{1}A_1^*\\ \vdots\\CX_{1}(A_1^*)^{n_1-1}
\end{bmatrix}=\ker X_1,
\quad
\ker\begin{bmatrix}CX_2\\ CX_{2}A_2^*\\ \vdots\\CX_{2}(A_2^*)^{n_2-1}
\end{bmatrix}=\ker X_2
$$
and
$$
\textup{rng}[X_j^*C^*\, A_jX_j^*C^*\, \cdots\,A_j^{n_j-1}X_j^*C^*]=
\textup{rng}X_j^*\quad\textrm{for}\ j=1,2.
$$
\end{lem}

\begin{lem}
\label{lem:dec10c8}
If $g\in H_2^{n_1}$ and $v\in\CC^{n_1}$, then (since $\sigma(A_1)
\subset \Pi_+$),
\begin{equation}
\label{eq:dec10i8}
(A_1-\lambda I_{n_1})^{-1}(g-v)\in H_2^{n_1}\Longleftrightarrow
v=\frac{1}{2\pi}\int_{-\infty}^\infty(A_1-it I_{n_1})^{-1}g(it)dt.
\end{equation}
If $g\in (H_2^{n_2})^\perp$ and $v\in\CC^{n_2}$, then (since $\sigma(A_2)
\subset \Pi_-$),
\begin{equation}
\label{eq:dec10j8}
(A_2-\lambda I_{n_2})^{-1}(g-v)\in (H_2^{n_2})^\perp\Longleftrightarrow
v=-\frac{1}{2\pi}\int_{-\infty}^\infty(A_2-it I_{n_1})^{-1}g(it)dt.
\end{equation}
\end{lem}

\begin{remark}
\label{rem:dec10a8}
Since $(C_{21}, A_1)$ is observable and $\sigma(A_1)\subset \Pi_+$, the matrix
$$
Q_1=-\int_0^\infty e^{-tA_1^*}C_{21}^*C_{21}e^{-tA_1}dt
=-\frac{1}{2\pi}\int_{-\infty}^\infty(A_1^*+i\mu I_{n_1})^{-1}C_{21}^*C_{21}
(A_1-i\mu I_{n_1})^{-1}d\mu
$$
is a negative definite solution of the Lyapunov equation
\begin{equation}
\label{eq:dec10k8}
A_1^*Q_1+Q_1A_1+C_{21}^*C_{21}=0
\end{equation}
and the mvf
\begin{equation}
\label{eq:dec10l8}
\widetilde{b}_2(\lambda)=I_q+C_{21}Q_1^{-1}(A_1^*+\lambda I_{n_1})^{-1}C_{21}^*
\end{equation}
is inner with respect to $\Pi_+$.
Similarly, since $(C_{12}, A_2)$ is observable and $\sigma(A_2)\subset \Pi_-$,
the matrix
$$
Q_2=\int_0^\infty e^{tA_2^*}C_{12}^*C_{12}e^{tA_2}dt
=\frac{1}{2\pi}\int_{-\infty}^\infty(A_2^*+i\mu I_{n_2})^{-1}C_{12}^*C_{12}
(A_2-i\mu I_{n_2})^{-1}d\mu
$$
is a positive definite solution of the Lyapunov equation
\begin{equation}
\label{eq:dec10m8}
A_2^*Q_2+Q_2A_2+C_{12}^*C_{12}=0
\end{equation}
and the mvf
\begin{equation}
\label{eq:dec10n8}
\widetilde{b}_1(\lambda)=I_p+C_{12}(A_2-\lambda I_{n_2})^{-1}Q_2^{-1}C_{12}^*
\end{equation}
is inner with respect to $\Pi_+$.

Let $F_{12}(\lam)=C_{12}(A_2-\lambda I_{n_2})^{-1}$,
$$
\widetilde{\cM}_2=\{F_{12}u: u\in \CC^{n_2}\}\quad\text{and}\quad
{\cM}_2=\{F_{12}[X_{21}\quad X_{22}]u: u\in \CC^n\}
$$
and let ${\cM}_2$ and $\widetilde{\cM}_2$ be equipped with the
(normalized standard) inner product
$$
\langle F_{12}u, F_{12}v\rangle_{nst}=v^*Q_2u
=v^*\left\{\frac{1}{2\pi}\int_{-\infty}^\infty F_{12}(i\mu)^*F_{12}(i\mu)
d\mu\right\}u.
$$
Then $\widetilde{\cM}_2$ is a RKHS with RK
$$
F_{12}(\lambda)Q_2^{-1}F_{12}(\omega)^*=\frac{I_p-\widetilde{b}_1(\lambda)
 \widetilde{b}_1(\omega)^*}{\lambda+\overline{\omega}},
$$
i.e.,
$$
\widetilde{\cM}_2=\cH(\widetilde{b}_1).
$$
Moreover,  $\cM_2$ is a subspace of $\widetilde{\cM}_2$ that is invariant
under the backwards shift operator $R_0$. Therefore, by the Beurling-Lax
theorem,
$$
\cM_2=\cH(b_1)\subseteq \cH(\widetilde{b}_1)=\widetilde{\cM}_2.
$$
Similar considerations based on
$F_{21}(\lam)=C_{21}(A_1-\lambda I_{n_1})^{-1}$,
$$
\widetilde{\cM}_1=\{F_{21}u: u\in \CC^{n_1}\}\quad\text{and}\quad
{\cM}_1=\{F_{21}[X_{11}\quad X_{12}]u: u\in \CC^n\}
$$
imply that
$$
\cM_1=\cH_*(b_2)\subseteq \cH_*(\widetilde{b}_2)=\widetilde{\cM}_1.
$$
\end{remark}

The subsequent analysis will make use of the left coprime factorizations
\begin{equation}
\label{eq:dec10r8}
C_{21}(A_1-\lambda I_{n_1})^{-1}[X_{11}\quad X_{12}]=b_2^{-1}\varphi_2
\end{equation}
and
\begin{equation}
\label{eq:dec10s8}
C_{12}(A_2-\lambda I_{n_2})^{-1}[X_{21}\quad X_{22}]=(b_1^\#)^{-1}\varphi_1
\end{equation}
where $b_2$ and $b_1$ are inner with respect to $\Pi_+$,
$\varphi_2\in{\cR}\cap H_2^{q\times n}$ and $\varphi_1\in
(H_2^\perp)^{p\times n}$. If $\textup{rank} X_j=n_j$ for $j=1,2$, then
$\widetilde{b}_j(\lambda)=b_j(\lambda)$.

\begin{lem}
\label{lem:dec10d8}
If  assumptions (B1a)-(B4a) are in force, then
(since $(C_{21},A_1)$ is observable and $\sigma(A_1)\subset \Pi_+$)
\begin{equation}
\label{eq:dec10p8}
\begin{split}
\{P_+(A_1^*+\lambda I_{n_1})^{-1}h: h\in (H_2^{n_1})^\perp\}&=
\{P_+(A_1^*+\lambda I_{n_1})^{-1}C_{21}^*g: g\in{\cH}_*(\widetilde{b}_2)\}\\
&=\{(A_1^*+\lambda I_{n_1})^{-1}x: x\in\CC^{n_1}\}.
\end{split}
\end{equation}
Similarly, (since $(C_{12}, A_2)$ is observable
and $\sigma(A_2)\subset \Pi_-$),
\begin{equation}
\label{eq:dec10q8}
\begin{split}
\{P_-(A_2^*+\lambda I_{n_2})^{-1}h: h\in H_2^{n_2}\}&=
\{P_-(A_2^*+\lambda I_{n_2})^{-1}C_{12}^*g: g\in\cH(\widetilde{b}_1)\}\\
&=\{(A_2^*+\lambda I_{n_2})^{-1}x: x\in\CC^{n_2}\}.
\end{split}
\end{equation}
\end{lem}

The rest of the development for $\Omega_+=\Pi_+$ is pretty much the
same as for $\Omega_+=\DD$, with the appropriate changes of $M$ and
$N$ and with (B1a)--(B4a) in place of (B1)--(B4) and leads to the
following conclusions:

\begin{thm}\label{thm:0.2a}
Let (B1a)--(B4a) be in force, let $\nu_-(P)=\kappa$ and let
\begin{equation}
\label{nua}
    \nu=\rank(P^2+C^*C)-\rank P.
\end{equation}
Then there are unitary matrices $U\in \dC^{p\times p}$, $V\in
\dC^{q\times q}$, such that $\cS_\kappa(M, N,C,P)$ if and only if
$s$ belongs to $\cS_\kappa^{\ptq}$ and is of the form
$s=T_W[\varepsilon]$, where
\begin{equation}\label{epsilon1}
    \varepsilon=U\left[\begin{array}{cc}
                                  \wt \varepsilon & 0 \\
                                        0 & I_\nu \\
                                      \end{array}  \right]V^*,\quad
\textrm{and}\quad \wt\varepsilon \in \cS^{(p-\nu)\times(q-\nu)}.
\end{equation}

If $\varepsilon\in\cS^{\ptq}$, then $T_W[\varepsilon]\in
\cS_\kappa^{\ptq}$ if and only if
\begin{enumerate}
    \item [(a)]
    the factorization
    $w^\#_{11}+\varepsilon w^\#_{12}
    =( \wt\varphi_{11}^\#+\varepsilon\wt\varphi_{12}^\#)b_1^{-1}$ is
    coprime over $\Omega_+$  and
\vspace{2mm}
    \item [(b)] the factorization
    $w_{21}\varepsilon+ w_{22}=b_2^{-1}(\varphi_{21}\varepsilon
+\varphi_{22})$ is
    coprime over $\Omega_+$.
\end{enumerate}
\end{thm}

\begin{thm}\label{thm:0.1a}
Let the data set $(M, N,C,P)$ satisfy the assumptions (B1a)-(B3a), let
$P$ be invertible,
$\kappa_1=\nu_-(P)\le\kappa$, and let the mvf's $W$, $b_1$, $b_2$ be defined
as in this section. Then:
\begin{enumerate}
\item[\rm(I)] $s\in \cS_\kappa(M, N,C,P)$
if and only if $s\in\cS_\kappa^{p\times q}$ and is of the form
$s=T_W[\varepsilon]$ for some $\varepsilon\in
\cS_{\kappa-\kappa_1}^{p\times q}$.
\item[\rm(II)] If  $\varepsilon\in
\cS_{\kappa-\kappa_1}^{p\times q}$ and $\theta_\ell$, $\theta_r$,
$\varepsilon_\ell$, $\varepsilon_r$ are a choice of its
Kre\u{\i}n-Langer factorizations
as in \eqref{eq:0.9},
  then $T_W[\varepsilon]\in
\cS_\kappa^{p\times q}$  if and only if the factorizations
\begin{equation}\label{Reg1a}
    \theta_\ell w^\#_{11}+\varepsilon_\ell w^\#_{12}
    =(\theta_\ell
    \wt\varphi_{11}^\#+\varepsilon_\ell\wt\varphi_{12}^\#)b_1^{-1},
\end{equation}
\begin{equation}\label{Reg2a}
    w_{21}\varepsilon_r+ w_{22}\theta_r=b_2^{-1}(\varphi_{21}\varepsilon_r
+\varphi_{22}\theta_r)
\end{equation}
are coprime over $\Omega_+$.
\end{enumerate}
\end{thm}

\end{document}